\documentclass[12pt]{article}
  
\usepackage{amsfonts}

\usepackage{hyperref}

\usepackage{xcolor}
 \definecolor{darkgreen}{rgb}{0,0.4,0}
 \hypersetup{
   colorlinks   = true, 
   urlcolor     = blue, 
   linkcolor    =blue, 
   citecolor   = darkgreen 
 }

\usepackage{latexsym}

\usepackage{amssymb}

\usepackage{graphicx}
  \usepackage{mathabx}

\usepackage{tikz-cd}

\newcommand{\D}{\displaystyle}

\newcommand{\breath}{\medskip}

\newtheorem{thm}{Theorem}[section]

\newcounter{claimcount}[thm]

\newcounter{subclaimcount}[claimcount]

\newtheorem{prop}[thm]{Proposition}

\newtheorem{lemma}[thm]{Lemma}

\newtheorem{cor}[thm]{Corollary}

\newcommand{\dfn}{\sf\em}

\newcommand{\Theorem}[2]{\begin{thm}{\sf #1}  #2 \end{thm}}

\newcommand{\Proposition}[2]{\begin{prop}{\sf #1}  #2 \end{prop}}

\newcommand{\Lemma}[2]{\begin{lemma}{\sf #1}  #2 \end{lemma}}

\newcommand{\Corollary}[2]{\begin{cor}{\sf #1}  #2 \end{cor}} 

\newcommand{\thmfont}[1]{{\sl #1}}    

\newcommand{\example}[1]{
        \refstepcounter{thm}
                     \begin{list}{}
 			{\setlength{\leftmargin}{0em}
 			\setlength{\rightmargin}{0em}}
        \item {\bf Example \thethm.} #1
                   \hfill$\diamondsuit$  \end{list}  
 			}

\newcommand{\bthmlist}{
 \begin{list}{{\bf (\alph{enumii})}}{\usecounter{enumii}}
 			{\setlength{\leftmargin}{0em}
 			\setlength{\itemsep}{0em}
 			\setlength{\parsep}{0em}
 			\setlength{\rightmargin}{0em}}}

\newcommand{\ethmlist}{\end{list}}

\renewcommand{\thesubclaimcount}{\arabic{claimcount}\Alph{subclaimcount}}

\newcommand{\Claim}[1]{\refstepcounter{claimcount}
 
               \noindent{\bf Claim \theclaimcount: \ }\thmfont{ #1}}

\newcommand{\Subclaim}[1]{\refstepcounter{subclaimcount}
 
                \noindent{\bf Claim \thesubclaimcount: \ }\thmfont{ #1}}

\newcommand{\claim}{\Claim}

\newcommand{\subclaim}{\Subclaim}

\newcommand{\bprf}[1][Proof.]{\begin{list}{}
 			{\setlength{\leftmargin}{1em}
 			\setlength{\rightmargin}{0em}
 			\setlength{\listparindent}{1em}}
                         \item  {\em \hspace{-1.5em}  #1   }}

\newcommand{\eprf}{\end{list}}

\newcommand{\bthmprf}{\setcounter{claimcount}{0}\bprf}

\newcommand{\bclaimprf}{\bprf}

\newcommand{\bsubclaimprf}{\bprf}

\newcommand{\ethmprf}{ \hfill$\Box$ 
 \eprf
 
 \breath
 
 }

\newcommand{\eclaimprf}{ \hfill $\Diamond$~{\scriptsize {\tt Claim~\theclaimcount}}\eprf}

\newcommand{\esubclaimprf}{ \hfill $\triangledown$~{\scriptsize 
 {\tt Claim~\thesubclaimcount}}\eprf}

\newcommand{\beq}{\begin{eqnarray*}}

\newcommand{\eeq}{\end{eqnarray*}}

\newcommand{\beqn}{ \begin{equation} }

\newcommand{\eeqn}{ \end{equation} }

\newcommand{\bitem}{\begin{itemize}}

\newcommand{\eitem}{\end{itemize}}

\newcommand{\bquote}{\begin{quotation}}

\newcommand{\equote}{\end{quotation}}

\newcommand{\If}{\mbox{\ if \ }} 

\newcommand{\And}{\mbox{\ and \ }}

\newcommand{\dI}{{\mathbb{I}}}

\newcommand{\dL}{{\mathbb{L}}}

\newcommand{\dN}{{\mathbb{N}}}

\newcommand{\dR}{{\mathbb{R}}}

\newcommand{\dZ}{{\mathbb{Z}}}

\newcommand{\barf}{{\overline{f}}}

\newcommand{\barmu}{{\overline{\mu }}}

\newcommand{\barnu}{{\overline{\nu }}}

\newcommand{\barphi}{{\overline{\phi }}}

\newcommand{\barsB}{{\overline{\mathcal{ B}}}}

\newcommand{\barsO}{{\overline{\mathcal{ O}}}}

\newcommand{\barsQ}{{\overline{\mathcal{ Q}}}}

\newcommand{\barsR}{{\overline{\mathcal{ R}}}}

\newcommand{\barsS}{{\overline{\mathcal{ S}}}}

\newcommand{\barsU}{{\overline{\mathcal{ U}}}}

\newcommand{\bI}{{\mathbf{ I}}}

\usepackage{amsbsy}

\newcommand{\sA}{{\mathcal{ A}}}

\newcommand{\sB}{{\mathcal{ B}}}

\newcommand{\sC}{{\mathcal{ C}}}

\newcommand{\sD}{{\mathcal{ D}}}

\newcommand{\sE}{{\mathcal{ E}}}

\newcommand{\sF}{{\mathcal{ F}}}

\newcommand{\sG}{{\mathcal{ G}}}

\newcommand{\sH}{{\mathcal{ H}}}

\newcommand{\sI}{{\mathcal{ I}}}

\newcommand{\sJ}{{\mathcal{ J}}}

\newcommand{\sK}{{\mathcal{ K}}}

\newcommand{\sL}{{\mathcal{ L}}}

\newcommand{\sM}{{\mathcal{ M}}}

\newcommand{\sN}{{\mathcal{ N}}}

\newcommand{\sO}{{\mathcal{ O}}}

\newcommand{\sP}{{\mathcal{ P}}}

\newcommand{\sQ}{{\mathcal{ Q}}}

\newcommand{\sR}{{\mathcal{ R}}}

\newcommand{\sS}{{\mathcal{ S}}}

\newcommand{\sT}{{\mathcal{ T}}}

\newcommand{\sU}{{\mathcal{ U}}}

\newcommand{\sV}{{\mathcal{ V}}}

\newcommand{\sX}{{\mathcal{ X}}}

\newcommand{\sY}{{\mathcal{ Y}}}

\newcommand{\sZ}{{\mathcal{ Z}}}

\newcommand{\gA}{{\mathfrak{ A}}}

\newcommand{\gB}{{\mathfrak{ B}}}
\newcommand{\gN}{{\mathfrak{ N}}}

\newcommand{\gC}{{\mathfrak{ C}}}

\newcommand{\gE}{{\mathfrak{ E}}}

\newcommand{\gF}{{\mathfrak{ F}}}

\newcommand{\gG}{{\mathfrak{ G}}}

\newcommand{\gM}{{\mathfrak{ M}}}

\newcommand{\gO}{{\mathfrak{ O}}}

\newcommand{\gP}{{\mathfrak{ P}}}

\newcommand{\gR}{{\mathfrak{ R}}}

\newcommand{\gU}{{\mathfrak{ U}}}

\newcommand{\ga}{{\mathfrak{ a}}}

\newcommand{\go}{{\mathfrak{ o}}}

\newcommand{\gr}{{\mathfrak{ r}}}

\newcommand{\alp }{\alpha}

\newcommand{\bet }{\beta}

\newcommand{\del }{\delta}

\newcommand{\eps }{\epsilon}

\newcommand{\kap }{\kappa}

\newcommand{\lam }{\lambda}

\newcommand{\sig }{\sigma}

\newcommand{\Ups }{\Upsilon}

\newcommand{\hb}{{\widehat{b}}}

\newcommand{\hf}{{\widehat{f}}}

\newcommand{\hg}{{\widehat{g}}}

\newcommand{\hs}{{\widehat{s}}}

\newcommand{\hdN}{{\widehat{\mathbb{ N}}}}

\newcommand{\hsB}{{\widehat{\mathcal{ B}}}}

\newcommand{\hsF}{{\widehat{\mathcal{ F}}}}

\newcommand{\hsG}{{\widehat{\mathcal{ G}}}}

\newcommand{\hsK}{{\widehat{\mathcal{ K}}}}

\newcommand{\hsO}{{\widehat{\mathcal{ O}}}}

\newcommand{\hsQ}{{\widehat{\mathcal{ Q}}}}

\newcommand{\hsR}{{\widehat{\mathcal{ R}}}}

\newcommand{\hsS}{{\widehat{\mathcal{ S}}}}

\newcommand{\hmu }{{\widehat{\mu}}}

\newcommand{\hphi}{{\widehat{\phi }}}

\newcommand{\fD}{{\mathsf{ D}}}

\newcommand{\fT}{{\mathsf{ T}}}

\newcommand{\lb}{\left}

\newcommand{\rb}{\right}

\newcommand{\maketall}{\rule[-0.5em]{0em}{1em}}

\newcommand{\implies}{\ensuremath{\Longrightarrow}}

\newcommand{\map}{{\longrightarrow}}

\newcommand{\goto}{{\rightarrow}}

\newcommand{\into}{{\map}}

\newcommand{\seilpmi}{{\Longleftarrow}}

\newcommand{\statement}[1]{\lb(  \maketall 
      \begin{minipage}{40em}
       \begin{tabbing}
         #1 
       \end{tabbing}
      \end{minipage}  \rb)}

\newcommand{\oo}{{\infty}}

\newcommand{\x}{\times}

\newcommand{\compl}[1]{#1^\complement}

\newcommand{\symdif}{{\bigtriangleup}}

\newcommand{\union}{\cup}

\newcommand{\Union}{\bigcup}

\newcommand{\intsct}{\cap}

\newcommand{\Intsct}{\bigcap}

\newcommand{\disj}{\sqcup}

\newcommand{\Disj}{\bigsqcup}

\newcommand{\set}[2]{{\left\{ #1 \; ; \; #2 \right\} }}

\newcommand{\norm}[2]{{\left\| #1 \right\|_{{#2}} }   }

\newcommand{\choice}[1]{{\lb\{ \begin{array}{rcl}
                                 #1 
                               \end{array}  \rb.  }}

\newcommand{\eeequals}[1]{\raisebox{-0.9ex}{$\overline{\overline{{\scriptscriptstyle{\mathrm{#1}}}}}$}}

\newcommand{\leeeq}[1]{\raisebox{-1ex}{${{\D\leq} \atop {\scriptscriptstyle{\mathrm{#1}}}}$}}

\newcommand{\lt}[1]{\raisebox{-1ex}{${{\D<} \atop {\scriptscriptstyle{\mathrm{#1}}}}$}}

\newcommand{\grt}[1]{\raisebox{-1ex}{${{\D>} \atop {\scriptscriptstyle{\mathrm{#1}}}}$}}

\newcommand{\geeeq}[1]{\raisebox{-1ex}{${{\D\geq} \atop {\scriptscriptstyle{\mathrm{#1}}}}$}}

\newcommand{\subseteeeq}[1]{\raisebox{-1.3ex}{$\stackrel{\D\subseteq}{\scriptscriptstyle{\mathrm{#1}}}$}}

\newcommand{\iiiff}[1]{\Leftarrow\!\!\!\!\lefteqn{\eeequals{#1}}\Rightarrow}

\newcommand{\goesto}[2]{{ -\!\!\!-\!\!\!-\!\!\!-\!\!\!\!\!\!\!\!\!\!\!
  ^{{\scriptscriptstyle #2}}_{{\scriptscriptstyle #1}}
   \!\!\!\!\!\!\!\!\!\longrightarrow }}

\newcommand{\Real}{\dR}

\newcommand{\Natur}{\dN}

\newcommand{\Zahl}{\dZ}

  \newcommand{\Co}{\sC_{0}}

 \newcommand{\Cb}{\sC_{\mathrm{b}}}

\newcommand{\drho}{\ \mathrm{d}\rho}

\newcommand{\dnu}{\ \mathrm{d}\nu}

\newcommand{\dmu}{\ \mathrm{d}\mu}

\newcommand{\dhmu}{\ \mathrm{d}\hmu}

\newcommand{\Int}{\mathrm{int}}

\newcommand{\Cl}{\mathrm{clos}}

\newcommand{\sssubset}[1]{\raisebox{-1ex}{${{\D\subseteq} \atop {\scriptscriptstyle{\mathrm{#1}}}}$}}

  \newcommand{\bone}{\boldsymbol{1}}

\newcommand{\restr}{\upharpoonleft}

\newcommand{\supp}{\mathrm{supp}}

\newcommand{\clop}{\mathfrak{Clp}}

\newcommand{\hgB}{{\widehat{\gB}}}

\newcommand{\CB}{\sC_\gB}

\newcommand{\GB}{\sG_\gB}

\newcommand{\Borel}{\gB{\scriptstyle\!\go\gr}}

\newcommand{\bargO}{\overline{\gO}}

 \newcommand{\hgG}{\widehat{\mathfrak{G}}}

\newcommand{\GhB}{\sG_{\hgB}}

\newcommand{\ChB}{\sC_{\hgB}}

\newcommand{\StPr}{\mathbf{StPr}}

\newcommand{\Cred}{\mathbf{Cont}}

\newcommand{\CmpCred}{\mathbf{CmpCnt}}

\newcommand{\GlPr}{\mathbf{GlPr}}

\newcommand{\bsC}{{\boldsymbol{\mathcal{C}}}}

\newcommand{\hbsC}{{\widehat{\boldsymbol{\mathcal{C}}}}}

\newcommand{\CpInt}{\mathbf{CmpInt}}

\newcommand{\LCPr}{\mathbf{LCPr}}

\newcommand{\INT}{\mathbf{Int}}

\newcommand{\bsI}{\boldsymbol{\mathcal{I}}}

 \usepackage{color}
\usepackage{xcolor}
\usepackage{marginnote}

\definecolor{darkblue}{rgb}{0,0,0.75}

\usepackage{esint}
 
 \newcommand{\integral}{\fint}
 \newcommand{\Integral}{\fint}

\begin{document}

\title{Measure and integration on Boolean algebras of regular open subsets in a topological space\thanks{This research was supported by Labex MME-DII (ANR11-LBX-0023-01).  We thank  David Fremlin,
Joel David Hamkins, Robert Furber, Tamar Lando, Vincenzo Marra, Dana Scott and Reem Yassawi for very helpful suggestions. We are especially grateful to KPS Bhaskara Rao for his many insightful comments.   They are not responsible for any  errors.}}

\author{Marcus Pivato\thanks{THEMA, Universit\'e de Cergy-Pontoise, France, {\tt marcuspivato@gmail.com}.} \ 
and Vassili Vergopoulos\thanks{Paris School of Economics, and University Paris 1 Panth\'eon-Sorbonne}.}
\maketitle

\begin{abstract}
 The regular open subsets of a topological space form a Boolean algebra, where the {\em join} of two regular open sets is the interior of the closure of their union.  A  {\em content} is a  finitely additive measure on this Boolean algebra, or  on one of its subalgebras.  We develop a theory of integration for such contents.  We then explain the  relationship between contents,  residual  charges, and Borel measures.  We show that a content can be represented by a normal Borel measure, augmented with a {\em liminal structure}, which specifies how two or more regular open sets share the measure of their common boundary.    In particular,  a content on a locally compact Hausdorff space  can be  represented by a normal Borel measure and a liminal structure on the Stone-\v{C}ech compactification of that space.   We also show how contents can be represented by Borel measures on the Stone space of the underlying Boolean algebra of regular open sets.  Finally, we show that these constructions are functorial.


\end{abstract}

\section{Introduction}

Let $\sS$ be a topological space, equipped with a Borel  measure $\mu$ having full support.  It is well-known that the apparent ``size'' of a subset $\sS$, as seen from a topological perspective, might greatly differ from its ``size'' from a measurable perspective.  For example, suppose $\sS$ is the unit interval $[0,1]$ equipped with the Lebesgue measure.  It is easy to construct a subset $\sO_n\subset\sS$ which is open and dense, but such that $\mu[\sO_n]<\frac{1}{n}$.  Thus, if $\sC_n$ is the complement $\sO_n$, then $\sC_n$ is nowhere dense, but $\mu[\sC_n]>1-\frac{1}{n}$.  If
$\sC=\Union_{n=1}^\oo\sC_n$, then $\sC$ is meager in $[0,1]$, but $\mu[\sC]=1$.  We can then construct a measurable function $f$ which is zero on the (co-meager) complement of $\sC$, but whose integral is arbitrarily large.

 Real analysis  is well-acquainted with  these sorts of pathologies, and  works around them.  But they can still be inconvenient in some applications of measure theory.  Note that the function $f$ constructed in the previous paragraph is discontinuous;  it is not possible to obtain this sort of pathological behaviour with a continuous function.  Likewise, the open sets $\sO_1,\sO_2,\ldots$ are not {\em regular} ---they are not the interiors of their own closures.   (Indeed, the only {\em regular} open dense subset of a topological space $\sS$ is $\sS$ itself.)   \ This suggests that, by confining our attention to regular open sets and continuous functions, we can develop a more ``well-behaved''  variant of measure theory on topological spaces.

Unfortunately, the family $\gR(\sS)$ of regular open subsets of $\sS$ is not closed under unions or complementation, so we cannot even {\em define} a classical measure if we confine its domain to $\gR(\sS)$.  But $\gR(\sS)$
{\em is}  a Boolean algebra under slightly different operations, and we can define a finitely additive, real-valued function on $\gR(\sS)$ with respect to this Boolean algebra structure.  We  call such a structure a {\em content}, to distinguish it from an ordinary measure.

  In some cases,  we might restrict attention to some Boolean subalgebra $\gB$ within $\gR(\sS)$.  For example, if $\sS=\Real^N$, then $\gB$ could be the Boolean algebra of regular open sets with piecewise smooth boundaries.    A content plays the role of a (finitely additive) measure on $\gB$.   But in most applications, we also need to compute the {\em integrals} of real-valued functions.  So such  a content would not be very useful, unless it came with a theory of integration.  One  goal  of this paper is to develop such a theory.   The other goal is to explicate the relationship between contents and  classical  Borel measures.

\paragraph{Applications.}  We focus on {\em finitely} additive measures because
 it is usually impossible to define a countably additive measure on $\gR(\sS)$ \cite[4XF(h,i), p.124]{Fremlin74}.   A reader who is used to countably additive measures may question the utility of  finitely additive ones.  But finitely additive measures
(or ``charges'') are more powerful and versatile than might first appear, as demonstrated  in the well-known monograph
\cite{BhaskaraRaoBook}.  Indeed,  Maharam \cite{Maharam76} quotes Bochner as saying, ``finitely additive measures are more interesting, and perhaps more important, than countably additive ones.''  Furthermore, there is a long tradition  
in decision theory and the philosophy of probability of working exclusively with finitely additive measures, as exemplified by the
seminal work of Savage  \cite{Savage}.  Dubins and Savage \cite{DubinsSavage65} argued that countable additivity is merely a convenient ``regularity condition'', while de Finetti \cite{deFinetti70,deFinetti95} went so far as to explicitly reject Kolmogorov's axiom of countably additive probability as not only unnecessary, but empirically meaningless.   

A decision-maker is called a {\em subjective expected utility maximizer} if she always chooses the option which maximizes expected utility, according to {\em some} utility function, and {\em some} probability measure over the possible states of the world.   Psychologically speaking, we can interpret the utility function as a representation of the decision-maker's ``desires'', while the probability measure represents her ``beliefs''.  Starting with the results of Savage and  de Finetti,  the literature of decision theory is replete with theorems showing that a decision-maker is a subjective expected utility maximizer if and only if her choices satisfy certain axioms of rationality or consistency.  But in these theorems, beliefs are usually represented with a {\em finitely} additive probability measure.   It is technically straightforward to obtain a {\em countably} additive measure, but this requires extra axioms which lack any clear interpretation ---either as a description of actual economic behaviour, or as a normative standard  of ideal rationality.   Thus, these so-called ``technical'' hypotheses are  not very appealing to either economists or philosophers.  Fortunately, finitely additive measures are sufficient  for most purposes in decision theory.

  In fact, the present paper grew out of an ongoing research project in decision theory.  Most papers in decision theory assume that the space  $\sS$   of possible states of the world is an abstract measurable space, and any {\em measurable} function can be interpreted as a feasible ``plan of action'' by the decision-maker.  But in many contexts,  $\sS$   has a natural topology, and most measurable functions are physically meaningless;  only {\em continuous} functions  on $\sS$   represent feasible plans of action.  Furthermore, most papers in decision theory assume that any measurable subset of  $\sS$   can play the role of an {\em event} ---that is, a piece of information the decision-maker can observe and upon which she can conditionalize her plans.  But in many contexts, decision-makers face considerable informational limitations;  they observe the world ``through a glass, darkly'', using imprecise measurement instruments which collapse a continuum of possible inputs into a finite set of possible outputs.    
    
     In two companion papers, we have developed a version of decision theory which takes these limitations seriously \cite{SEU_continuous,ImperfectPerception}.    In these papers, feasible plans of action correspond to   {\em continuous} functions  on $\sS$ .  The only observable events are those in a Boolean subalgebra of regular open subsets of  $\sS$ .  To even {\em define} a subjective expected utility representation in such a context, we need a theory of measure and integration for such Boolean algebras ---this is the origin of the material in Sections \ref{S:notation}, \ref{S:integration}, and \ref{S:measurable} of the present paper (described in detail below).\footnote{In \cite{SEU_continuous,ImperfectPerception}, we refer to contents as {\em credences}, since they represent the decision-maker's beliefs.}     But it is also necessary to understand the relationship between a representation based on contents, and more classical representations based on Borel probability measures ---this is the origin of the material in Sections
     \ref{S:liminal}, \ref{S:liminal.compactification}, and \ref{S:gleason}.  Interestingly, in \cite{SEU_continuous}, the liminal structures introduced in Section \ref{S:liminal} have an appealing interpretation in terms of measurement instruments with random errors.

The paper \cite{Intertemporal} adapts the results of \cite{SEU_continuous,ImperfectPerception} to setting of intertemporal choice.  Here, $\sS=[0,\oo)$ represents a time continuum, and a continuous function from $\sS$ into some other space $\sX$ represents a {\em trajectory}.  Preferences over such trajectories are called {\em intertemporal preferences};  \cite{Intertemporal} shows that if these satisfy certain axioms, then they admit a {\em discounted utility integral} representation. Formally, this looks like a subjective expected utility representation, but now the measure represents the agent's intertemporal discounting rather than her beliefs.  In this setting, the integral representations introduced in Section \ref{S:liminal} of the current paper can be simplified to Riemann-Stieltjes integrals.

Finally, \cite{Spectral} applies the results of \cite{SEU_continuous,ImperfectPerception} to setting where the agent has preferences over a Riesz space or a real commutative Banach algebra $\sA$.  The state space $\sS$ arises as the spectrum of  $\sA$ (via the Kakutani and Gelfand representation theorems, respectively).  The Boolean algebra of regular open subsets of $\sS$ is isomorphic to the Boolean algebra of {\em bands} of $\sA$ (in the case of a Riesz space) or {\em regular ideals}  of $\sA$ (in the case of a Banach algebra); these play the role of the {\em events} observable by the agent.

 \paragraph{Organization and main results.}  
    Section \ref{S:literature} briefly reviews prior literature.   In Section \ref{S:notation}, we define terminology, give key examples and state some basic results about contents.
 In particular, we show that, for any Baire space $\sS$,
 there is a natural bijective correspondence between  contents on $\gR(\sS)$ and  {\em residual charges} on $\sS$  ---that is, finitely additive measures that annihilate all meager subsets (Proposition \ref{baire.algebra}) .
  
    In Section \ref{S:integration}, we develop a theory of integration for contents.  The main result is
    Theorem \ref{from.probability.to.expectation0}: on any topological space $\sS$  and  any
 Boolean subalgebra $\gB$ of $\gR(\sS)$,  there is  a unique {\em integrator}: an operator which returns a real number $\integral_{\sB} g \dmu$ for any ``integrable''  real-valued function $g$  on $\sS$,  any content $\mu$ on $\gB$, and  any  $\sB$ in $\gB$.  This theorem is crucial in applications;
 for example, if we interpret $ \mu $ as a probability, then $\integral_{\sB} g \dmu$ is the ``expected value'' of $g$, conditional on  event $\sB$.  But its usefulness would be limited if we lacked technical tools to analyze such integrals.  Thus,  Section \ref{S:measurable} develops such tools.  The most important one is a ``change of variables'' theorem for integrators (Proposition \ref{change.of.variables}). 
  
  Section \ref{S:liminal} and the following  sections contain the most interesting and important results of the paper;  they  
explore the relationship between contents and Borel measures.  Na\"ively, we might think that a content on $\gR(\sS)$
is equivalent to a Borel measure on $\sS$.  But this is not the case (see Nonexample \ref{X:lebesgue.is.not.content}).  
The problem is that the join operation $\vee$ on $\gR(\sS)$ treats the boundaries of regular sets differently
than the ordinary union operation $\union$.  Thus, a content can allocate nonzero mass to these boundaries in a way that  cannot be fully described by a Borel measure.  To describe these boundary allocations, we need an extra informational component, which we call a {\em liminal structure}.  The main result of Section \ref{S:liminal} says that,
if $\sS$ is a locally compact $T_4$ space, then any content $\mu$ on $\gR(\sS)$ can be represented in a unique way as a combination of a Radon measure $\nu$ and a liminal structure (Theorem \ref{content.vs.Riesz.charge2}).  Furthermore, if $g:\sS\into\Real$ is continuous and decays to zero at $\oo$, then for any  $\sR\in\gR(\sS)$, the integral
 $\integral_{\sR} g \dmu$  is the sum of a Lebesgue integral with respect to $\nu$, and another integral with respect to the liminal structure.   If $\sS$ is a {\em compact} Hausdorff space, then this statement holds for {\em all} continuous real-valued functions on $\sS$ (Corollary \ref{content.vs.Riesz.charge3}).   We also have a similar result in the case
 when $\sS$ is a $T_4$ space which is {\em not} necessarily locally compact, but in this case, instead of a Radon measure,
 the representation uses a normal charge $\nu$ defined on the Boolean algebra generated by open sets (Theorem \ref{content.vs.Riesz.charge}).  If $g:\sS\into\Real$  is bounded and continuous, then 
$\integral_{\sR} g \dmu$ is the sum of an integral with respect to $\nu$ and an integral with respect to the
  liminal structure.
 
   Section \ref{S:liminal.compactification} extends the representation theorems of Section \ref{S:liminal} to locally compact Hausdorff spaces.  The main result (Theorem \ref{content.vs.Riesz.charge4}) says that, if $\sS$ is a locally compact Hausdorff space,
   and $\barsS$ is a compactification of $\sS$, then any content on $\sS$ can be represented in a unique way as a combination of a Borel measure $\barnu$ on $\barsS$ and a liminal structure on $\barsS$.  Furthermore, if $g:\sS\into\Real$ is  continuous and can be extended continuously to $\barsS$, then for any  $\sR\in\gR(\sS)$, the integral
 $\integral_{\sR} g \dmu$  is the sum of Lebesgue integrals with respect to $\barnu$ and the  liminal structure.  In particular,
 if $\barsS$ is the Stone-\v{C}ech compactification, this holds for {\em all} bounded continuous real-valued functions on $\sS$.
 
    Given any Boolean algebra $\gB$, one can construct a unique topological space $\sigma(\gB)$ such that $\gB$ is isomorphic to the algebra of clopen subsets of $\sigma(\gB)$;  this is called the {\em Stone space} of $\gB$.  In particular, if $\sS$ is a  compact Hausdorff space, and $\sS^*$ is the Stone space of  $\gR(\sS)$, then there is a canonical continuous surjection from $\sS^*$ to $\sS$;  this is called the {\em Gleason cover} of $\sS$.  This  suggests  another way to represent contents, which we develop in Section \ref{S:gleason}.  First, we generalize the Gleason construction in two ways:  given any {\em locally} compact Hausdorff space $\sS$, and any large enough {\em subalgebra} $\gB$ of $\gR(\sS)$, we construct
 a canonical continuous surjection from the Stone space $\sS^*$ of $\gB$ to the Stone-\v{C}ech compactification of $\sS$ (Proposition \ref{generalized.gleason}).    Any $\sB\in\gB$ corresponds to a unique clopen subset $\sB^*\subseteq\sS^*$,
 and any bounded continuous  $g:\sS\into\Real$ has  a unique lift  $g^*:\sS^*\into\Real$.   We then  show that any content  $\mu$ on $\gB$ is
 represented by a unique Borel measure $\mu^*$ on $\sS^*$, such that $\mu^*[\sB^*]=\mu[\sB]$ for all $\sB\in\gB$.
 Furthermore, if $g$ is any bounded continuous real-valued function $\sS$  that is ``measurable'' relative to $\gB$, then $ \integral_{\sB} g \dmu  = \int_{\sB^*} g^*\dmu^*$ (Theorem \ref{generalized.gleason}). 
 
   Section \ref{S:cat} gives a functorial formulation to the constructions in  Sections  \ref{S:liminal.compactification} and \ref{S:gleason}.  We define a {\dfn content space} to be a locally compact Hausdorff space $\sS$ equipped with a content on
  some Boolean subalgebra of $\gR(\sS)$.  These spaces form a category.
   If we apply Theorem \ref{content.vs.Riesz.charge4} to the Stone-\v{C}ech compactification, we obtain
 a faithful functor from  the category of content spaces to the subcategory of 
  {\em compact} content spaces  (Proposition \ref{stone.cech.compactification.functor}).
  Meanwhile, the construction in Theorem \ref{generalized.gleason} yields a faithful functor from the category of content spaces
  to the category of Stonean spaces equipped with Borel measures (Proposition \ref{stone.space.functor}).
  Finally the Gleason covering map from  Proposition \ref{generalized.gleason} is a natural transformation from the first
  of these functors to the second (Proposition \ref{gleason.natural.transform}).
 The following diagram shows the logical dependencies between these sections.

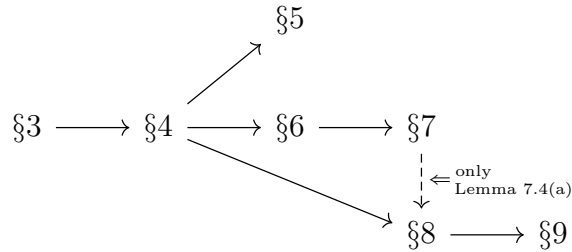
\[
\begin{tikzcd}
 ~ & ~ &\S\ref{S:measurable} \\
\S\ref{S:notation} \arrow{r} & \S\ref{S:integration}  \arrow{ur} \arrow{r} \arrow{drr} 
 						& \S\ref{S:liminal} \arrow{r} & \S\ref{S:liminal.compactification}\arrow[dashed]{d}{\Longleftarrow \; {\mbox{\tiny  only Lemma  \ref{stone.cech.boolean.isomorphism} }}} \\
						& & & \S\ref{S:gleason}\arrow{r} &  \S\ref{S:cat}
\end{tikzcd}
\]

\section{Prior literature\label{S:literature}}

The theory of measure and integration on topological spaces is  well-established;  see \cite{Wheeler}, \cite{Bogachev} and  \cite{Fremlin4I,Fremlin4II} for excellent surveys.  But the present paper has surprisingly little connection with this theory.  The existing literature deals with measures defined on the Borel or Baire sigma-algebras of a topological space (or suitable sub-algebras thereof), whereas we are concerned with the Boolean algebra of regular open sets. 

Sections \ref{S:liminal} and \ref{S:liminal.compactification} of the present paper show how to represent the 
non-classical integration  operator $\integral$ from Section \ref{S:integration} in terms of a classical integral with respect to
a normal Borel measure.   This is reminiscent of the Riesz and Alexandroff Representation Theorems, and hence the literature on {\em strict topologies} \cite{Wheeler}.  This literature investigates the following question:  given a topological space $\sS$ and a topological vector space $\sH$ of real-valued continuous functions on $\sS$, is there some space $\sM$ of measures on $\sS$ (or on some superspace, e.g. the Stone-\v{C}ech compactification of $\sS$) such that for every bounded linear function $I:\sH\into\Real$, there is a unique $\mu\in\sM$ such that $I(h)=\int_\sS h\dmu$ for all $h\in\sH$?      But our results in  Sections \ref{S:liminal} and \ref{S:liminal.compactification} are unrelated to this literature for a simple reason:  we are also interested in integrals of the form $\integral_\sR h$  where $\sR$ is any  regular open {\em subset} of $\sS$.  To represent such an integral, we need not only a Borel measure, but also a ``liminal structure'' which describes the distribution of measure along the boundary of $\sR$.  As far as we know, there is no analogous structure in the literature on strict topologies.

The Boolean algebra of regular open subsets was introduced  by Tarski \cite{Tarski37}.  The first analysis of finitely additive measures on Boolean algebras was by Horn and Tarski \cite{HornTarski48},  but this paper did not specifically consider the algebra of regular open sets.  The literature inspired by  \cite{HornTarski48} has focused mainly on identifying necessary and sufficient conditions for abstract
Boolean algebras to support finitely additive measures with particular properties; see \cite{DzamonjaPlebanek08} and \cite{Borodulin-NadziejaDzamonja13} for recent results and a review of this literature.

As we will soon see (Proposition \ref{baire.algebra}), there is a close relationship between contents and
{\em residual measures} ---that is, measures on the sigma-algebra of Baire-property subsets of a topological space  which vanish on all meager subsets of that space.   Residual measures were studied by  Armstrong and Prikry \cite{ArmstrongPrikry78},  Flachsmeyer and Lotz
\cite{FlachsmeyerLotz78,FlachsmeyerLotz80,FlachsmeyerLotz80b,Lotz82} and Zindulka \cite{Zindulka2000}.\footnote{ Flachsmeyer and Lotz called them {\em hyperdiffusive} measures.}   
 But these papers focused on {\em countably} additive measures, whereas we are interested in {\em finitely} additive ones (which we call residual {\em charges}).  This is important, because countably additive residual measures are much harder to construct than finitely additive ones.   Nevertheless, there are parallels between the finitely additive and countably additive cases.  For example, Armstrong and Prikry \cite[Proposition 7]{ArmstrongPrikry78} observe that any residual measure on a space $\sS$ can be represented by a measure on the {\em Gleason cover} of $\sS$ (i.e. the Stone space of $\gR(\sS)$).   Theorem \ref{gleason.covering} of the present paper makes a similar statement for contents defined on an arbitrary subalgebra $\gB$ of $\gR(\sS)$.   Thus, by combining Proposition \ref{baire.algebra} and Theorem \ref{gleason.covering} in the special case when $\gB=\gR(\sS)$, we obtain a version of Armstrong and Prikry's  Proposition 7 for residual charges.

   Let $\gO(\sS)$ be the lattice of open subsets of a topological space $\sS$.  A 
  content on $\gR(\sS)$ might seem superficially similar to a {\em valuation} defined on $\gO(\sS)$.\footnote{ A {\dfn valuation} is a function $\nu:\gO(\sS)\into\Real_+$ such
  that $\nu[\emptyset]=0$ and for any $\sO,\sQ\in\gO(\sS)$, we have  $\nu[\sO \union \sQ] = \nu[\sO]+\nu[\sQ]-\nu[\sO\intsct\sQ]$, and also $\nu[\sO]\leq\nu[\sQ]$ if $\sO\subseteq\sQ$.}  But there are two important differences:  first, contents are only defined on {\em regular} open sets; second, the additivity property for a content is defined with respect to a special join operation $\vee$ (see Section \ref{S:notation}), whereas the additivity property of a valuation is defined with respect to the standard set-theoretic operations of union and intersection.  These differences have two important consequences.    First, there is now an integration theory for  valuations  \cite{Edalat95,Edalat95b,EdalatNegri98,Howroyd2000,LawsonLu03,Lawson04}.\footnote{ To be precise, this theory establishes a relationship between  the Lebesgue integral with respect to a Borel measure on a space $\sX$ and
  an integral with respect to a valuation on its {\em upper space} $\sU(\sX)$.}
    But as far as we know, there is no comparable theory of integration for contents on $\gR(\sS)$ or its subalgebras;  we shall develop one in Section \ref{S:integration}.  
  Second, under fairly general conditions, a valuation on $\gO(\sS)$ can be uniquely extended to a charge defined on the Boolean algebra generated by $\gO(\sS)$ \cite{Smiley44,HornTarski48}, or even to a Borel measure \cite{Lawson82,Edalat95b,AlvarezEdalatSaheb2000,Alvarez-Manilla02,KeimelLawson05}.  But as we shall see in 
  Section \ref{S:liminal}, the corresponding results for contents on $\gR(\sS)$ are much more subtle, because they
  involve not only a Borel measure but also a {\em liminal structure}, which, roughly speaking, describes the way the content deals with the {\em boundaries} of regular open sets.

\section{Regular open algebras and contents\label{S:notation}}
\setcounter{equation}{0}

Throughout this paper, let $\sS$ be a topological space.  For any subset $\sA\subseteq\sS$, let
$\Int(\sA)$ denote its interior, let $\Cl(\sA)$ denote its closure, and let $\partial\sA$ denote its boundary.  An open subset $\sO\subseteq\sS$ is
 {\dfn regular} if $\sO=\Int[\Cl(\sO)]$.  
 For example, the interior of any closed subset of $\sS$ is a regular set.
   Clearly, the intersection of two regular open subsets is a regular open subset. Given any two regular open subsets $\sQ,\sR\subseteq\sS$, we define $\sQ\vee\sR:=\Int[\Cl(\sQ\union\sR)]$; this is  the smallest regular open subset of $\sS$ containing both $\sQ$ and $\sR$.
  For example, if $\sS=\Real$, then $(0,1)\vee(1,2) = (0,2)$.
Meanwhile, we define $\neg\sD:=\Int[\sS\setminus\sD]$, which is 
 another regular open subset. The set  $\gR(\sS)$  of all regular open subsets of $\sS$ forms a  Boolean algebra under  the operations $\intsct$, $\vee$, and $\neg$ \cite[Theorem 314P]{Fremlin}.
A subcollection $\gB\subseteq\gR(\sS)$ is a {\dfn Boolean subalgebra} if $\gB$ is closed under the operations  $\intsct$, $\vee$, and $\neg$.  

\example{\label{X:algebra} (a) A subset $\sE\subseteq\Real$ is {\dfn elementary} if
 $\sE:=(a_1,b_1)\disj (a_2,b_2)\disj\cdots\disj(a_N,b_N)$ for some
  $-\oo\leq a_1<b_1<a_2<b_2<\cdots<a_N<b_N\leq \oo$.   Any elementary subset is open and regular.  Let $\gE$ be the
  collection of all elementary subsets of $\Real$;  then $\gE$ is a Boolean subalgebra of the algebra of
  regular subsets of $\Real$.  (Here, $\disj$ indicates disjoint union.)

\item (b) Suppose $\sS$ is a differentiable manifold.  A subset $\sH\subseteq\sS$ is a {\dfn smooth hypersurface} if there
  is a differentiable function $\phi:\sS\into\Real$ such that $\sH:=\phi^{-1}\{r\}$ for some $r\in\Real$,
  and such that  d$\phi(h)\neq 0$ for all $h\in\sH$.  We will say that a regular open subset $\sR\subseteq\sS$ has a  {\dfn piecewise smooth boundary}
  if there is a finite collection $\sH_1,\sH_2,\ldots,\sH_N$ of smooth hypersurfaces such that
  $\partial\sR = (\sH_1\intsct\partial\sR)\union\cdots\union(\sH_N\intsct\partial\sR)$.
  Let $\gB_{\mathrm{smth}}$ be the set  of regular open sets with piecewise smooth boundaries; then $\gB_{\mathrm{smth}}$ is a Boolean subalgebra of $\gR(\sS)$.

\item (c)  Suppose $\sS$ is a topological vector space.  A subset $\sH\subseteq\sS$ is a {\dfn hyperplane} if there
  is a continuous linear function $\phi:\sS\into\Real$ such that $\sH:=\phi^{-1}\{r\}$ for some $r\in\Real$.
   A regular open subset $\sR\subseteq\sS$ is a {\dfn polyhedron}
  if there is a finite collection $\sH_1,\sH_2,\ldots,\sH_N$ of hyperplanes such that
   $\partial\sR = (\sH_1\intsct\partial\sR)\union\cdots\union(\sH_N\intsct\partial\sR)$.
  (Heuristically,  $\sH_n\intsct\partial\sR$ is a ``face''  of $\sR$.  Note that we
  do not require $\sR$  to be convex, or even connected.)\   Let $\gB_{\mathrm{poly}}$ be the set of regular polyhedra;
  then $\gB_{\mathrm{poly}}$ is a Boolean subalgebra of $\gR(\sS)$. 
  
  \item  (d) Let $\sS\subseteq\Real^N$ be an open set.  Let $\gB_{\mathrm{jor}}(\sS)$ be the set of all regular open subsets of $\sS$ whose boundaries have Lebesgue measure zero.  This is a Boolean subalgebra of $\gR(\sS)$, which is sometimes called the {\dfn Jordan algebra}. 
    (Note that $\gB_{\mathrm{poly}}(\Real^N)\subset \gB_{\mathrm{smth}}(\Real^N)\subset\gB_{\mathrm{jor}}(\Real^N)$.)
   }

\noindent A  {\dfn content} on $\gB$ is a function $\mu:\gB\into[0,\oo)$, such that for any finite collection $\{\sB_n\}_{n=1}^N$ of disjoint elements
of $\gB$,
 we have
\beqn
\label{prob.measure}
\mu\lb[\bigvee_{n=1}^N \sB_n \rb] \quad = \quad 
\sum_{n=1}^N \mu[\sB_n].
\eeqn
Note an important difference from the usual definition of a measure:   additivity is defined with respect to the operation $\vee$, rather than ordinary union.  We say that  $\mu$ has {\dfn full support} if $\mu[\sB]>0$ for all nonempty $\sB\in\gB$.  

For simplicity, throughout this paper we will assume that $\mu[\sS]=1$ ---in other words, the content acts like a ``probability measure''.  But all our results easily generalize to the case when $0<\mu[\sS]<\oo$, through a simple rescaling.

\example{\label{X:elementary.lebesgue.content}
Let $\sS:=(0,1)$, and let $\gE$ be the Boolean algebra of elementary regular open subsets of $(0,1)$, from Example \ref{X:algebra}(a).   Let $\phi:\sS\into\Real$ be continuous and nondecreasing.  
For any $\sE\in\gE$, if $\sE=(a,b)$ for some $a<b$, then
define $\mu[\sE]:= \phi(b)-\phi(a)  $.   If $\sE=\sE_1\disj\cdots\disj\sE_N$ for some
disjoint open intervals $\sE_1,\ldots,\sE_N$, then define
$\mu[\sE]:=\mu[\sE_1]+\cdots+\mu[\sE_N]$.   Then $\mu$ is a content on $\gE$.   If $\phi$ is strictly increasing,
then $\mu$ has full support. 
}

 If  $\phi(s)=s$ for all $s\in(0,1)$, then  the content in Example \ref{X:elementary.lebesgue.content} is  the restriction of
the Lebesgue measure $\lam$ to $\gE$.   It is tempting to extrapolate  that we can
obtain a content on all of $\gR[0,1]$ by restricting $\lam$ to regular open sets.  But
this is not the case.

\refstepcounter{thm}\label{X:lebesgue.is.not.content}
\paragraph{Nonexample \thethm.}
{\em The Lebesgue measure  restricted to $\gR[0,1]$ is not  a content. }
 
\bthmprf[Proof sketch.\footnote{We are grateful to  Joel David Hamkins for showing us this construction.}]
Let $\lam$ be the Lebesgue measure, and let
 $\sK$ be a ``fat'' Cantor set ---that is, a closed, nowhere dense subset of $[0,1]$ with $\lam[\sK]>0$.  Thus, $\sU:=[0,1]\setminus\sK$ is an open dense subset of $[0,1]$ with $\lam[\sU]<1$.
The set $\sU$ is a countable disjoint union of open intervals, and  is
not regular.  However, it can be divided into two pieces, $\sL$ and $\sR$,
which are defined by taking the ``left half'' and ``right half'' of each of open intervals
comprising $\sU$.  These {\em are} regular open sets, and $\sR=\neg\sL$.   Thus, $\sL\vee\sR=[0,1]$.
But clearly, $\lam[\sL]+\lam[\sR]=\lam[\sU]<1=\lam[0,1]$.   Thus, $\lam$ is not
finitely additive on $\gR[0,1]$.  
(For a complete proof of a more general result, see Proposition \ref{X:lebesgue.is.not.content2} below.)
\ethmprf 
 
 \noindent Some contents are very different from the ones in Examples \ref{X:elementary.lebesgue.content}.

\example{\label{X:ultrafilter.example2} 
Let $\gB\subseteq\gR(\sS)$ be a Boolean subalgebra.  An {\dfn ultrafilter} is a collection
$\gU\subseteq\gB$ such that:
(a) If $\sU,\sV\in\gU$, then $\sU\intsct\sV\in\gU$;  \ (b) If $\sU\in\gU$ and $\sU\subseteq\sV$, then $\sV\in\gU$; \
(c) $\emptyset\not\in\gU$; \ and   
(d)  For any $\sB\in\gB$, either $\sB\in\gU$, or $\neg\sB\in\gU$, but not both.
Given any ultrafilter $\gU$, we can define a function $\del_\gU:\gB\into\{0,1\}$ as follows: for any $\sB\in\gB$, $\del_\gU[\sB]=1$ if $\sB\in\gU$, whereas $\del_\gU[\sB]=0$ if $\sB\not\in\gU$.   It is easy to verify that $\del_\gU$ is a content.   

An ultrafilter $\gU$ is {\dfn fixed} if there is some point $s\in\sS$ such that
$s\in\sU$ for all $\sU\in\gU$.     To obtain such an ultrafilter, let $\gF_s:=\{\sR\in\gR$; \ $s\in\sR\}$.
Then $\gF_s$ is a {\em filter} (i.e. it satisfies properties (a), (b) and (c) above.)  The Ultrafilter Theorem says that
there is an ultrafilter $\gU\supseteq\gF_s$ \cite[Theorem 2A1O]{FremlinII}. 
  It is easy to see that $\gU$ is fixed at $s$. 
 In this case, $\del_\gU$ can be interpreted as a ``point mass'' at $s$.  But there is a complication:  $\gU$ is not unique.  There are many different ultrafilters fixed at $s$, and each one defines a different content ---a different version of the ``point mass'' at $s$.  We will return to this in Example
 \ref{X:charge.vs.measure1}.
 
 If $\gU$ is not fixed, then it is {\dfn free}.  A free ultrafilter can behave like  a point mass ``at infinity''
 (see Example \ref{X:ultrafilter.example} below), or a point mass on the metaphorical ``boundary'' of the space $\sS$ (Example \ref{X:bad.content}).  But it can also behave like a point mass in the interior
 of $\sS$ ``approached from one side'' (Example \ref{X:charge.vs.measure1b}).
 }

 Let $\Borel(\sS)$ be the Borel sigma-algebra on $\sS$ ---that is, the smallest sigma-algebra containing all open and closed subsets of $\sS$.      Observe that $\gR(\sS)\subseteq\Borel(\sS)$ as sets, but the Boolean algebra operations are different.
Nevertheless, the next result says that contents on $\gR(\sS)$  often  have a nice representation on $\Borel(\sS)$. 
 Recall that a subset $\sN\subseteq\sS$ is {\dfn nowhere dense} if $\Int[\Cl(\sN)]=\emptyset$.
A subset $\sM\subseteq\sS$ is {\dfn meager} if it is a countable union of nowhere dense sets. 
Let $\gM(\sS)$ be the set of all meager subsets of $\sS$;  then $\gM(\sS)$ is an ideal under
the standard set-theoretic operations. (That is: the union of two meager sets is meager, and the intersection
of a meager set with any other set  is meager.)
 Let $\gB$ be a Boolean algebra of subsets of  $\sS$ (for example, $\Borel(\sS)$). 
  A {\dfn probability charge} on $\gB$ is a function $\nu:\gB\into[0,1]$ 
such that (1) $\nu[\sS]=1$ and  (2)
 $\nu[\sA\disj\sB]=\nu[\sA]+\nu[\sB]$ for any disjoint $\sA,\sB\in \gB$.    Suppose that
 $\gM(\sS)\subseteq\gB$.  We will 
 say that $\nu$ is a {\dfn  residual  charge} if $\nu[\sM]=0$ for all $\sM\in\gM(\sS)$.

 A topological space $\sS$ is a {\dfn Baire space} if the intersection of any countable family
of open dense sets is dense.  In particular, any locally compact Hausdorff space is Baire, and any
completely metrizable space is Baire \cite[Corollary 25.4]{Willard}.

\Proposition{\label{residual.Borel.charge}}
{
Let $\sS$ be a  Baire  space.  Then
there is a bijective correspondence between the contents on $\gR(\sS)$ and
the  residual  charges on $\Borel(\sS)$.  To be precise, if $\nu$ is a  residual  charge on $\Borel(\sS)$, then we
get a content  by  restricting $\nu$ to $\gR(\sS)$.  Every content arises
in this fashion, and no two  residual  charges produce the same content.
}

In fact, Proposition \ref{residual.Borel.charge} is a special case of a more general result.\footnote{We are grateful to Robert Furber for pointing out this result to us.}
 A subset $\sB\subseteq\sS$ has the {\dfn Baire property} if $\sB=\sO\symdif \sM$ for some
open $\sO\subseteq\sS$ and meager $\sM\subseteq\sS$.  Let $\gB\!\ga(\sS)$ be the collection of all subsets with the Baire property; then $\gB\!\ga(\sS)$ is  a sigma-algebra under the standard set-theoretic operations.\footnote{Note that $\gB\!\ga(\sS)$ is completely unrelated to the Baire sigma-algebra.}
  
\Proposition{\label{baire.algebra}}
{
Let $\sS$ be a Baire space.  Then
there is a bijective correspondence between the contents on $\gR(\sS)$ and
the  residual  charges on $\gB\!\ga(\sS)$.  To be precise, if $\nu$ is a  residual  charge on $\gB\!\ga(\sS)$, then we
can obtain a content  by simply restricting $\nu$ to $\gR(\sS)$.  Every content arises
in this fashion, and no two  residual  charges produce the same content.
} 
\bthmprf
For any $\sB_1,\sB_2\in\gB\!\ga(\sS)$, write $\sB_1\sim\sB_2$ if $\sB_1\symdif\sB_2\in\gM(\sS)$;
this is an equivalence relation on $\gB\!\ga(\sS)$.  Let $\gA:=\gB\!\ga(\sS)/\sim$; then
the Boolean algebra operations on $\gB\!\ga(\sS)$ factor through to $\gA$ (because 
$\gM(\sS)$ is an ideal), making $\gA$  a Boolean algebra.\footnote{$\gA$ is called the {\dfn category algebra} of $\sS$.}  Thus, a  residual  charge is equivalent to a finitely additive function on the Boolean algebra $\gA$.

Recall that $\gR(\sS)$ is a subset (but not a subalgebra) of $\gB\!\ga(\sS)$.
Let $\pi:\gB\!\ga(\sS)\into\gA$ be the quotient map, and let $\phi$ be the restriction of $\pi$ to $\gR(\sS)$.   Then $\phi$ is a Boolean algebra homomorphism from $\gR(\sS)$ to $\gA$ \cite[514I(b), p.44]{Fremlin5I}.  Furthermore, if $\sS$ is Baire, then
 $\phi$ is bijective, and hence, an isomorphism \cite[514I(f)]{Fremlin5I}.
 In other words, every $\sim$-equivalence class in $\gB\!\ga(\sS)$ contains a unique representative
from $\gR(\sS)$.  

Thus, if $\mu$ is any content on $\gR(\sS)$, then we obtain a finitely additive function
$\mu\circ\phi^{-1}:\gA\into[0,1]$, and from there, a  residual  charge $\mu\circ\phi^{-1}\circ\pi:\gB\!\ga(\sS)\into[0,1]$.  Conversely, given any  residual  charge $\nu$ on $\gB\!\ga(\sS)$, we obtain
a finitely additive function $\nu^*$ on $\gA$, and from there, a content $\nu^*\circ\phi$ on $\gR(\sS)$.
\ethmprf

\bthmprf[Proof of Proposition \ref{residual.Borel.charge}.]
  $\Borel(\sS)$ is a sigma-subalgebra of  $\gB\!\ga(\sS)$, because $\gB\!\ga(\sS)$ contains every open subset of $\sS$,
 and hence, the sigma-algebra generated by open sets, which is $\Borel(\sS)$.  
   Thus, every (residual) charge on $\gB\!\ga(\sS)$  defines a (residual) charge on $\Borel(\sS)$.
Conversely, if we define the relation $\sim$ as in the proof of Proposition \ref{baire.algebra}, then every $\sim$-equivalence class in $\gB\!\ga(\sS)$ contains an element of  $\Borel(\sS)$, because $\gR(\sS)\subseteq\Borel(\sS)$, and as we have already noted
in the proof of Proposition \ref{baire.algebra}, every $\sim$-equivalence class in $\gB\!\ga(\sS)$ contains an element of  $\gR(\sS)$.  Thus, every {\em residual} charge on  $\Borel(\sS)$ admits a unique extension to a residual charge on $\gB\!\ga(\sS)$.
Thus, there is a bijective correspondence between the set of residual charges on $\gB\!\ga(\sS)$  and the set of residual charges on $\Borel(\sS)$.  The result now follows from Proposition \ref{baire.algebra}.
\ethmprf

\noindent Let $s\in\sS$, and let $\gU$ be an ultrafilter fixed at $s$.
The ``point mass'' content $\del_\gU$ from Example \ref{X:ultrafilter.example2} seems to contradict Proposition \ref{baire.algebra}, since the singleton $\{s\}$ is obviously meager.
  But this is misguided.   If $\gU$ is an ultrafilter fixed at $s$, then $\del_\gU$ is equivalent to  residual  charge which gives probability 1 to every open neighbourhood of $s$, but gives {\em probability zero} to $\{s\}$ itself.  This is possible because charges are only {\em finitely} (and not {\em countably}) additive.

\breath

This raises the question:  why do we confine our attention to {\em finitely} additive contents?  The set $\gR(\sS)$ is actually a {\em $\sigma$-Boolean algebra} under the following operations.\footnote{\label{complete.boolean.algebra.footnote}
Indeed, $\gR(\sS)$ is a {\em complete} Boolean algebra:  it is closed even under the uncountably infinite versions of
the operations $\bigvee$ and $\bigwedge$; see \cite[Theorem 314P]{Fremlin} or \cite[Proposition 2.3, p.45]{Walker}.}  Given any countable collection $\{\sR_n\}_{n=1}^\oo$ of regular subsets, we define
\beqn
\bigvee_{n=1}^\oo \sR_n \quad = \quad \Int\lb[\Cl\lb(\Union_{n=1}^\oo \sR_n\rb)\rb]
\quad\And\quad
\bigwedge_{n=1}^\oo \sR_n\quad:=\quad \Int\lb[\Intsct_{n=1}^\oo \Cl(\sR_n)\rb].
\label{sigma.regular}
\eeqn
 However, if $\sS$ is a perfect,
 second-countable  Hausdorff space, then there is no nontrivial countably additive content  on  $\gR(\sS)$ \cite[4XF(h,i), p.124]{Fremlin74}.

\section{Integration\label{S:integration}}
 \setcounter{equation}{0}

Let $\gB$  be a Boolean subalgebra of $\gR(\sS)$, and 
let $\sR\in\gB$.  A {\dfn $\gB$-partition} of $\sR$ is a  collection $\{\sB_1,\ldots,\sB_N\}$ (for some $N\in\Natur$) of disjoint elements of $\gB$ such that $\sR=\sB_1\vee\cdots\vee\sB_N$.   (For instance, if $\sS=\Real$, then $\{(0,1),(1,2)\}$ is an  $\gR(\sS)$-partition of $(0,2)$.)     Let $\sM$ be a set of contents on $\gB$.
   Let $\sH$ be a vector space of bounded, real-valued functions on $\sS$.    
   Assume $\sH$ contains $\bone$, the constant function with value 1.
 An {\dfn integrator}
     is a function $\integral:\gB\x\sH\x\sM\into\Real$ with the following properties.
\begin{description}
   \item[(I1)]  For any $\sB\in\gB$ and $\mu\in\sM$, the function $\integral(\sB,\bullet,\mu):\sH\into\Real$
 is a bounded linear functional  and is {\em nondecreasing}:
 for any $f,g\in\sH$, if $f(s)\leq g(s)$ for all $s\in\sS$, then $\integral(\sB,f,\mu)\leq \integral(\sB,g,\mu)$.
Also,  $\integral(\sB,\bone,\mu)=\mu[\sB]$,
while  $\integral(\sB,h,\mu)=0$ if $ h_{\restr\sB}=0$.
 
   \item[(I2)] For every non-negative $ h\in\sH$ and $\mu\in\sM$, the function $\integral(\bullet,h,\mu):\gB\into[0,\oo)$  is a  content on $\gB$.  In other words,
       for any   $\gB$-partition $\{\sB_n\}_{n=1}^N$ of $\sB$,
           \beqn
    \label{conditional.expectation0}
     \mbox{$\integral$}(\sB,h,\mu) \quad=\quad  \sum_{n=1}^N \mbox{$\integral$}(\sB_n,h,\mu) .
\eeqn
     \item[(I3)] For any $\sB\in\gB$ and $h\in\sH$, the function $\integral(\sB,h,\bullet):\sM\into\Real$
 is  linear.  In other words, for any $\mu_1,\mu_2\in\sM$ and $r_1,r_2\in\Real$,
 if $r_1\,\mu_1+r_2\,\mu_2\in\sM$, then  $\integral(\sB,h,r_1\,\mu_1+r_2\,\mu_2)=
 r_1\,\integral(\sB,h,\mu_1)+r_2\,\integral(\sB,h,\mu_2)$.
     
\end{description}
We will indicate $\integral(\sB,h,\mu)$ with the expression 
\[\Integral_\sB h \dmu.\]
The slash should prevent confusion with the standard Lebesgue integral.

\example{Let $\sS=(0,1)$, let $\gE$  be the Boolean algebra of elementary regular open subsets of $(0,1)$,
 and let $\sH$ be the set of all bounded, continuous, real-valued functions on $\sS$.
  For any continuous, increasing function $\phi:\sS\into\Real$, let
   $\mu_\phi$  be the content on $\gE$ from Example \ref{X:elementary.lebesgue.content}.
 Let $\sM$ be the set of all contents obtained in this way.
 We obtain an integrator  as follows:  for any $\sE=(a_1,b_1)\disj
(a_2,b_2)\disj
\cdots \disj (a_N,b_N)$ in $\gE$, 
and  $h\in\sH$,  let 
\[
 \Integral_\sE h \dmu_\phi  \quad:=\quad  \int_{a_1}^{b_1} h(x) \  \mathrm{d} \phi(x) 
\ + \ 
\int_{a_2}^{b_2} h(x) \  \mathrm{d} \phi(x) 
 \ + \ \cdots \ + \ \int_{a_N}^{b_N} h(x) \  \mathrm{d} \phi(x) ,
\]  
where the integrals on the right-hand side  are   Riemann-Stieltjes  integrals.}

    We will show that there is a unique integrator, defined for all contents on $\gB$ and 
 a  large space of real-valued functions.  To explain this, we need some notation. 
 Let $\gB$  be a Boolean subalgebra of $\gR(\sS)$.  
We will say that a function $f:\sS\into\Real$ is {\dfn $\gB$-simple} if there is a finite $\gB$-partition $\gP$ of $\sS$ such that
$f$ is constant on each cell of $\gP$.\footnote{The behaviour of $f$ on the boundaries of the cells will be  irrelevant, so it is left arbitrary. }
  We then say that $f$ is {\dfn subordinate} to $\gP$.     We will say that a bounded function $h:\sS\into\Real$ is {\dfn $\gB$-integrable} if there exists a sequence $\{f_n\}_{n=1}^\oo$ of $\gB$-simple functions such that $\lim_{n\goto\oo}\norm{h-f_n}{\oo}=0$.  Let $\sL_\gB(\sS)$ be the set of all $\gB$-integrable functions.   We will soon see that $\sL_\gB(\sS)$ is a vector space, and closed in the uniform norm.

  We will be particularly interested in a certain class of $\gB$-integrable functions. 
 Let $\sC(\sS,\Real)$ denote the vector space of all continuous, real-valued functions on
 $\sS$. Let  $\Cb(\sS,\Real)$ be the Banach space of {\em bounded}, continuous, real-valued functions, with the uniform norm $\norm{\cdot}{\oo}$.   
 Let $\gE$ be the Boolean algebra of elementary regular open subsets of $\Real$, as defined in Example \ref{X:algebra}(a). 
 Let $\gB\subseteq\gR(\sS)$ be a Boolean subalgebra of $\gR(\sS)$.
  A function $f:\sS\rightarrow\Real$ is   {\dfn $\gB$-comeasurable} 
  if  $\Int\lb(f^{-1}\lb[\Cl(\sE)\rb]\rb)\in\gB$ for all $\sE\in\gE$.  
 Equivalently,  $f$ is {\dfn $\gB$-comeasurable} if $\Int\lb(f^{-1}(-\oo,r]\rb)\in\gB$ and 
 $\Int\lb(f^{-1}[r,\oo)\rb)\in\gB$ for all $r\in\Real$.    We will soon see that any $\gB$-comeasurable function is $\gB$-integrable.
 But first, here are some examples.

\example{\label{X:comeasurable} (a)  If $f:\sS\into\Real$ is  continuous, then $f$ is $\gR(\sS)$-comeasurable.

\item (b)   Let $\sS$ be a differentiable manifold, and let
$\gB_{\mathrm{smth}}$ be the Boolean algebra  of regular open sets with piecewise smooth boundaries, from Example \ref{X:algebra}(b).
 If $f:\sS\into\Real$ is any differentiable
  function such that d$f(s)\neq 0$ for all $s\in\sS$.  Then for all $r\in\Real$, the set $f^{-1}\{r\}$ is a smooth hypersurface, so that
   $\Int\lb(f^{-1}(-\oo,r]\rb)$ and 
 $\Int\lb(f^{-1}[r,\oo)\rb)$ have smooth boundaries; thus,   $f$ is $\gB_{\mathrm{smth}}$-comeasurable.
  
  More generally, let $f:\sS\into\Real$ be a {\em Morse function} ---that is, a smooth function where all critical points are nondegenerate.\footnote{A point $s\in\sS$ is {\dfn critical} if d$f(s)=0$.  It is {\dfn nondegenerate} if there exists a coordinate chart in a neighbourhood of $s$ such that the Hessian matrix of $f$ at $s$ (with respect to this chart) is nonsingular.}  Then for all $r\in\Real$, if $f^{-1}\{r\}$ does not contain a critical point, then it is a smooth hypersurface, whereas if  $f^{-1}\{r\}$ {\em does} contain one or more critical points, then it is a finite union of smooth hypersurfaces which meet only at these critical points (by the Morse Lemma). Either way,   $\Int\lb(f^{-1}(-\oo,r]\rb)$ and 
 $\Int\lb(f^{-1}[r,\oo)\rb)$ have piecewise smooth boundaries; thus,   $f$ is $\gB_{\mathrm{smth}}$-comeasurable.

  Not every differentiable function is $\gB_{\mathrm{smth}}$-comeasurable.  To see why some condition like d$f\neq0$ is required, suppose $\sS=\Real^2$, and define $f:\Real^2\into\Real$ by $f(x,y):=y^4\,\sin(x/y)$ for all $(x,y)\in\Real^2$.
  Then $f$ is differentiable everywhere on $\Real^2$ (with d$f(x,0)=0$ for all $x\in\Real$).
  However, if $\sR:= \Int\lb(f^{-1}(-\oo,0]\rb)$, then $\sR\not\in\gB_{\mathrm{smth}}$, because
  $\partial\sR=f^{-1}\{0\} = \{(x,y)\in\Real^2$; \ $y=0$ or $y=x/n\pi$ for some $n\in\Zahl\}$.  This set is
  an infinite union of lines passing through the origin, which ``converge'' to the horizontal line $y=0$; hence it cannot be represented as  a finite union of smooth curves.

\item (c)  Let $\sS$ be a topological vector space, and let $\gB_{\mathrm{poly}}$ be the Boolean algebra  of regular open polyhedra, from Example \ref{X:algebra}(c).
      A function $f:\sS\into\Real$ is {\dfn affine} if
  $f = f_0+r$ for some continuous linear function $f_0:\sS\into\Real$ and some constant $r\in\Real$.
  We say $f$ is  {\dfn piecewise affine}
  if there is a  collection $\sP_1,\ldots,\sP_N$ of disjoint regular open polyhedra such
  that $\sS=\sP_1\vee\cdots\vee\sP_N$, and a collection
  $f^1,\ldots,f^N:\sS\into\Real$ of affine functions, such that 
  $f_{\restr\sR_n}=f^n_{\restr\sR_n}$ for all $n\in[1\ldots N]$.  Any piecewise affine function
  is $\gB_{\mathrm{poly}}$-comeasurable. 

\item  (d) The sum of two  $\gB$-comeasurable functions is not necessarily $\gB$-comeasurable.
To see this, let $\sS=\Real$, and let $\gE$ be the Boolean algebra  of elementary sets from Example \ref{X:algebra}(a).
Let $f(x):=-2x$ and let $g(x):=2x+ x^2\,\sin(1/x)$, for all $x\in\Real$.  Then both $f$ and $g$ are $\gE$-comeasurable (because they are continuous and monotone).  But if $h:=f+g$, then  $h(x)=x^2\,\sin(1/x)$ for all $x\in\Real$.
This function is {\em not} $\gE$-comeasurable: if $\sR:=\Int\lb(h^{-1}(-\oo,0]\rb)$, then $\sR$ is an infinite union of
open intervals, so $\sR\not\in\gE$.  (Indeed,  $\partial\sR=\{0\}\union \{1/n\pi$; \ $n\in\Zahl\}$, which has a cluster point at
$0$.)
 }

 Let $\CB(\sS)$ be the set of all $\gB$-comeasurable functions in $\Cb(\sS,\Real)$.  This set is not necessarily closed under addition,  as shown by Example \ref{X:comeasurable}(d).
   So, let $\GB(\sS)$ be the  closed  linear subspace of $\Cb(\sS,\Real)$ spanned by $\CB(\sS)$;
    then $\GB(\sS)$ is a Banach space under the uniform norm $\norm{\cdot}{\oo}$.
 (If $\gB=\gR(\sS)$, then  $\GB(\sS)=\CB(\sS)=\Cb(\sS,\Real)$, because every continuous function is
 $\gR(\sS)$-comeasurable.)  For any subset $\sB\subseteq\sS$, let $\GB(\sB):=\{g_{\restr\sB}$; \ $g\in\GB(\sS)\}$.  This is a linear subspace
    of $\Cb(\sB,\Real)$.\footnote{Even when $\gB=\gR(\sS)$, it is {\em not}  generally true that
    $\GB(\sB)=\Cb(\sB,\Real)$, for two reasons.  First, some functions in  $\Cb(\sB,\Real)$ cannot
     be extended to functions in $\Cb(\barsB,\Real)$, where $\barsB$ is the closure of
    $\sB$.  (For example, let $\sS:=\Real$, let $\sB:=(0,1)$, and let
      $g(x):=\sin(1/x)$ for all $x\in\sB$.)  Second, not all functions in $\Cb(\barsB,\Real)$ can be extended
      to $\Cb(\sS,\Real)$, unless $\sS$ is a normal space, which we will not assume in general.}   
    Likewise, let $\sL_\gB(\sB):=\{g_{\restr\sB}$; \ $g\in\sL_\gB(\sS)\}$.

\Proposition{\label{density.of.simple.functions2}}
{For any $\sB\in\gB$, $\sL_\gB(\sB)$ is a uniformly closed vector space, and
 $\GB(\sB)\subseteq \sL_\gB(\sB)$.}

\bthmprf[ Proof of Proposition \ref{density.of.simple.functions2}.]
 It suffices to consider the case $\sB=\sS$; the general case is similar.
 Clearly, $\sL_\gB(\sS)$ is closed under scalar multiplication.  
To see that $\sL_\gB(\sS)$ is a vector space,  suppose $g=g_1+g_2$ for some $g_1,g_2\in\sL_\gB(\sS)$.
Let $\eps>0$.
 There exist $f_1,f_2\in\sF_\gB(\sS)$ with $\norm{f_n-g_n}{\oo}<\eps/2$ for  $n\in\{1,2\}$. Let $f:=f_1+f_2$.   Then $f\in\sF_\gB(\sS)$ (it is subordinate to any $\gB$ partition that is a common refinement of the partitions for $f_1$ and $f_2$).  
  Clearly,   $\norm{f-g}{\oo}<\eps$.  This works for any $\eps>0$.
 Thus, $g\in\sL_\gB(\sS)$.

 To see that $\sL_\gB(\sS)$ is uniformly closed, let $g$ be in the uniform closure of $\sL_\gB(\sS)$, and let $\eps>0$.   There exists $ h\in\sL_\gB(\sS)$ with $\norm{g-h}{\oo}<\eps/2$, and then
 there exists $f\in\sF_\gB(\sS)$ with $\norm{f-h}{\oo}<\eps/2$.  
 Thus, $\norm{f-g}{\oo}<\eps$.   This works for any $\eps>0$. Thus, $g\in\sL_\gB(\sS)$. 
 
 To see that $\GB(\sS)\subseteq\sL_\gB(\sS)$, 
it suffices to show that $\CB(\sS)\subseteq\sL_\gB(\sS)$, because $\GB(\sS)$ is defined to be the smallest closed vector space containing $\CB(\sS)$, and we have just shown that 
$\sL_\gB(\sS)$ is a closed vector space.    So 
let $g\in\CB(\sS)$.  For any    $\eps>0$, we must show that there exists $f\in\sF_\gB(\sS)$ with   $\norm{f-g}{\oo}<\eps$.   (In fact, we will even construct $f$ such that $f(s)\leq g(s)$ for all $s\in\sS$.) 
   Since $g$ is bounded, there is some $M\in\Natur$ that $|g(s)|< M$ for all
$s\in\sS$.
Fix $N\in\Natur$ such that ${1/N}<\eps$.   For all $m\in[-MN\ldots MN)$, let $\sI_m:=[\frac{m}{N},\frac{m+1}{N}]$
(a closed interval in $\Real$).  Let $\sC_m:=g^{-1}(\sI_m)$ and let $\sB_m:=\Int(\sC_m)$. 
{Then $\sB_m\in\gB$, because $g$ is 
 $\gB$-comeasurable.}

\Claim{$\D\bigvee_{m=-MN}^{MN-1} \sB_m=\sS$.}
\bclaimprf
Let $\sB^*:=\D\Union_{m=-MN}^{MN-1} \sB_m$.
We must show that {$\sB^*$}  is dense in $\sS$.
Suppose not.  Let $s\in\sS$ be a point not in the closure of {$\sB^*$}.
Then there is some open neighbourhood $\sD$ of $s$ which does not intersect
{$\sB^*$}.  Now, for any $m$, the set $g^{-1}(\frac{m}{N},\frac{m+1}{N})$ is a
subset of $\sB_m$ (because it is an open subset of $\sC_m$,  because $g$ is continuous).  Thus, $\sD$ cannot intersect
$g^{-1}({\frac{m}{N},\frac{m+1}{N}})$.  Thus, for all $s'\in\sD$ we must have $g(s')={\frac{m}{N}}$ for some $m\in[-MN\ldots MN)$.  In particular, $g(s)={\frac{m_0}{N}}$ for some $m_0\in[-MN\ldots MN)$.
By making $\sD$ small enough, we can ensure that $g(s')={\frac{m_0}{N}}$ for all $s'\in\sD$ (because $g$ is continuous).    But then $\sD$ is an open subset of $\sC_{m_0}$; hence $\sD\subseteq\sB_{m_0}$.  Contradiction.
\eclaimprf

\noindent Let $\sP_{-MN}:=\sB_{-MN}$ and, for any $m\in (-MN\ldots MN)$, define $\sP_m:=\sB_m\intsct(\neg\sB_{m-1})$. Let $\sM:=\{m\in[-MN\ldots MN)$; \ $\sP_m\neq\emptyset\}$, and define  $\gP:=\{\sP_m\}_{m\in\sM}$; then $\gP$ is a  $\gB$-partition of $\sS$.
 We define $f\in\sF_\gB(\sS)$ as follows:
 for all $m\in{\sM}$, and all ${s}\in\sP_m$, define $f(s):=\frac{m}{N}$.  
 Meanwhile, for all $s\in\sS$ not in $\Union_{{m\in\sM}} \sP_m$, define $f(s):=g(s)$.
 Thus, $f\in\sF_\gB(\sS)$, and  $f(s)\leq g(s)$ for all $s\in\sS$.  Finally, for any $s\in\sS$,
 $|f(s)-g(s)|\leq\frac{1}{N}< \eps$.  Thus,   $\norm{f-g}{\oo}<\eps$, as desired.
\ethmprf

 \noindent   Here is the main result of this section.

\Theorem{\label{from.probability.to.expectation0}}
{
Let $\sS$ be any topological space,  let $\gB$ be any Boolean subalgebra of $\gR(\sS)$,  and let $\sM(\gB)$ be the set of all contents on $\gB$.    There exists a unique integrator $\integral$ on $\gB\x\sL_\gB(\sS)\x\sM(\gB)$.
For any $f\in\sF_\gB(\sS)$, if $f$ is subordinate to the $\gB$-partition $\gP=\{\sP_1,\ldots,\sP_N\}$, with $r_1,\ldots,r_N\in\Real$ such that $f(p)=r_n$ for all
 $p\in\sP_n$ and  all $n\in[1\ldots N]$, then for any  $\mu\in\sM(\gB)$ and   $\sB\in\gB$, 
 \beqn
 \label{integral.of.simple.function}
 \Integral_\sB f \dmu \quad=\quad \sum_{n=1}^N r_n \, \mu[\sP_n\intsct\sB].
 \eeqn
  Meanwhile, for any $g\in\sL_\gB(\sB)$ and   $\sB\in\gB$, 
   \beqn
\label{integral.defn}
\Integral_{\sB} g \dmu\quad = \quad \sup_{f\in\sF_g(\sB)} \  \Integral_\sB f\dmu,
\eeqn
where $\sF_g(\sB):=\{f\in\sF_\gB(\sB)$; \ $f(b)\leq g(b)$, for all $b\in\sB\}$.
}

The proof  of Theorem \ref{from.probability.to.expectation0} proceeds as follows.  Let $\gA$ be
the Boolean algebra of subsets of $\sS$ generated by all open subsets and all nowhere-dense subsets.
Then $\gR(\sS)\subset\gA$ (although the Boolean algebra operations are not the same).  Furthermore,
 $\mu$ can be extended from $\gB$ to a finitely additive measure $\nu$ on $\gA$ in a natural way.
 For any $g\in\sL_\gB(\sB)$ and   $\sB\in\gB$,
 we then define 
 \beqn
 \label{new.integral.defn}
 \Integral_{\sB} g \dmu\quad:=\quad \int_{\sB} g \dnu,
\eeqn
where the $\nu$-integral on the right-hand side defined in the standard way \cite[\S11.2]{AliprantisBorder}, 
\cite[Ch.4]{BhaskaraRaoBook}.  It is easily verified
that this construction satisfies  properties (I1) - (I3).  There remains the issue
of uniqueness.  If $\gB=\gR(\sS)$, then the extension $\nu$ is unique.  If
 $\gB\subset\gR(\sS)$, then  $\nu$ is generally not unique, but it is ``unique enough'' that
   equation (\ref{integral.of.simple.function})  holds for any $\gB$-simple function $f$; thus
 equation (\ref{integral.defn}) holds for any $g\in\sL_\gB(\sB)$.

\bthmprf[Proof of Theorem \ref{from.probability.to.expectation0}.]
Let $\gN$ be the set of all nowhere-dense subsets of $\sS$, 
let $\gO$ be the set of all open subsets of $\sS$, and
let $\gA:=\{\sO\union\sN$; \  $\sO\in\gO$ and $\sN\in\gN\}$.

\Claim{\label{from.probability.to.expectation0.C1}$\gA$ is a Boolean algebra.}
\bclaimprf
 (1) \ Let $\sA_1,\sA_2\in\gA$; then $\sA_1=\sO_1\union\sN_1$ and  $\sA_2=\sO_2\union\sN_2$  for some
 $\sO_1,\sO_2\in\gO$ and $\sN_1,\sN_2\in\gN$.  
 Let $\sO:=\sO_1\union\sO_2$ and $\sN:=\sN_1\union\sN_2$; then
 $\sO\in\gO$ and $\sN\in\gN$, and $\sA_1\union\sA_2 =\sO\union\sN$, so that $\sA_1\union\sA_2$ is
 also an element of $\gA$.  Thus, $\gA$ is closed under finite unions.
 
 \item (2) \  Let $\sN'_1:=(\sN_1\intsct\sO_2)$,  $\sN'_2:=(\sN_2\intsct\sO_1)$, and $\sN_0:=\sN_1\intsct\sN_2$.
 Then $\sN_0,\sN'_1,\sN'_2\in\gN$.  Thus, if $\sN:=\sN'_1\union\sN'_2\union \sN_0$, then $\sN\in\gN$.
  Let $\sO:=\sO_1\intsct\sO_2$; then $\sO\in\gO$.   But
\[
\sA_1\intsct\sA_2 \ = \ (\sO_1\union\sN_1)\intsct(\sO_2\union\sN_2)
\ = \ \sO\union \sN'_1\union\sN'_2\union\sN_0
\ = \ \sO\union\sN,
\]
so that $\sA_1\intsct\sA_2$ is
 also an element of $\gA$.  Thus, $\gA$ is closed under finite intersections.

\item (3) \  Let $\sC\subset\sS$ be any closed set. If $\sO:=\Int(\sC)$, and  $\sN:=\partial\sC$, then $\sO\in\gO$,
$\sN\in\gN$, and $\sC=\sO\union\sN$, so $\sC\in\gA$.  Thus, $\gA$ contains all closed sets.
In particular, $\gA$ contains the complement of any element of $\gO$.

\item (4) \  Let $\sN\in\gN$.  Let $\sM:=\Cl(\sN)$; then $\sM\in\gN$ also.
Let $\sN':=\sM\setminus\sN$; then $\sN'\in\gN$ also, since  $\sN'\subseteq\sM$.
Let $\sO:=\compl{\sM}$; then $\sO$ is open,
and $\compl{\sN}=\sO\union\sN'$.  Thus, $\compl{\sN}\in\gA$.

\item (5) \  Finally, let $\sA\in\gA$ be arbitrary; then $\sA=\sO\union\sN$ for some  $\sO\in\gO$ and $\sN\in\gN$.
Thus, $\compl{\sA}=\compl{\sO}\intsct\compl{\sN}$.  But $\compl{\sO}\in\gA$ and $\compl{\sN}\in\gA$ by steps 
(3) and (4).  Thus  $\compl{\sO}\intsct\compl{\sN}\in\gA$, by step (2).  Thus, $\gA$ is closed under complementation.
\eclaimprf

\Claim{\label{from.probability.to.expectation0.C2}
$\gA=\{(\sR\setminus\sN^-)\disj\sN^+$; \  $\sR\in\gR(\sS)$ and $\sN^-,\sN^+\in\gN$ are disjoint$\}$. 
}
\bclaimprf
Let $\gA':=\{(\sR\setminus\sN^-)\disj\sN^+$; \  $\sR\in\gR(\sS)$ and $\sN^-,\sN^+\in\gN$ are disjoint$\}$.  
Clearly, $\gA'\subseteq\gA$, because $\gR(\sS)\subseteq\gO$.  We must show the opposite inclusion.
So, let $\sA\in\gA$; then $\sA=\sO\union\sN$ for some  $\sO\in\gO$ and $\sN\in\gN$.
Let $\sR:=\Int[\Cl(\sO)]$; then $\sR\in\gR(\sS)$ and $\sO\subseteq\sR$.
Let $\sM:=\sR\setminus\sO$; then $\sM\subseteq\Cl(\sO)\setminus\sO=\partial\sO\in\gN$, and thus, $\sM\in\gN$.
Let   $\sN':=\sN\intsct\sM$ and   $\sN^+:=\sN\setminus\sR$; then
$\sO\union\sN=\sO\disj\sN'\disj\sN^+$.
Finally, let
$\sN^-:=\sM\setminus\sN'$.  Then $\sN^-$ is disjoint from $\sN^+$ (because $\sN^-\subseteq\sR$ and
$\sR$ is disjoint from $\sN^+$), and $\sO\disj\sN' = (\sR\setminus\sM)\disj\sN' = \sR\setminus\sN^-$.
Thus,
$\sA =\sO\union\sN
=\sO \disj \sN'\disj\sN^+
= (\sR\setminus\sN^-) \disj\sN^+$,
as desired.
\eclaimprf

\newcommand{\sW}{\mathcal{W}}

\Claim{ If $\sO\in\gO$ and $\sN\in\gN$,  then $\Cl(\sO\setminus\sN)=\Cl(\sO)$.
\label{from.probability.to.expectation0.C3a}}
\bclaimprf
Clearly, 
$\Cl(\sO\setminus\sN)\subseteq \Cl(\sO)$.   To show the reverse inclusion,
 let $x\in\Cl(\sO)$.  Then for every open neighbourhood $\sU$ around $x$, we have $\sO\intsct\sU\neq\emptyset$.
But $\Int(\sN)=\emptyset$ because  $\sN\in\gN$.
Thus, $\sO\intsct\sU\not\subseteq\sN$.  Thus, $\sU\intsct(\sO\setminus\sN)\neq\emptyset$.
Thus,  every open neighbourhood of $x$ intersects $\sO\setminus\sN$, so
$x\in\Cl(\sO\setminus\sN)$, as desired.
\eclaimprf

\Claim{If $\sF\subseteq\sS$ is closed and $\sN\in\gN$,  then $\Int(\sF\union\sN)=\Int(\sF)$.
\label{from.probability.to.expectation0.C3b}}
\bclaimprf
 Let $\sO:=\compl{\sF}$.  Then $\sO\in\gO$ and $\Cl(\sO)=\compl{\Int(\sF)}$.  Also, $\sO\setminus\sN= \compl{(\sF\union\sN)}$; thus $\Cl(\sO\setminus\sN)=\compl{[\Int(\sF\union\sN)]}$.
 Thus, Claim \ref{from.probability.to.expectation0.C3a} implies that
$\compl{[\Int(\sF\union\sN)]}=\compl{\Int(\sF)}$.  Thus, $\Int(\sF\union\sN)=\Int(\sF)$.
\eclaimprf

\Claim{\label{from.probability.to.expectation0.C3}
For any $\sA\in\gA$, there exists a unique $\sR\in\gR(\sS)$ and disjoint $\sN^-,\sN^+\in\gN$ such that
$\sN^-\subseteq\sR$ and  $\sA=(\sR\setminus\sN^-)\disj\sN^+$.  Furthermore, 
{\rm(a)} 
$\sR= \Int\lb[\Cl(\sA)\rb]$; \ {\rm(b)} 
$\sN^+=\sA\setminus\sR$; \ and {\rm(c)} 
$ \sN^-=\sR\setminus\sA$.
}
\bclaimprf
By Claim \ref{from.probability.to.expectation0.C2}, we know that $\sA=(\sR\setminus\sN^-)\disj\sN^+$,
for some  $\sR\in\gR(\sS)$ and disjoint $\sN^-,\sN^+\in\gN$.  By replacing
$\sN^-$ with $\sN^-\intsct\sR$ if necessary, we can assume without loss of generality that $\sN^-\subseteq\sR$.
Statements (b) and (c) follow immediately.
It remains to prove (a).     First, note that
\beqn
\label{from.probability.to.expectation0.C3.e2}
\Cl(\sA) \ = \ \Cl\lb((\sR\setminus\sN^-)\disj\sN^+\rb)
\ = \ 
 \Cl\lb(\sR\setminus\sN^-\rb)\union\Cl\lb(\sN^+\rb)
\ \eeequals{(*)} \ 
 \Cl(\sR)\union\Cl\lb(\sN^+\rb),
 \eeqn
 where $(*)$ is by Claim \ref{from.probability.to.expectation0.C3a}.
 Thus, 
\[
 \Int\lb[\Cl(\sA)\rb]
\quad\eeequals{(*)}\quad
 \Int\lb[\Cl(\sR)\union\Cl\lb(\sN^+\rb) \rb]
 \quad\eeequals{(\dagger)}\quad
 \Int\lb[\Cl(\sR)\rb]
 \quad\eeequals{(\diamond)}\quad
 \sR,
 \]
 which proves (a).
Here, $(*)$ is by  (\ref{from.probability.to.expectation0.C3.e2}),
$(\dagger)$ is by  Claim \ref{from.probability.to.expectation0.C3b},
and $(\diamond)$ is because $\sR\in\gR(\sS)$.

Now suppose that also 
$\sA=(\sR_1\setminus\sN_1^-)\disj\sN_1^+$
for some other $\sR_1\in\gR(\sS)$ and disjoint $\sN_1^-,\sN_1^+\in\gN$
with $\sN^-_1\subseteq\sR_1$.  By repeating the above arguments, we deduce that 
(a1) $\sR_1=\Int\lb[\Cl(\sA)\rb]$,
(b1) $\sN_1^+=\sA\setminus\sR_1$,
and (c1) $\sN_1^-=\sR_1\setminus\sA$.
  
Combining (a) and (a1) yields $\sR_1=\sR$.  Substituting this into (b1) and (c1) and
combining with (b) and (c), we get
$\sN_1^+=\sN^+$ and $\sN_1^-=\sN^-$.
\eclaimprf

\Claim{\label{from.probability.to.expectation0.C4}Let $\sA_1,\sA_2\in\gA$, and let $\sA=\sA_1\union\sA_2$.
Write
$\sA_1=(\sR_1\setminus\sN_1^-)\disj\sN_1^+$, 
$\sA_2=(\sR_2\setminus\sN_2^-)\disj\sN_2^+$,
and $\sA=(\sR\setminus\sN^-)\disj\sN^+$,
 as in Claim \ref{from.probability.to.expectation0.C3}. Then 
 \bthmlist
 \item $\sR=\sR_1\vee\sR_2$.
 \item If $\sA_1$ and $\sA_2$ are disjoint, then $\sR_1$ and $\sR_2$ are disjoint.
 \ethmlist
 }
 \bclaimprf
 (a)
 Let $\sN'_1:=\sN_1^-\setminus\sR_2$ and $\sN'_2:=\sN_2^-\setminus\sR_1$;
then $\sN'_1\subseteq\sR_1$, \ $\sN'_2\subseteq\sR_2$, and
 $\sN'_1,\sN'_2\in\gN$.  Furthermore,
\beq
\sA
&=&
\sA_1\union\sA_2 \quad=\quad
(\sR_1\setminus\sN_1^-)\union\sN_1^+\union (\sR_2\setminus\sN_2^-)\union\sN_2^+
\\ &=&
\lb((\sR_1\union\sR_2) \setminus (\sN'_1\union\sN'_2) \rb)\union(\sN_1^+\union\sN_2^+).\\
\mbox{Thus,}\quad
\Cl[\sA]
& = & 
\Cl\lb[\lb((\sR_1\union\sR_2) \setminus (\sN'_1\union\sN'_2)\rb) \union(\sN_1^+\union\sN_2^+)\rb]
\\ &=&
\Cl\lb[(\sR_1\union\sR_2) \setminus (\sN'_1\union\sN'_2)\rb] \union\Cl(\sN_1^+\union\sN_2^+)
\\ &\eeequals{(*)}&
\Cl(\sR_1\union\sR_2) \union \Cl(\sN_1^+\union\sN_2^+),
\eeq
where $(*)$ is by Claim \ref{from.probability.to.expectation0.C3a}.
  Thus,
\beqn
\Int\lb[\Cl(\sA)\rb]
  \  =  \  
\Int\lb[\Cl(\sR_1\union\sR_2) \union \Cl(\sN_1^+\union\sN_2^+)\rb]
 \ \eeequals{(*)} \ 
\Int\lb[\Cl(\sR_1\union\sR_2)\rb], 
\label{from.probability.to.expectation0.C4.e1}
\eeqn
where $(*)$ is  by Claim \ref{from.probability.to.expectation0.C3b}.
 Thus,
\[
\sR \quad\eeequals{(*)}\quad \Int\lb[\Cl(\sA)\rb]
\quad\eeequals{(\dagger)}\quad 
\Int\lb[\Cl(\sR_1\union\sR_2)\rb] 
\quad=\quad
\sR_1\vee\sR_2,
\]
as claimed.  Here,
 $(*)$ is by Claim \ref{from.probability.to.expectation0.C3}(a), 
and $(\dagger)$ is by equation (\ref{from.probability.to.expectation0.C4.e1}).

\item (b) First note that
\beq
\emptyset
 &=&
 \sA_1\intsct\sA_2 
 \ = \ \lb((\sR_1\setminus\sN_1^-)\disj\sN_1^+\rb) \intsct \lb((\sR_2\setminus\sN_2^-)\disj\sN_2^+\rb)
 \ \supseteq \ 
 (\sR_1\setminus\sN_1^-) \intsct (\sR_2\setminus\sN_2^-)
\\ &=&
(\sR_1\intsct\sR_2) \setminus(\sN_1^-\union\sN_2^-).
\eeq
However, $\sR_1\intsct\sR_2$ is open, while $\sN_1^-\union\sN_2^-$ has empty interior (because it is in $\gN$).
Thus, if $(\sR_1\intsct\sR_2) \setminus(\sN_1^-\union\sN_2^-)=\emptyset$, then $\sR_1\intsct\sR_2 =\emptyset$. \eclaimprf

\noindent 
Now, let $\mu$ be a content on some Boolean algebra $\gB\subseteq\gR(\sS)$.
Let $\mu'$ be an extension of $\mu$ to a content defined on the Boolean algebra $\gR(\sS)$;  this exists by the
Horn-Tarski Extension Theorem \cite[Theorem 1.22]{HornTarski48}.
We can define a function $\nu:\gA\into\Real_+$ as follows: for any $\sA\in\gR$, if $\sA=(\sR\setminus\sN^-)\disj\sN^+$
as in Claim \ref{from.probability.to.expectation0.C3}, then define $\nu(\sA):=\mu'(\sR)$.  This is well-defined by the
uniqueness of the representation in Claim \ref{from.probability.to.expectation0.C3}.
In particular, $\nu(\sR) =\mu'(\sR)$ for all $\sR\in\gR(\sS)$, and $\nu(\sN)=0$ for all $\sN\in\gN$.

\Claim{\label{from.probability.to.expectation0.C5} $\nu$ is a finitely additive measure on $\gA$.}
\bclaimprf
Let $\sA_1,\sA_2\in\gA$ be disjoint, and let $\sA=\sA_1\disj\sA_2$.
 We must show that $\nu[\sA]=\nu[\sA_1]+\nu[\sA_2]$.
 Write
 $\sA_1=(\sR_1\setminus\sN_1^-)\disj\sN_1^+$, 
$\sA_2=(\sR_2\setminus\sN_2^-)\disj\sN_2^+$,
and $\sA=(\sR\setminus\sN^-)\disj\sN^+$,
 as in Claim \ref{from.probability.to.expectation0.C3}.
Then
\[
\nu[\sA] \quad \eeequals{(*)} \quad  \mu'[\sR] 
\quad \eeequals{(\dagger)} \quad \mu'[\sR_1\vee\sR_2]
\quad \eeequals{(\diamond)} \quad  \mu'[\sR_1]+\mu'[\sR_2]
\quad \eeequals{(*)} \quad  \nu[\sA_1]+\nu[\sA_2],
\]
as desired.
Here, both $(*)$ are by the definition of $\nu$, $(\dagger)$ is because  $\sR=\sR_1\vee\sR_2$ by
 Claim \ref{from.probability.to.expectation0.C4}(a), and $(\diamond)$ is because $\sR_1$ and $\sR_2$ are disjoint
 by Claim \ref{from.probability.to.expectation0.C4}(b).
\eclaimprf

\noindent Now, for any $g\in\sL_\gB(\sB)$ and   $\sB\in\gB$,
  define $\integral_{\sB} g \dmu$ by equation (\ref{new.integral.defn}).
  It is easily verified that this definition satisfies properties (I1) - (I3) in the definition of ``integrator''.

\Claim{\label{from.probability.to.expectation0.C6}Equation {\rm(\ref{integral.of.simple.function})} holds for every $\gB$-simple function $f$.}
\bclaimprf
Let $f$ be $\gB$-simple function,  subordinate to a $\gB$-partition $\gP=\{\sP_1,\sP_2,\ldots,\sP_N\}$
as described prior to equation (\ref{integral.of.simple.function}).
Note that $\partial\sP_1,\ldots,\partial\sP_N\in\gN$, so 
$\nu[\partial\sP_1]=\cdots=\nu[\partial\sP_N]=0$.
Thus, we can modify the values of $f$ on $\partial\sP_1\union\ldots\union\partial\sP_N$ without affecting
the value of $\int f\dnu$. In particular, we can assume without loss of generality that
$f(s)=0$ for all $s\in \partial\sP_1\union\ldots\union\partial\sP_N$.  Thus,
for any $\sB\in\gB$, we have
\[
\Integral_\sB f\dmu  \quad\eeequals{(*)}\quad
\int_\sB f\dnu \quad\eeequals{(\dagger)}\quad \sum_{n=1}^N r_n \, \nu[\sP_n\intsct\sB]
\quad\eeequals{(\diamond)}\quad  \sum_{n=1}^N r_n \, \mu[\sP_n\intsct\sB],
\]
as claimed.  Here, $(*)$ is by defining formula (\ref{new.integral.defn}), $(\dagger)$ is because 
$f$ is zero on $\partial\sP_1\union\ldots\union\partial\sP_N$,
and $(\diamond)$ is because $\nu$ agrees with  $\mu'$ on $\gR(\sS)$, while $\mu'$ agrees with
$\mu$ on $\gB$.
\eclaimprf

\Claim{\label{from.probability.to.expectation0.C7}
Equation {\rm(\ref{integral.defn})} holds for all $g\in\sL_\gB(\sS)$.}
\bclaimprf
  Let  $\{f_n\}_{n=1}^\oo$ be a sequence of $\gB$-simple functions that converge 
  uniformly to $g$ (such a sequence exists by the definition of $\sL_\gB(\sS)$).
Then for any $\sB\in\gB$,
\beqn
\label{from.probability.to.expectation0.e1}
\int_\sB g\dmu \quad=\quad\lim_{n\goto\oo} \int_\sB f_n\dmu,
\eeqn
because  the function $ h\mapsto\integral_\sB h\dmu$ is continuous with respect to the uniform norm, by (I1). 
In particular, we can construct the sequence $\{f_n\}_{n=1}^\oo$ so that it converges to $g$ from below;
then (\ref{from.probability.to.expectation0.e1}) yields equation (\ref{integral.defn}).
\eclaimprf

\noindent {\em Well-defined.}  The extension $\mu'$ is not necessarily unique;  we must show that the above definition is independent
of the choice of $\mu'$.   Let $\mu'_1$ and $\mu'_2$ be two different contents on $\gR(\sS)$
extending $\mu$, and let $\nu_1$ and $\nu_2$ be the resulting finitely additive measures on $\gA$ obtained via
Claim \ref{from.probability.to.expectation0.C5}.   Claim \ref{from.probability.to.expectation0.C6} implies that
$\nu_1$ and $\nu_2$ yield the same integral for any $\gB$-simple function.
Then Claim \ref{from.probability.to.expectation0.C7} implies that
$\nu_1$ and $\nu_2$ yield the same integral for any element of $\sL_\gB(\sS)$.
\ethmprf

\noindent The next result says that the integration operator is ``strictly increasing''  when restricted to $\CB(\sS)$, as long
as $\mu$ has full support.
For any $\sB\in\gB$, we define $\sH(\sB):=\{h_{\restr\sB}$; \ $h\in\sH\}$.

\Proposition{\label{strict.monotonicity}}
{
Suppose  $\mu$ has full support. 
Let $\sB\in\gB$   be nonempty, and let $f\in\CB(\sB)$.
If $f(b)> 0$ for all $b\in\sB$, then
$\integral_\sB f\dmu> 0$. 
If $f(b)< 0$ for all $b\in\sB$, then
$\integral_\sB f\dmu< 0$. 
}
\bthmprf By linearity, it is sufficient to consider the case where $f(b)> 0$ for all $b\in\sB$. 
 For all $n\in\Natur$, let $\sC_n:=f^{-1}[0,\frac{1}{n}]$. Note that $\Intsct_{n=1}^\oo(\sB\intsct \sC_n)=\sB\intsct f^{-1}\{0\}=\emptyset$ (because $0< f(b)$ for all $b\in\sB$).
There must then be some $n\in\Natur$ such that $\sB\intsct\sC_n\neq\sB$.  Let $\sQ:=\sB\intsct\Int(\sC_n)$; then 
 $\sQ\subsetneq\sB$ and $\sQ\in\gB$ (because $f$ is $\gB$-comeasurable).  Thus, if  $\sP:=\sB\intsct(\neg \sQ)$, then $\emptyset\neq \sP\subseteq\sB$ and $\sP\in\gB$.

\Claim{ $f(p)\geq \frac{1}{n}$ for all $p\in\sP$.}
\bclaimprf
(by contradiction) Let $p\in\sP$, and suppose $f(p)<\frac{1}{n}$.  Then there is some open neighbourhood
$\sO\subseteq\sS$ containing $p$ such that $f(o)<\frac{1}{n}$ for all $o\in\sO$ (because $f$ is continuous).
Thus, $\sO\subseteq\Int(\sC_n)$.  Thus, since $p\in\sP\subseteq\sB$ and $p\in\sO$, we deduce that $p\in\sB\intsct\Int(\sC_n)=\sQ$.
But $p\in\neg\sQ$.  Contradiction.
\eclaimprf
Let $\kap:\sB\into\Real$ be the constant function with value $\frac{1}{n}$. Then 
\beqn
\label{strict.monotonicity.e1}
 \Integral_\sP f\dmu\quad\geeeq{(*)}\quad\Integral_\sP \kap\dmu\quad=\quad  \frac{1}{n}\cdot\mu[\sP], 
 \eeqn
 where $(*)$ is by Claim 1.
      Meanwhile, $f(q)>0$ for all $q\in\sQ$.  Thus, a similar argument yields
  \beqn
\label{strict.monotonicity.e2}
\Integral_\sQ f\dmu\quad\geq\quad 0.
\eeqn
 Thus,
\beq
\Integral_\sB f\dmu &\eeequals{(*)}& \Integral_\sQ f\dmu + \Integral_\sP f\dmu
\quad \geeeq{(\dagger)} \quad\frac{1}{n}\cdot \mu[\sP]  
\quad \grt{(\diamond)} \quad 
0,
\eeq
Here $(*)$ is by property (I2) (because
${\sB}=\sP\vee\sQ$ and $f\in\CB(\sB)$),
while $(\dagger)$ is by inequalities (\ref{strict.monotonicity.e1}) and (\ref{strict.monotonicity.e2}).
Finally, $(\diamond)$ is because $\sP\neq\emptyset$ and $\mu$ has full support.
\ethmprf

\example{\label{X:ultrafilter.example} 
The hypothesis of {full support} is needed for Proposition \ref{strict.monotonicity}.  To see this,
let $\sS=\Natur$ with the discrete topology; then {every} subset of $\sS$ is regular, and
 $\vee$ is just the  union operation.  Let $\gB=\wp(\Natur)$ (the power set of $\Natur$), and  let $\gU$ be a free ultrafilter in $\wp(\Natur)$ (see Example \ref{X:ultrafilter.example2}).
In this case, conditions (c) and (d) take the form:
(c) No finite subset of $\Natur$ is in $\gU$; and
(d)  For any $\sR\subseteq\Natur$, either $\sR\in\gU$, or $(\Natur\setminus\sR)\in\gU$, but not both.

 Define $\del_\gU$ as in Example \ref{X:ultrafilter.example2}.
Then $\del_\gU$ is a finitely additive probability measure  on the power set of $\Natur$, and thus, a content.
For any $N\in\Natur$, property (c) implies that $[1\ldots N]\not\in\gU$;  thus,
property (d) implies that $(N\ldots \oo)\in\gU$.  Thus, $\del_\gU[1\ldots N]=0$ and
$\del_\gU(N\ldots \oo)=1$. 

Now let $f(n)=\frac{1}{n}$ for all $n\in\Natur$.  Then $f\in\Cb(\Natur,\Real)$,  and $f$ is $\gB$-comeasurable.
For any $N\in\Natur$, if
$\sR:=[1\ldots N]$ and $\sQ:=(N\ldots \oo)$, then
$\integral_\sR f\dmu =0$ (because $\del_\gU[\sR]=0$) while $\integral_\sQ[f]\dmu\leq 1/N$, because $f(q)\leq 1/N$ for all
$q\in\sQ$).  Thus,
 \[
 0\quad\leq\quad 
\Integral_\Natur f\dmu \quad = \quad \Integral_\sR f\dmu + \Integral_\sQ f\dmu
\quad \leq \quad 0 + \frac{1}{N}\quad=\quad \frac{1}{N}.
\]
Letting $N\goto\oo$, we obtain $\integral_\Natur f\dmu  \ = \ 0$, despite the fact that $f(n)>0$ for all $n\in\Natur$.}

Metaphorically speaking, the content in Example \ref{X:ultrafilter.example} is like a ``point mass at infinity''.  Later, we will  make this metaphor precise in Theorem \ref{content.vs.Riesz.charge4}.

\section{Measurable functions and change of variables\label{S:measurable}}
\setcounter{equation}{0}

In classical measure theory, a measurable function from a space $\sX$ to a space $\sY$
 can be used to ``push forward'' a measure from $\sX$ to $\sY$.  Conversely, it can be used to ``pull back''
 a  measurable  real-valued function from $\sY$ to $\sX$, and thereby convert an integral computation on $\sY$
into an integral computation on $\sX$,  via  ``change of variables''.  We will now develop an analogous theory for contents and their associated integrators.

 Let $\sX$ and $\sY$ be two topological spaces.
Let $\gA\subseteq\gR(\sX)$ and $\gB\subseteq\gR(\sY)$ be Boolean subalgebras of the algebras of regular
sets on $\sX$ and $\sY$.  A function $\phi:\sX\into\sY$ is {\dfn measurable} with respect to $\gA$ and $\gB$ if
$\phi^{-1}(\sB)\in\gA$ for all $\sB\in\gB$.  For example, if $\phi$ is a continuous, open function, then
$\phi$ is measurable with respect to $\gR(\sX)$ and $\gR(\sY)$ \cite[Appendix 4A2B, item (f)(iii), p.453]{Fremlin4II}.  
Unfortunately, not all continuous functions are
measurable with respect to the algebras of regular sets.\footnote{For example, let
$\sX=\sY=\Real$ with the standard topology, let $\phi(x)=x^2$, and let $\sB:=(0,\oo)$.  Then $\sB\in\gR(\Real)$,
but $\phi^{-1}(\sB)=(-\oo,0)\disj(0,\oo)\not\in\gR(\Real)$.}  Thus, we must introduce a weaker notion.

 For any function $\phi:\sX\into\sY$ and subset $\sB\subseteq\sY$, we define $\phi^{\leftarrow}(\sB):=\Int\lb(\phi^{-1}\lb[\Cl(\sB)\rb]\rb)$.  We say that $\phi$
  is {\dfn comeasurable} with respect to $\gA$ and $\gB$ if
$\phi^{\leftarrow}(\sB)\in\gA$ for all $\sB\in\gB$.     In particular,
if $\sY=\Real$ and $\gB=\gE$ (the algebra of elementary sets), then this is the
definition of ``comeasurable'' given for real-valued functions in Section \ref{S:notation}.

 \Lemma{\label{measurability.lemma}}
 {
 Let $\sX$ and $\sY$ be  topological spaces, and let $\gA\subseteq\gR(\sX)$ and   $\gB\subseteq\gR(\sY)$
  be Boolean subalgebras of the algebras of regular sets.
\bthmlist
\item Any $(\gA,\gB)$-measurable function from $\sX$ to $\sY$ is $(\gA,\gB)$-comeasurable.
\item Suppose $\phi:\sX\into\sY$ is both continuous and open.  Then for any  $\sB\in\gR(\sY)$ , we have
$\Int\lb(\phi^{-1}[\Cl(\sB)]\rb)=\phi^{-1}(\sB)$.  Thus, $\phi$ is $(\gA,\gB)$-measurable if and only if
$\phi$ is $(\gA,\gB)$-comeasurable.  Furthermore, if $\phi$ is $(\gA,\gB)$-(co)measurable, then the function
$\phi^{-1}:\gB\into\gA$ is a Boolean algebra homomorphism.
\ethmlist
}
\bthmprf
(a) Suppose $\phi:\sX\into\sY$ is measurable.
Let $\sB\in\gB$.  Let $\sC:=\neg\sB$;  then $\sC\in\gB$ also.
Thus, if we define $\sD:=\phi^{-1}(\sC)$, then $\sD\in\gA$, because $\phi$ is measurable by hypothesis.
However, $\Cl(\sB)=\compl{\sC}$.  Thus,
\beq
\phi^{-1}\lb[\Cl(\sB)\rb]&=&\phi^{-1}(\compl{\sC})\quad=\quad\compl{\phi^{-1}(\sC)}\quad=\quad\compl{\sD}.\\
\mbox{Thus,}\quad \Int\lb(\phi^{-1}\lb[\Cl(\sB)\rb]\rb) &=& \Int\lb(\compl{\sD}\rb) \quad=\quad \neg\sD,
\eeq
which is an element of $\gA$, as desired.

\item (b) \ Suppose $\phi$ is open and continuous. Then $\phi^{-1}[\Int(\sB)]=\Int[\phi^{-1}(\sB)]$ and 
$\phi^{-1}[\Cl(\sB)]=\Cl[\phi^{-1}(\sB)]$, for any $\sB\subseteq\sY$.  Thus, if  $\sB\in\gR(\sY)$ , then
\beqn
\label{measurability.lemma.e0}
 \phi^{\leftarrow}(\sB)
 \ \ = \ \  
\Int\lb(\phi^{-1}[\Cl(\sB)]\rb)
 \ \ = \ \ 
\phi^{-1}\lb(\Int[\Cl(\sB)]\rb)
 \ \ = \ \ 
\phi^{-1}(\sB), 
\eeqn
as claimed, where the last step is because $\Int[\Cl(\sB)]=\sB$ because $\sB$    is a regular open set.  
Given the identity (\ref{measurability.lemma.e0}), the conditions for 
$\phi$ to be $(\gA,\gB)$-measurable and to be $(\gA,\gB)$-comeasurable are logically equivalent.

Now suppose that  $\phi$ is $(\gA,\gB)$-(co)measurable.  To see that 
$\phi^{-1}:\gB\into\gA$ is a Boolean algebra homomorphism, first recall that
\beqn
\label{measurability.lemma.e1}
\phi^{-1}(\sB_1\intsct\sB_2)\quad=\quad \phi^{-1}(\sB_1)\intsct \phi^{-1}(\sB_2),
\qquad\mbox{for any $\sB_1,\sB_2\subseteq\sY$.}
\eeqn
 Meanwhile, if $\sB\in\gB$, then
\begin{eqnarray}
\nonumber
\phi^{-1}\lb(\neg\sB\rb)
&=&
\phi^{-1}\lb[\Int\lb(\compl{\sB}\rb)\rb]
\quad=\quad
\Int\lb(\phi^{-1}\lb[\compl{\sB}\rb]\rb)
\\&=&
\Int\lb(\compl{\lb[\phi^{-1}(\sB)\rb]}\rb)
\quad=\quad
\neg \phi^{-1}(\sB).
\label{measurability.lemma.e2}
\end{eqnarray}
Finally, for any $\sB_1,\sB_2\in\gB$, de Morgan's law yields
$\sB_1\vee\sB_2 = \neg[(\neg\sB_1)\intsct(\neg\sB_2)]$.  Thus,
an  application of equations (\ref{measurability.lemma.e1}) and (\ref{measurability.lemma.e2}) yields
$\phi^{-1}(\sB_1\vee\sB_2)=\phi^{-1}(\sB_1)\vee \phi^{-1}(\sB_2)$.
\ethmprf

\noindent To illustrate Lemma \ref{measurability.lemma}, let
 $\phi:\sX\into\sY$ be {\em any} continuous function.  Then it is easy to see that $\phi$ is $(\gR(\sX),\gR(\sY))$-comeasurable.  Thus, if $\phi$ is also open, then Lemma \ref{measurability.lemma}(b) says that 
$\phi$ is  $(\gR(\sX),\gR(\sY))$-measurable, and 
$\phi^{-1}:\gR(\sY)\into\gR(\sX)$ is a Boolean algebra homomorphism.  The next example is more involved.

\example{Let $\sX$ and $\sY$ be differentiable manifolds, and let $\gA\subset\gR(\sX)$ and $\gB\subset\gR(\sY)$ be the  Boolean algebras of regular open sets with piecewise smooth boundaries, as defined in Example \ref{X:algebra}(b).    Let $\phi:\sX\into\sY$ be a {\em submersion} ---that is, a differentiable
function such that the derivative $\fD_x\phi$ is a linear surjection from the tangent space $\fT_x\sX$ to the tangent space $\fT_{\phi(x)}\sY$ for all 
 for all $x\in\sX$.  (This implies that dim$(\sX)\geq$dim$(\sY)$.)   Then the Open Mapping Theorem implies that $\phi$ is also open.
Suppose $\sB\in\gB$;  then  $\partial\sB = (\sH_1\intsct\partial\sB)\union\cdots\union(\sH_N\intsct\partial\sB)$
for some smooth hypersurfaces $\sH_1,\ldots,\sH_N\subseteq\sY$.  
For all $n\in[1\ldots N]$, let $\sH'_n:=\phi^{-1}(\sH_n)$; then 
$\sH'_n$ is a smooth hypersurface in $\sX$.\footnote{{\em Proof.} Suppose $\sH_n=\psi_n^{-1}\{0\}$ for some smooth function $\psi_n:\sX\into\Real$.  Then $\sH'_n:=(\psi_n\circ\phi)^{-1}\{0\}$.
If d$\psi_n$ is never zero and $\phi$ is a submersion, then the Chain Rule implies that d$(\psi_n\circ\phi)$ is also never zero.}
  Meanwhile, if $\sA:=\phi^{-1}(\sB)$, then
 $\partial\sA = \phi^{-1}(\partial\sB)$, because $\phi$ is open and continuous.  Thus, 
\beq
\partial \sA
&=& \phi^{-1}(\partial\sB)
\quad=\quad \phi^{-1}\lb[(\sH_1\intsct\partial\sB)\union\cdots\union(\sH_N\intsct\partial\sB)\rb]
\\&=&
(\phi^{-1}[\sH_1]\intsct \phi^{-1}[\partial\sB])\union\cdots\union(\phi^{-1}[\sH_N]\intsct \phi^{-1}[\partial\sB])
\\ &=&
(\sH'_1\intsct\partial\sA)\union\cdots\union(\sH'_N\intsct\partial\sA).
\eeq
Thus, $\sA$ has a piecewise smooth boundary;  in other words, $\sA\in\gA$.  This shows that
$\phi$ is $(\gA,\gB)$-measurable.  Thus,
Lemma \ref{measurability.lemma}(b) says that $\phi^{-1}:\gB\into\gA$ is a Boolean algebra homomorphism.
}

Continuing the notation of Lemma \ref{measurability.lemma}(b),
suppose that $\mu$ is a content on $\gA$.   Let $\phi:\sX\into\sY$ be an open, continuous,
 $(\gA,\gB)$-measurable function.   The {\dfn image} (or ``push-forward'') $\phi(\mu)$ is the function
$\nu:\gB\into[0,1]$ defined by setting $\nu[\sB]:=\mu[\phi^{-1}(\sB)]$ for
all $\sB\in\gB$.   Since $\phi^{-1}:\gB\into\gA$ is a Boolean algebra homomorphism by Lemma \ref{measurability.lemma}(b), it is immediate
that $\nu$ is a content on $\sY$.
The next result is the analog of the ``Change of Variables'' theorem for integration with respect to contents. 

\Proposition{\label{change.of.variables}}
{
 Let $\sX$ and $\sY$ be  topological spaces, and let $\gA\subseteq\gR(\sX)$ and   $\gB\subseteq\gR(\sY)$
  be Boolean subalgebras of the algebras of regular sets.  Let $\phi:\sX\into\sY$ be an open, continuous,
  $(\gA,\gB)$-measurable function, and let  $g:\sY\into\Real$ be any function.  
 \bthmlist
 \item   If $g\in\sF_\gB(\sY)$, then  $g\circ \phi\in\sF_\gA(\sX)$. 
  If $g\in\sL_\gB(\sY)$, then  $g\circ \phi\in\sL_\gA(\sX)$.  
 If $g\in\sG_\gB(\sY)$, then  $g\circ \phi\in\sG_\gA(\sX)$.  Finally,  if $g\in\CB(\sY)$, then $g\circ \phi\in\sC_\gA(\sX)$.
 \item Let $\mu$ be a content on $\gA$, and let $\nu:=\phi(\mu)$. For any $\sB\in\gB$, if $\sA:=\phi^{-1}(\sB)$,
 then $  \integral_{\sA} g\circ \phi \dmu = \integral_{\sB} g \dnu$.
    \ethmlist
    }
 \bthmprf (a)    Let $g\in\sF_\gB(\sY)$ be subordinate to some $\gB$-partition $\{\sB_1,\ldots,\sB_N\}$ of  $\sY$.  For all $n\in[1\ldots N]$, let $\sA_n:=\phi^{-1}(\sB_n)$;
 then $\sA_n\in\gA$ because $\phi$ is $(\gA,\gB)$-measurable.  Furthermore, $\{\sA_1,\ldots,\sA_N\}$
 is a $\gA$-partition of $\sX$ (because $\phi^{-1}:\gB\into\gA$ is a Boolean algebra homomorphism, by
  Lemma \ref{measurability.lemma}(b)),   and  $g\circ\phi$ is an
$\gA$-simple function subordinate to this partition.   Thus,   $g\circ \phi\in\sF_\gA(\sX)$.

  Now suppose $g\in\sL_\gB(\sY)$.   Then $g$ is a limit (in the uniform norm) of a sequence
$\{f_n\}_{n=1}^\oo$ where  $f_n\in\sF_\gB(\sY)$.  The transformation
$\Cb(\sY,\Real)\ni h\mapsto h\circ\phi\in\Cb(\sX,\Real)$ is continuous in the uniform norm.
Thus, $g\circ\phi=\lim_{n\goto\oo}
f_n\circ\phi$ (in the uniform norm).  By the previous paragraph, $f_n\circ\phi \in\sF_\gA(\sX)$ for all $n\in\Natur$.
Thus,  $g\circ\phi \in\sL_\gA(\sX)$.
 
Next, suppose  $g\in\CB(\sY)$.
 Let $r\in\Real$, and let $\sB:=\Int\lb[  g^{-1}[r,\oo)\rb]$.  Then $\sB\in\gB$, because $g$ is $\gB$-comeasurable.
 Thus,
 \beq
 \Int\lb[(g\circ \phi)^{-1}[r,\oo)\rb]
 &=&
 \Int\lb[ \phi^{-1}\lb( g^{-1}[r,\oo)\rb)\rb]
\  \ \eeequals{(*)} \ \
  \phi^{-1}\lb(\Int\lb[  g^{-1}[r,\oo)\rb]\rb)
   \  = \ 
  \phi^{-1}(\sB) \  \in  \ \gA,
  \eeq
  where $(*)$ is because $\phi$ is open and continuous, and the last step is
  because $\phi$ is $(\gA,\gB)$-measurable.  Likewise, 
  $\Int\lb[(g\circ \phi)^{-1}(-\oo,r]\rb]\in\gA$ for all $r\in\Real$.  Thus,
  $g\circ \phi$ is $\gA$-comeasurable.  Meanwhile, $g\circ\phi\in\Cb(\sX,\Real)$, because
  $g\in\Cb(\sY,\Real)$.  Thus,  $g\circ \phi\in\sC_\gA(\sX)$.

 Finally, suppose $g\in\GB(\sY)$.  There are two cases.   First
  suppose  $g=g_1+\cdots+g_N$ for some $g_1,\ldots,g_N\in\CB(\sY)$.
Then $g\circ\phi = g_1\circ\phi+\cdots+g_N\circ\phi$.  
But for all $n\in[1\ldots N]$,  $g_n\circ\phi\in\sC_\gA(\sX)$, by  the previous paragraph.  Thus, $g\circ\phi
\in\sG_\gA(\sX)$, as claimed.

 Now suppose  $g$ is an arbitrary element of $\GB(\sY)$.   Then $g$ is a limit (in the uniform norm) of a sequence
$\{g_n\}_{n=1}^\oo$ where each $g_n$ is as in  the previous paragraph.  The transformation
$\Cb(\sY,\Real)\ni h\mapsto h\circ\phi\in\Cb(\sX,\Real)$ is continuous in the uniform norm.
Thus, $g\circ\phi=\lim_{n\goto\oo}
g_n\circ\phi$ (in the uniform norm).  By   the previous paragraph, $g_n\circ\phi \in\sG_\gA(\sX)$ for all $n\in\Natur$.
Thus,  $g\circ\phi \in\sG_\gA(\sX)$, because $\sG_\gA(\sX)$ is closed in the uniform norm.

\item (b)  \  Let  $f\in\sF_\gB(\sB)$ be  subordinate to the $\gB$-partition $\{\sB_1,\ldots,\sB_N\}$.  For all $n\in[1\ldots N]$, let $r_n\in\Real$ be the value of $f$ on $\sB_n$, and let $\sA_n:=\phi^{-1}(\sB_n)$.
 Then as we argued in the proof of part (a), $f\circ\phi$ is an
$\gA$-simple function subordinate to the $\gA$-partition $\{\sA_1,\ldots,\sA_N\}$.  Clearly, $f\circ\phi$  takes the value
$r_n$ on $\sA_n$.   Thus, we have:
 \begin{eqnarray}
\nonumber
 \Integral_\sA (f\circ\phi)\dmu
 &\eeequals{(*)}&
 \sum_{n=1}^N r_n\,\mu[\sA_n]
 \quad=\quad
  \sum_{n=1}^N r_n\,\mu[\phi^{-1}(\sB_n)]
\\ &\eeequals{(\dagger)}&
    \sum_{n=1}^N r_n\,\nu[\sB_n]
 \quad\eeequals{(*)}\quad
 \Integral_\sB f\dnu.
  \label{change.of.variables.e1}
    \end{eqnarray}
Here, both $(*)$ are by   equation  (\ref{integral.of.simple.function}), while $(\dagger)$ is 
by the definition of $\nu=f(\mu)$.

Now, let $g\in\sL_\gB(\sB)$;  then $g\circ\phi\in\sL_\gA(\sA)$, by part (a).
  Let
$ \sF_g(\sB):=\{f\in\sF_\gB(\sB)$; \ $f(b)\leq g(b)$, for all $b\in\sB\}$, and let 
$ \sF_{g\circ\phi}(\sA):=\{f\in\sF_\gA(\sA)$; \ $f(a)\leq g\circ\phi(a)$, for all $a\in\sA\}$.
Then $f\circ\phi\in\sF_{g\circ\phi}(\sA)$ for every $f\in\sF_g(\sB)$.  Thus,
 \begin{eqnarray}
\nonumber
    \Integral_{\sA} g\circ\phi \dmu &\eeequals{(*)}& \sup_{f\in\sF_{g\circ\phi}(\sA)} \  \Integral_\sA f\dmu
    \quad\geeeq{(\dagger)}\quad
     \sup_{f\in\sF_g(\sB)} \  \Integral_\sA (f\circ \phi)\dmu
    \\ &\eeequals{(\diamond)}& 
      \sup_{f\in\sF_g(\sB)} \  \Integral_\sB f\dnu
\quad\eeequals{(*)}\quad
    \Integral_{\sB} g \dnu.
     \label{change.of.variables.e2}
    \end{eqnarray}
    Here, both $(*)$ are by  equation  (\ref{integral.defn}),
    $(\dagger)$ is because we have just observed that  $\{f\circ \phi$; \ $f\in\sF_g(\sB)\}\subseteq\sF_{g\circ\phi}(\sA)$,
    and $(\diamond)$ is by applying  (\ref{change.of.variables.e1}) to each $f\in\sF_g(\sB)$.  Meanwhile
   \beqn
 \label{change.of.variables.e3}
-\Integral_{\sA} g\circ\phi \dmu\quad=\quad
\Integral_{\sA} -g\circ\phi \dmu\quad \geeeq{(*)}\quad  \Integral_{\sB} -g \dnu
\quad=\quad
-\Integral_{\sB} g \dnu,
  \eeqn
  where $(*)$ is obtained  like  inequality (\ref{change.of.variables.e2}). 
  Multiplying both sides of  (\ref{change.of.variables.e3}) 
   by $-1$, we get
   \beqn
 \label{change.of.variables.e4}
 \Integral_{\sA} g\circ\phi \dmu
 \quad\leq\quad
 \Integral_{\sB} g \dnu.
\eeqn
Combining inequalities (\ref{change.of.variables.e2}) and (\ref{change.of.variables.e4}), we
get $\integral_{\sA} g\circ\phi \dmu=\integral_{\sB} g \dnu$, as desired.
 \ethmprf

\noindent Clearly, the composition of any two measurable functions is measurable.   
  Likewise, part (b) of the next result shows that the composition of two comeasurable functions is comeasurable, as long as it satisfies an auxiliary condition;  roughly speaking, the composite function must ``preserve negations''.
Meanwhile, part (a) says that  the composition of a comeasurable function with a measurable function is comeasurable.

\Proposition{\label{comeasurability.lemma}}
{
Let $\sX,\sY,\sZ$ be three topological spaces, and let $\gA\subseteq\gR(\sX)$, \  $\gB\subseteq\gR(\sY)$
and $\gC\subseteq\gR(\sZ)$  be Boolean subalgebras of the algebras of regular sets.
 Let $\phi:\sX\into\sY$ and $\psi:\sY\into\sZ$ be functions.
\bthmlist
\item  If $\phi$ is $(\gA,\gB)$-comeasurable, and $\psi$ is $(\gB,\gC)$-measurable, then $\psi\circ \phi$ is $(\gA,\gC)$-comeasurable.

\item   Suppose $\phi$ is continuous and $(\gA,\gB)$-comeasurable, while $\psi$ is continuous and $(\gB,\gC)$-comeasurable.
If $(\psi\circ \phi)^{\leftarrow}(\neg\sC) \subseteq \neg (\psi\circ \phi)^{\leftarrow}(\sC)$ for all
$\sC\in\gC$,  then $\psi\circ \phi$ is $(\gA,\gC)$-comeasurable.
\ethmlist
}
\bthmprf 
(a) \   Let $\sC\in\gC$.   Let $\sD:=\neg\sC$; then
 $\sD\in\gC$ also.  Thus, if we define $\sB:=\psi^{-1}(\sD)$, then $\sB\in\gB$, because $\psi$ is $(\gB,\gC)$-measurable, by hypothesis.
 But \beq
 \sD&=& \compl{\Cl(\sC)}.\\
\mbox{Thus,}\quad \sB&=&\psi^{-1}(\sD)\quad=\quad \psi^{-1}\lb[\compl{\Cl(\sC)}\rb]\quad=\quad\compl{\lb(\psi^{-1}\lb[\Cl(\sC)\rb]\rb)}.\\
\mbox{Thus,} \quad \psi^{-1}\lb[\Cl(\sC)\rb]&=&\compl{\sB}\quad=\quad\Cl(\sE),
\quad\mbox{where $\sE:=\neg\sB$ so that $\sE\in\gB$.}\\
\mbox{Thus,}\quad (\psi\circ \phi)^{-1}\lb[\Cl(\sC)\rb]&=& \phi^{-1}\lb(\psi^{-1}\lb[\Cl(\sC)\rb]\rb)\quad=\quad \phi^{-1}\lb[\Cl(\sE)\rb].\\
\mbox{Thus,}\quad \Int\lb((\psi\circ \phi)^{-1}\lb[\Cl(\sC)\rb]\rb)&=&\Int\lb(\phi^{-1}\lb[\Cl(\sE)\rb]\rb),
\eeq
which is an element of $\gA$,
 because $\phi$ is $(\gA,\gB)$-comeasurable, by hypothesis.

 \item (b) \  The proof uses the following claim.
 
 \Claim{\label{comeasurability.lemma.C1}
 For any $\sD\in\gC$, \ $\phi^\leftarrow \circ \psi^\leftarrow(\sD) \subseteq (\psi\circ \phi)^\leftarrow (\sD)$.}
 \bclaimprf
 Let $\sB:=\psi^\leftarrow(\sD)$.  Then $\sB\in\gB$ because $\psi$ is comeasurable.
 But $\sB= \Int\lb(\psi^{-1}\lb[\Cl(\sD)\rb]\rb)\subseteq \psi^{-1}\lb[\Cl(\sD)\rb]$, which is a closed set
 (because $\psi$ is continuous).  Thus, $\Cl(\sB)\subseteq\psi^{-1}\lb[\Cl(\sD)\rb]$.
 Thus, 
 \beq\phi^\leftarrow(\sB)&=&\Int\lb(\phi^{-1}\lb[\Cl(\sB)\rb]\rb)\quad\subseteq \quad
 \Int\lb(\phi^{-1}\lb[\psi^{-1}\lb(\Cl[\sD]\rb)\rb]\rb)
 \\ &=&
  \Int\lb((\psi\circ\phi)^{-1}\lb[\Cl(\sD)\rb]\rb)
\quad  =\quad (\psi\circ\phi)^\leftarrow(\sD).
  \eeq
  In other words, $\phi^\leftarrow \circ \psi^\leftarrow(\sD) \subseteq (\psi\circ \phi)^\leftarrow (\sD)$.
 \eclaimprf

 \noindent
  Now, let $\sC\in\gC$, and let $\sD:=\neg\sC$.  Then $\sD\in\gC$ also, and
 \beqn
 \label{comeasurability.lemma.e1}
  \phi^\leftarrow \circ \psi^\leftarrow(\sD) \quad\subseteeeq{(*)}\quad (\psi\circ \phi)^\leftarrow (\sD)
  \quad=\quad (\psi\circ \phi)^\leftarrow (\neg\sC)
  \quad\subseteeeq{(\dagger)}\quad \neg(\psi\circ \phi)^\leftarrow (\sC),
  \eeqn
  where $(*)$ is by Claim \ref{comeasurability.lemma.C1}, and $(\dagger)$ is by the hypothesis on $\psi\circ\phi$.
 Thus,
 \beqn
 \label{comeasurability.lemma.e2}
 (\psi\circ \phi)^\leftarrow (\sC)\quad\eeequals{(*)}\quad
\neg \neg(\psi\circ \phi)^\leftarrow (\sC)
\quad\subseteeeq{(\dagger)}\quad 
\neg\lb[\phi^\leftarrow \circ \psi^\leftarrow(\sD)\rb].
\eeqn
 Here, $(*)$ is because $(\psi\circ \phi)^\leftarrow (\sC)\in\gR(\sX)$ (because $\psi\circ\phi$ is continuous), and the
 negation operator $\neg$ is an involution on $\gR(\sX)$.  Meanwhile,
$(\dagger)$  is by negating both sides of (\ref{comeasurability.lemma.e1}) (thereby reversing the direction of inclusion).
But $\psi^\leftarrow(\sD)\in\gB$ because $\psi$ is $(\gB,\gC)$-measurable;  thus, 
$\phi^\leftarrow \circ \psi^\leftarrow(\sD)\in\gA$ because $\phi$ is
 $(\gA,\gB)$-comeasurable.  Thus, equation (\ref{comeasurability.lemma.e2}) implies that
 $ (\psi\circ \phi)^\leftarrow (\sC)\in\gA$, because $\gA$ is closed under negation.

This argument holds for all $\sC\in\gC$; thus,  $\psi\circ \phi$ is $(\gA,\gC)$-comeasurable.
 \ethmprf
  
  \noindent To see how the hypothesis of Proposition \ref{comeasurability.lemma}(b) could fail, suppose
  $(\psi\circ \phi)^{-1}(\partial\sC)$ contained an open subset $\sO$.  (Clearly, this could only happen if $\psi\circ \phi$ was not an open function.)  \ Then  $\sO\subseteq (\psi\circ \phi)^\leftarrow (\sC)$ and $\sO\subseteq (\psi\circ \phi)^\leftarrow (\neg\sC)$; so that $(\psi\circ \phi)^\leftarrow (\neg\sC)$
 and $(\psi\circ \phi)^\leftarrow (\sC)$ would be non-disjoint, and hence  $(\psi\circ \phi)^{\leftarrow}(\neg\sC) \not\subseteq \neg (\psi\circ \phi)^{\leftarrow}(\sC)$.

   If $\phi:\sX\into\sY$ is merely {\em co}measurable, but not measurable, then we cannot use $\phi$ to ``push forward'' a content $\mu$   from $\sX$ to $\sY$ as in Proposition \ref{change.of.variables}.  Nevertheless,  Proposition \ref{comeasurability.lemma}(a) still allows us to  ``push forward''  {\em integration} with respect $\mu$.  To see this, suppose $\sZ=\Real$ and $\gC=\gE$ (the Boolean algebra of 
   elementary functions from Example \ref{X:algebra}(a)).    
 If  $\phi:\sX\into\sY$ is $(\gA,\gB)$-comeasurable, and $g:\sY\into\Real$ is $(\gB,\gE)$-measurable,
 then Proposition \ref{comeasurability.lemma}(a) says that $g\circ\phi$ is $(\gA,\gE)$-comeasurable;
 hence $\integral_{\sA} g\circ\phi \dmu$ is well-defined for any $\sA\in\gA$.  This sort of computation plays a key role in  the  companion papers \cite{SEU_continuous,ImperfectPerception}, where $\mu$ represents the ``beliefs'' of a rational agent,
   $g$ represents  this agent's ``utility function'', and $\integral_{\sA} g\circ\phi \dmu$ is
 interpreted as ``expected utility''.

\refstepcounter{thm}
\paragraph{Remark \thethm.}\label{clopen.algebra.remark} The fact that the function $\phi$ is both continuous {\em and} open
is crucial in Lemma \ref{measurability.lemma} and Proposition \ref{change.of.variables}.  But
a continuous function $\phi$   need not  be open to induce a Boolean algebra homomorphism.
To see this, note that any clopen subset of a topological space $\sX$ is a regular open set.  The collection of all clopen sets of $\sX$ is a Boolean subalgebra of $\gR(\sX)$, in which the Boolean operations of $\gR(\sX)$ agree with the standard set-theoretic operations of union, intersection, and complementation.  Let $\sX$ and $\sY$ be two topological spaces,
and let $\clop(\sX)$ and $\clop(\sY)$ be their Boolean algebras of clopen sets.  If $\phi:\sX\into\sY$ is {\em any} continuous function, then it easy to verify that $\phi$ is measurable with respect to $\clop(\sX)$ and $\clop(\sY)$, and
$\phi^{-1}:\clop(\sY)\into\clop(\sX)$ is a Boolean algebra homomorphism.

\section{Liminal structures\label{S:liminal}}
\setcounter{equation}{0}

   Sections \ref{S:notation}, \ref{S:integration} and \ref{S:measurable}  considered   arbitrary Boolean subalgebras of $\gR(\sS)$,  and the integration of all real-valued functions that were integrable with respect to these subalgebras. 
  But  this section  and Section \ref{S:liminal.compactification}  only consider contents defined on the full Boolean algebra $\gR(\sS)$,  and the integration of {\em continuous} real-valued functions. 
    We are interested in
whether such contents, and their associated integrators, can be represented in terms of a traditional measure (either finitely  or countably additive) defined on some Boolean algebra or sigma-algebra of subsets of $\sS$.   Proposition  \ref{baire.algebra} already suggests one answer to this question.

\Proposition{\label{baire.integral}}
{
Suppose $\sS$ is a Baire space.  Let $\mu$ be a content on $\gR(\sS)$.
Let $\nu$ be the unique  residual  charge associated with $\mu$ by Proposition  \ref{baire.algebra}.
Then for any $\sB\in\gR(\sS)$, and any $g\in\Cb(\sS,\Real)$, we have
\[
  \Integral_{{\sB}} g \dmu\quad=\quad \int_{\sB} g\dnu.
  \]
 }
 \bthmprf  
 Let $\gB\!\ga(\sS)$ and $\gM(\sS)$ be as defined prior to  Proposition  \ref{residual.Borel.charge}.
 Recall that $\nu$ is a charge on $\gB\!\ga(\sS)$ such that $\nu(\sM)=0$ for all
 $\sM\in\gM(\sS)$, and such that $\nu[\sR]=\mu[\sR]$ for every $\sR\in\gR(\sS)$.
 Let $f: \sB\into\Real$ be a simple function, subordinate to a
 regular open partition $\gP=\{\sP_1,\sP_2,\ldots,\sP_N\}$  of $\sB$. For all $n\in[1\ldots N]$,
 suppose $f(p)=r_n$ for all  $p\in\sP_n$.   Then
 \beqn
 \label{baire.integral.e1}
 \int_{\sB} f\dnu
 \quad\eeequals{(*)}\quad
 \sum_{n=1}^N r_n \,\nu[\sP_n]
  \quad\eeequals{(\dagger)}\quad
 \sum_{n=1}^N r_n \,\mu[\sP_n]
  \quad\eeequals{(\diamond)}\quad
 \Integral_{\sB} f\dmu.
 \eeqn
 Here, $(*)$ is because  $\nu[\partial\sP_n]=0$ for all $n\in[1\ldots N]$, because $\partial\sP_n\in\gM(\sS)$.  Next, $(\dagger)$ is because $\nu[\sP_n]=\mu[\sP_n]$ for all
 $n\in[1\ldots N]$ by  construction,  and   $(\diamond)$ is by   equation  (\ref{integral.of.simple.function}).   
 
 Now let $g\in\Cb(\sS,\Real)$, and let $\sF_g(\sB)$ be the set of simple functions
 $f$ such that $f(b)\leq g(b)$ for all $b\in\sB$. Then
\[
\int_{\sB} g\dnu
\quad\eeequals{(*)}\quad
\sup_{f\in\sF_g(\sB) } \  \int_{\sB} f\dnu
\quad\eeequals{(\dagger)}\quad
\sup_{f\in\sF_g(\sB) } \  \Integral_{\sB} f\dmu
\quad\eeequals{(\diamond)}\quad
\Integral_{{\sB}} g \dmu, 
\]
as desired.  Here, $(*)$ is because $g$ can be uniformly approximated  from below by elements of $\sF_g(\sB)$,
$(\dagger)$ is by equation (\ref{baire.integral.e1}), and
$(\diamond)$ is by   equation (\ref{integral.defn})  from Theorem \ref{from.probability.to.expectation0} .
\ethmprf

\noindent The representation in Proposition \ref{baire.integral} is not entirely satisfactory, because
finitely additive charges can exhibit somewhat pathological behaviour.  If possible,
we would like to represent a content using a {\em countably additive} measure ---ideally, a Borel measure.     If we {\em must} use a charge, then we would like it to be well-behaved.
But  residual  charges can be badly behaved.  To see this,  let $\gU$ be an ultrafilter fixed at some point $s\in\sS$, and let $\del_\gU$ be the content in Example \ref{X:ultrafilter.example2}.  Let $\nu$ be
the  residual  charge associated to $\del_\gU$ by Proposition \ref{baire.algebra}.
Then for any regular open set $\sR$ containing $s$, we must have $\nu[\sR]=\del_\gU[\sR]=1$.
Since any open set is equal to a regular open set modulo some meager set, this implies that
$\nu[\sO]=1$ for any open neighbourhood of $s$.  But $\nu\{s\}=0$, because
the singleton $\{s\}$ is meager.  Thus, $\nu$ violates {\em normality} ---a basic ``continuity'' condition
for measures, which requires the measure of any set to be well-approximated from above by open sets 
and well-approximated from below by closed sets.   It would be better to represent contents using {\em normal} charges or Borel measures.  We will now construct 
such representations.

\breath

Let $\Borel(\sS)$ be the Borel sigma-algebra on $\sS$.   A {\dfn Borel probability measure} is a (countably additive) probability measure on $\Borel(\sS)$.
Since $\gR(\sS)$ is a subset of
$\Borel(\sS)$,  it is tempting to think that every Borel probability measure on $\sS$  defines a  content on $\gR(\sS)$.  
 Nonexample \ref{X:lebesgue.is.not.content} already showed that this is not the case.
Nevertheless, we might conversely hope that every content $\mu$ on $\gR(\sS)$ could be represented by a Borel measure $\nu$ in an essentially unique way, such that integration with respect to $\nu$ (in the classical sense) will be the same as integration with respect to $\mu$ (in the sense defined in Section \ref{S:integration}).
But as we shall now see, this is not the case either.

Let $\mu$ be a content on $\gR(\sS)$.
  Let $\nu$ be a Borel probability measure on $\gB(\sS)$.
  Let $\sF\subseteq\sC(\sS,\Real)$ be some collection of continuous functions.
  We will  say that $\nu$ satisfies the {\dfn Riesz representation property}  on $\sF$  if
  \beqn
\label{riesz1}
 \Integral_\sS f\dmu \quad=\quad \int_\sS f\dnu,
\qquad\mbox{for all $f\in\sF$.}
\eeqn
In this case, what is the relationship between $\nu$ and $\mu$?  What is the relationship between $\nu$ and $\bI$?
As the next examples show, this question does not have a simple answer.

\example{\label{X:charge.vs.measure1}
Let $\sS:=[-1,1]$.  Let $\gU_0\subseteq\gR[-1,1]$ be an ultrafilter fixed at $0$,
and define the content $\del_0:=\del_{\gU_0}$ as in 
  Example \ref{X:ultrafilter.example2}.
  The associated integrator is defined as follows:
for any $\sR\in\gR(\sS)$, and any $f\in\Cb(\sS,\Real)$,  
\[
\Integral_{\sR} f \dmu \quad=\quad\choice{   f(0) &\If &\sR\in\gU_0 \quad \mbox{(in particular, if $0\in\sR$)}; \\
0 &\If &\sR\not\in\gU_0 \quad \mbox{(in particular, if $0\in\neg\sR$)}.}
\]
 Heuristically, $\del_0$ is like a ``point mass'' at zero, but with an additional feature:  if the point 0 lies on the boundary between a regular set $\mathcal{R}$ and its negation $\neg\mathcal{R}$, then exactly {\em one} of $\mathcal{R}$ or $\neg\mathcal{R}$ can ``claim ownership'' of 0;  this decision is made by the ultrafilter $\mathfrak{U}_0$.  For example, exactly {\em one} of the following two statements is true:
 \bitem
 \item For all $\epsilon>0$, $\del_0[(0,\epsilon)]=1$ while $\del_0[(-\epsilon,0)]=0$.
 \item For all $\epsilon>0$, $\del_0[(0,\epsilon)]=0$ while $\del_0[(-\epsilon,0)]=1$.
 \eitem
  $\mathfrak{U}_0$ also decides the ``ownership'' of zero in more complicated cases.  For example, let 
\beq
\mathcal{E}_+ & := &  \bigsqcup_{n=1}^\infty \left(\frac{1}{2n+1},\frac{1}{2n}\right) \quad\mbox{and}\quad 
	\mathcal{O}_+ \quad := \quad  \bigsqcup_{n=1}^\infty \left(\frac{1}{2n},\frac{1}{2n-1}\right) \\
\mbox{while}\quad	
\mathcal{E}_- & := &  \bigsqcup_{n=1}^\infty \left(\frac{-1}{2n},\frac{-1}{2n+1}\right) \quad\mbox{and}\quad 
\mathcal{O}_- \quad := \quad  \bigsqcup_{n=1}^\infty \left(\frac{-1}{2n-1},\frac{-1}{2n}\right).
\eeq	
These are four disjoint regular open sets, with $ \mathcal{E}_+\vee\mathcal{O}_+\vee\mathcal{E}_-\vee\mathcal{O}_- =[-1,1] $.  Thus,  {\em one} of the four sets $\mathcal{E}_+$, $\mathcal{O}_+$, $\mathcal{E}_-$, and $\mathcal{O}_-$ gets $\del_0$-measure 1 (i.e. claims ``ownership'' of 0), while the other three get $\del_0$-measure 0 ---the ultrafilter $\mathfrak{U}_0$ decides which one.}

\example{\label{X:charge.vs.measure1b}
We will now refine Example \ref{X:charge.vs.measure1}.  Let
$\gF_+:=\{\sR\in\gR[-1,1]$; \ $(0,\eps)\subseteq\sR$ for some $\eps>0\}$.
Let $\gF_-:=\{\sR\in\gR[-1,1]$; \ $(-\eps,0)\subseteq\sR$ for some $\eps>0\}$.  Then
$\gF_+$ and $\gF_-$ are both free filters.
Let $\gU_+$ be a free ultrafilter containing $\gF_+$ and let $\gU_-$ be a free ultrafilter containing
$\gF_-$.  Define the contents $\del_\pm:\gR[-1,1]\into\Real$ using $\gU_\pm$ as in Example \ref{X:ultrafilter.example2}.  Fix some constants $\varphi_\pm\in[0,1]$ such that
$\varphi_-+\varphi_+=1$.   

 \quad Let $\gE$  be the Boolean algebra of elementary regular open subsets of $[-1,1]$, from Example \ref{X:algebra}(a).
 Let  $\lam'$  be the content on $\gE$ that agrees with the Lebesgue measure on $\gE$, as in Example \ref{X:elementary.lebesgue.content}.
   Use the Horn-Tarski Extension Theorem \cite[Theorem 1.22]{HornTarski48} to extend $\lam'$ to a content
 $\lam$ on all of $\gR[-1,1]$.
  
   For any $\sR\in\gR[-1,1]$, define
$\nu[\sR]:=\frac{1}{2}\lb(\lam[\sR]+\varphi_-\,\del_-[\sR]+\varphi_+\,\del_+[\sR]\rb)$.   Then $\nu$ is a content on
$\gR[-1,1]$.   

\quad If $0\in\neg \sR$ then $\nu[\sR]=\lam[\sR]/2$, whereas
if $0\in\sR$, then $\nu[\sR]=(\lam[\sR]+1)/2$.
  But if $0$ lies on the boundary between $\sR$ and $\neg\sR$,
then $\nu$'s behaviour is determined by $\gU_\pm$, and $\varphi_\pm$.
  In particular, for any $\eps>0$, we have
$\nu[(-\eps,0)] = ({\eps+\varphi_-})/{2}$ and
$\nu[(0,\eps)] = ({\eps+\varphi_+})/{2}$.

\quad Let $\lam_*$ be the Lebesgue measure (not to be confused with $\lam$), and 
let $\del_*$ be the Borel probability measure which assigns probability 1 to the singleton $\{0\}$ (not to be
confused with the contents  $\del_\pm$).  
Let $\nu:=\frac{1}{2}(\lam+\del_*)$.  Then $\nu$ is a  Radon measure on $[-1,1]$, and satisfies the  Riesz representation property (\ref{riesz1})  on $\sC(\sS,\Real)$   for {\em any} choice of ultrafilters $\gU_\pm$  and any choice of constants
 $\varphi_\pm$.  Thus, there are many different contents  which ``look like'' the same Radon measure on $\sS$,
  in the sense of equation (\ref{riesz1}).}

We will now develop a theory that explains these examples.
Let $\gA(\sS)$ be the Boolean algebra generated by the open subsets of $\sS$. 
 Thus, $\gA(\sS)$ contains all open subsets, all closed subsets, and all finite unions and intersections of such sets. 
For any $\sT\subseteq\sS$, let $\gA(\sT):=\{\sA\intsct\sT$; \ $\sA\in\gA(\sS)\}$;  this is the Boolean algebra  generated by  all relatively open and relatively closed subsets of $\sT$.
A {\dfn  charge} on $\gA(\sT)$ is a function $\nu:\gA(\sT)\into[0,1]$ which is finitely additive ---i.e. $\nu[\sA\disj\sB]=\nu[\sA]+\nu[\sB]$ for any disjoint $\sA,\sB\in\gA(\sT)$.   In particular, any Borel  measure on
$\sT$ restricts to a  charge on $\gA(\sT)$ in the obvious way.
We will say that $\nu$ is a {\dfn probability charge} on $\sT$ if $\nu[\sT]=1$.

  $\gR(\sS)$ is a subset of $\gA(\sS)$, but it is not a sub-{\em algebra}, because the join operator on $\gR(\sS)$ is not union.
Furthermore,  Examples \ref{X:charge.vs.measure1} and \ref{X:charge.vs.measure1b} show that not every charge on $\gA(\sS)$ determines a unique content on $\gR(\sS)$.  Nevertheless, if $\sS$ is a $T_4$ space, we will see that any content can be represented by a charge on $\gA(\sS)$, enriched with an auxiliary structure. 

\paragraph{Liminal charge structures.} ~
Let $\nu$ be a charge on $\gA(\sS)$.  For any $\sB\in\gA(\sS)$, let
$\nu_\sB$ be the restriction of $\nu$ to a charge on  $\gA(\sB)$.
A {\dfn liminal charge structure subordinate to $\nu$} is a collection $\{\rho_\sR\}_{\sR\in\gR(\sS)}$,
where for all $\sR\in\gR(\sR)$, $\rho_\sR$ is a charge on  $\gA(\partial\sR)$, absolutely continuous with respect to  $\nu$, such that
for any regular partition $\{\sR_1,\ldots,\sR_N\}$ of $\sS$,
\beqn
\label{charge.consistency.condition}
\rho_{\sR_1}+\cdots+\rho_{\sR_N}
\quad=\quad
\nu_{\partial\sR_1\union\cdots\union\partial\sR_N}.
\eeqn
In particular, for any $\sR\in\gR(\sS)$,  if $\sQ:=\neg\sR$ (so that $\partial\sQ=\partial\sR$), then $\rho_\sR+\rho_\sQ=\nu_{\partial\sR}$.  Heuristically,  $\rho_\sR$ and $\rho_\sQ$ describe the way that $\sR$ and $\sQ$ ``share'' the $\nu$-mass of their common boundary.  Given such a structure, we can define a function
$\mu:\gR(\sS)\into[0,1]$ by 
\beqn
\label{from.liminal.structure.to.content}
\mu[\sR]\quad:=\quad \nu(\sR)+\rho_\sR(\partial\sR),
\qquad\mbox{for all $\sR\in\gR(\sS)$.}
\eeqn
It is easy to verify that $\mu$ is a content.\footnote{{\em Proof sketch.} For any regular partition $\{\sR_1,\ldots,\sR_N\}$ of $\sS$, 
equations  (\ref{charge.consistency.condition})  and (\ref{from.liminal.structure.to.content}) together
yield $\mu[\sR_1]+\cdots+\mu[\sR_N]=1$.}  The first main result of this section establishes that, on any  $T_4$ topological space, {\em every} content arises in this fashion.\footnote{Recall that a Hausdorff topological space $\sS$ is {\dfn $T_4$} if, for any disjoint closed subsets $\sC_1,\sC_2\subset\sS$,  there exist disjoint open sets $\sO_1,\sO_2\subset\sS$ with $\sC_1\subseteq\sO_1$ and $\sC_2\subseteq\sO_2$.   For example, any metrizable space is $T_4$.}

A charge $\nu$ is {\dfn normal} if, for every $\sB\in\gA(\sS)$, we have
 $\nu[\sB]=\sup\{\nu[\sC]$; \ $\sC\subseteq\sB$ and $\sC$  closed in $\sS\}$
and  $\nu[\sB]=\inf\{\nu[\sO]$; \ $\sB\subseteq\sO\subseteq\sS$ and $\sO$ open in $\sS\}$. 
A liminal charge structure $\{\rho_\sR\}_{\sR\in\gR(\sS)}$ is {\dfn normal}  if
$\rho_\sR$ is a normal charge on $\partial\sR$ for all $\sR\in\gR(\sS)$.

\Theorem{\label{content.vs.Riesz.charge}}
{
Let $\sS$ be a $T_4$  space,  let $\mu$ be a content on $\gR(\sS)$.  
 There is a unique normal  probability charge $\nu$ on $\gA(\sS)$  that satisfies the Riesz representation property {\rm(\ref{riesz1})}  on $\Cb(\sS,\Real)$ .  Furthermore,
there is a unique normal liminal charge structure $\{\rho_\sR\}_{\sR\in\gR(\sS)}$ which is subordinate to $\nu$, such that  for any $\sR\in\gR(\sS)$, $\mu$ satisfies equation
{\rm(\ref{from.liminal.structure.to.content})}, and also
\beqn
\label{content.vs.Riesz.charge.parta}
\Integral_{\sR} f\dmu \quad=\quad  \int_{\sR} f \dnu + \int_{\partial\sR} f \drho_\sR,
\qquad\mbox{for all $f\in\Cb(\sS,\Real)$.}
\eeqn
}
In particular, if $\nu[\partial\sR]=0$, then   (\ref{from.liminal.structure.to.content}) and   (\ref{content.vs.Riesz.charge.parta}) say 
 $\mu[\sR]= \nu(\sR)$ and $\D\Integral_{{\sR}} f \dmu =  \int_{\sR} f \dnu$.

\bthmprf[Proof of Theorem \ref{content.vs.Riesz.charge}.]
 The integration operator $\integral_\sS\! \!\dmu$   is a 
  positive bounded linear functional on  $\Cb(\sS,\Real)$ .
Since $\sS$ is a normal Hausdorff space, 
a version of the Riesz Representation Theorem yields a unique normal probability charge $\nu$ on $\gA(\sS)$ 
that satisfies the Riesz representation property {\rm(\ref{riesz1})}  on $\Cb(\sS,\Real)$  
\cite[Theorem 14.9]{AliprantisBorder}.

 Let $\sR\in\gR(\sS)$.
Since $\nu$ is normal, 
for any $\eps>0$,  there is a closed set $\sK_\eps\subseteq\sR$ with
\beqn
\label{content.vs.Riesz.charge.e0}
\nu[\sK_\eps]\quad>\quad\nu[\sR]-\eps.
\eeqn
  Now, $\sK_\eps$ and $\compl{\sR}$ are disjoint closed sets,  and $\sS$ is T4, so
 Urysohn's Lemma  yields a continuous function
$\alp_\eps:\sS\into[0,1]$ such that 
\beqn
\label{content.vs.Riesz.charge.e0b}
\alp_\eps(\sK_\eps)\quad=\quad 1\quad \And \quad\alp_\eps(\compl{\sR})\quad=\quad 0.
\eeqn
Let $\bet_\eps:=1-\alp_\eps$.  Then for any $f\in\Cb(\sS,\Real)$, we have $f=\alp_\eps \, f+ \bet_\eps\, f$; thus
\beqn
\label{content.vs.Riesz.charge.e1}
\Integral_{\sR} f \dmu \quad=\quad \Integral_{\sR} \alp_\eps \, f+ \bet_\eps\, f \dmu \quad=\quad
\Integral_{\sR} \alp_\eps \, f \dmu+ \Integral_{\sR} \bet_\eps\, f \dmu.
\eeqn

\Claim{\label{content.vs.Riesz.charge.C1} For any $f\in\Cb(\sS,\Real)$, we have
$\D\Integral_{\sR} \alp_\eps\,f \dmu = \int_\sR \alp_\eps\,f \dnu$.}
\bclaimprf 
We have
\beq
\Integral_{\sS} \alp_\eps\,f \dmu &\eeequals{(*)}& \Integral_{\sR} \alp_\eps\,f \dmu+\Integral_{{\neg\sR}} \alp_\eps\,f \dmu
 \ \ \eeequals{(\dagger)} \ \  \Integral_{\sR} \alp_\eps\,f \dmu+\Integral_{{\neg\sR}} 0 \dmu
\\ &=&  \Integral_{\sR} \alp_\eps\,f \dmu+ 0, \\
\And \ \ 
\Integral_{\sS} \alp_\eps\,f \dmu &\eeequals{(\diamond)}& \int_\sS \alp_\eps\,f \dnu
\quad=\quad
 \int_\sR \alp_\eps\,f \dnu + \int_{\compl{\sR}} \alp_\eps\,f \dnu
\\&\eeequals{(\ddagger)}&
 \int_\sR \alp_\eps\,f \dnu + \int_{\compl{\sR}} 0\dnu
\quad=\quad
 \int_\sR \alp_\eps\,f \dnu +0.
 \eeq
 Combining these observations yields the claim.  Here $(*)$ is by equation (\ref{conditional.expectation0}), $(\diamond)$ is by the Riesz representation property {\rm(\ref{riesz1})}, and both
  $(\dagger)$ and $(\ddagger)$ are by  (\ref{content.vs.Riesz.charge.e0b}).
\eclaimprf

\Claim{\label{content.vs.Riesz.charge.C2}
For any $f\in\Cb(\sS,\Real)$, we have
$\D\lim_{\eps\goto0} \int_\sR \alp_\eps\,f \dnu=\int_\sR f \dnu$.}
\bclaimprf
\beq
\lb|\int_\sR f \dnu-\int_\sR \alp_\eps\,f \dnu\rb|
&=&
\lb|\int_\sR (1-\alp_\eps)\, f \dnu\rb|
\quad=\quad
\lb|\int_\sR \bet_\eps\, f \dnu\rb|\\
&\leq&
\norm{f}{\oo}\cdot \int_\sR \lb|\bet_\eps\rb| \dnu
\quad\leq\quad
\norm{f}{\oo}\cdot \nu[\sR\intsct\supp(\bet_\eps)] 
\\ &\leeeq{(*)} &
\norm{f}{\oo}\cdot \nu[\sR\setminus\sK_\eps] 
\quad\leeeq{(\diamond)}\quad
\norm{f}{\oo}\cdot \eps
\quad\goesto{\eps\goto0}{} \quad0,
\eeq
as desired.  Here, $(*)$ is by the defining properties (\ref{content.vs.Riesz.charge.e0b}) of the function $\alp_\eps$
(since $\bet_\eps=1-\alp_\eps$), while
 $(\diamond)$ is by inequality (\ref{content.vs.Riesz.charge.e0}).
\eclaimprf

\noindent Taking the limit as $\eps\goto0$ in equation (\ref{content.vs.Riesz.charge.e1}), and combining Claims \ref{content.vs.Riesz.charge.C1}
and \ref{content.vs.Riesz.charge.C2}, we obtain:
\beqn
\label{content.vs.Riesz.charge.e2}
\Integral_{\sR} f \dmu \quad=\quad 
\int_\sR f \dnu \ + \ \lim_{\eps\goto0} \ \Integral_{\sR} \bet_\eps\, f \dmu,
\qquad\mbox{for all $f\in\Cb(\sS,\Real)$.}
\eeqn

\Claim{\label{content.vs.Riesz.charge.C3}
There is a bounded positive linear functional
$\Phi_\sR:\Cb(\partial\sR,\Real)\into\Real$ with $\norm{\Phi_\sR}{\oo}\leq 1$,
such that
\beqn
\label{content.vs.Riesz.charge.e4}
\lim_{\eps\goto0}\ \Integral_{\sR} \bet_\eps\, f \dmu
\quad=\quad
\Phi_\sR(f_{\restr\partial\sR}),
\qquad\mbox{for all $f\in\Cb(\sS,\Real)$.}
\eeqn
}
\bclaimprf  The proof involves two subclaims.

\subclaim{\label{content.vs.Riesz.charge.C3A}$\lim_{\eps\goto0} \integral_{\sR} \bet_\eps\, f \dmu$ does not depend upon the
particular system of sets $\{\sK_\eps\}_{\eps>0}$ and functions $\{\alp_\eps\}_{\eps>0}$ that
we use in the above construction, as long as they satisfy the defining conditions
  {\rm(\ref{content.vs.Riesz.charge.e0})} and {\rm(\ref{content.vs.Riesz.charge.e0b})}.}
  \bsubclaimprf  
  Observe that
  equation (\ref{content.vs.Riesz.charge.e2}) can be rewritten:
  \[
  \lim_{\eps\goto0} \ \Integral_{\sR} \bet_\eps\, f \dmu\quad=\quad 
  \Integral_{\sR} f \dmu  - \int_\sR f \dnu,
  \]
   and the right hand side clearly does not depend upon  
   $\{\sK_\eps\}_{\eps>0}$ and  $\{\alp_\eps\}_{\eps>0}$. 
   \esubclaimprf
 \noindent 
 In the following argument, Claim \ref{content.vs.Riesz.charge.C3A} means that
   we can assume without loss of generality that $\{\sK_\eps\}_{\eps>0}$ and  $\{\alp_\eps\}_{\eps>0}$
   have whatever additional properties we require.

\subclaim{\label{content.vs.Riesz.charge.C3B}Let $f,g\in\Cb(\sS,\Real)$.  If $f_{\restr\partial\sR}=g_{\restr\partial\sR}$,
then $\D\lim_{\eps\goto0}\ \Integral_{\sR} \bet_\eps\, f \dmu=\lim_{\eps\goto0}\ \Integral_{\sR} \bet_\eps\, g \dmu$.}
\bsubclaimprf
Fix $\eps>0$.  For all $r\in\partial\sR$, there exists an open neighbourhood $\sB'_r\subseteq\sS$ around
$r$ such that $|f(b)-f(r)|<\frac{\eps}{2}$ for all $b\in\sB'_r$ (because $f$ is continuous).
Likewise, there exists an open neighbourhood $\sB''_r\subseteq\sS$ around
$r$ such that $|g(b)-g(r)|<\frac{\eps}{2}$ for all $b\in\sB''_r$ (because $g$ is continuous).
 Let $\sB^\eps_r:=\sB'_r\intsct\sB''_r$.  Then for all $b\in\sB^\eps_r$, 
\begin{eqnarray}
|f(b)-g(b)|
&=&
|f(b)-f(r)+f(r)-g(b)|
\quad = \quad 
|f(b)-f(r)+g(r)-g(b)| \nonumber \\
& \leq &
|f(b)-f(r)|+|g(r)-g(b)|
\quad<\quad
\frac{\eps}{2}+\frac{\eps}{2}\ \ = \ \ \eps.
\label{content.vs.Riesz.charge.e3}
\end{eqnarray}
Now let $\sB_\eps:=\Union_{r\in\partial\sR} \sB^\eps_r$.  Then $\sB_\eps$ is an open neighbourhood
of $\partial\sR$, and by combining the inequalites (\ref{content.vs.Riesz.charge.e3}) for all $r\in\partial\sR$, we
obtain
\beqn
|f(b)-g(b)|
\quad<\quad \eps,
\qquad\mbox{for all $b\in\sB_\eps$.}
\label{content.vs.Riesz.charge.e3b}
\eeqn
Let $\barsR$ denote the closure of $\sR$ in $\sS$.  Then $\barsR\setminus \sB_\eps$ is a closed
subset of $\sS$.  But $\barsR=\sR\union\partial\sR$, and
$\partial\sR\subseteq\sB_\eps$ by construction;  thus, $(\barsR\setminus \sB_\eps)\subseteq\sR$.
  By replacing $\sK_\eps$ with $\sK_\eps\union (\barsR\setminus\sB_\eps)$ if necessary,
we can assume without loss of generality that $(\sR\setminus\sK_\eps)\subseteq\sB_\eps$, for each
$\eps>0$.   Thus, 
\begin{eqnarray}
\lb|\maketall\Integral_{\sR} \bet_\eps\,f \dmu-\Integral_{\sR} \bet_\eps\,g \dmu\rb|
&=&
\lb|\maketall\Integral_{\sR} \bet_\eps\,f-\bet_\eps\,g \dmu\rb|
\quad=\quad
\lb|\maketall\Integral_{\sR} \bet_\eps\,(f-g) \dmu\rb| \nonumber\\
&\leq& \norm{\bet_\eps\,(f-g)}{\oo}
\quad\leeeq{(\diamond)}\quad
\sup_{s\in\supp(\bet_\eps)} |f(s)-g(s)| \nonumber
\\ &\leeeq{(*)}&
\sup_{b\in\sB_\eps} |f(b)-g(b)| 
\quad\leeeq{(\dagger)}\quad \eps.
\label{content.vs.Riesz.charge.e3a}
\end{eqnarray}
Here, $(\diamond)$ is because $\norm{\bet_\eps}{\oo}=1$,
 $(*)$ is because  $\supp(\bet_\eps)\subseteq\sR\setminus\sK_\eps\subseteq\sB_\eps$
 (by the defining conditions (\ref{content.vs.Riesz.charge.e0b})), and
 $(\dagger)$ is by inequality (\ref{content.vs.Riesz.charge.e3b}).

Letting $\eps\goto 0$ in  inequality (\ref{content.vs.Riesz.charge.e3a}), we obtain $\D\lim_{\eps\goto0}\ \Integral_{\sR} \bet_\eps\, f \dmu=\lim_{\eps\goto0}\ \Integral_{\sR} \bet_\eps\, g \dmu$. \esubclaimprf

\noindent  Let $\sC_\restr(\partial\sR,\Real):= \set{f_{\restr\partial\sR}}{f\in \Cb(\sS,\Real)}$.
 Claim \ref{content.vs.Riesz.charge.C3B} implies that we can define a function
$\Phi_\sR:\sC_\restr(\partial\sR,\Real)\into\Real$ by
equation (\ref{content.vs.Riesz.charge.e4}).   The linearity of $\Phi_\sR$ follows automatically from the linearity of
 the integration operator $\integral\dmu$  and the fact that $(r\,f+g)_{\restr\partial\sR} = r\,f_{\restr\partial\sR}+g_{\restr\partial\sR}$
for any $f,g\in  \Cb(\sR,\Real)$ and $r\in\Real$.  Likewise, $\Phi_\sR$ is positive because   $\integral\dmu$   is nondecreasing.

Next, recall that $\partial\sR$ is a closed subset of $\sS$,
and $\sS$ is $T_4$; thus, the Tietze Extension Theorem implies that 
$\sC_\restr(\partial\sR,\Real)=\Cb(\partial\sR,\Real)$. Thus, the function 
$\Phi_\sR$ is well-defined on all of $\Cb(\partial\sR,\Real)$.

Finally, to show that $\norm{\Phi_\sR}{\oo}\leq 1$, let $f\in\Cb(\partial\sR,\Real)$.
Suppose $\norm{f}{\oo}=M$, so we can think of $f$ as a function $f:\partial\sR\into[-M,M]$.
The Tietze Extension Theorem yields a continuous function $F:\sS\into[-M,M]$ such that
$F_{\restr\partial\sR}=f$.  Thus, for all $\eps>0$, we have $\norm{\bet_\eps\, F}{\oo}\leq M$,
and thus, $\lb|\integral_{\sR} \bet_\eps\, F \dmu\rb|  \leq M\cdot\mu(\sR) \leq M$ (because $\bI$ is compatible with $\mu$).  Thus,
\[ \lb|\Phi_\sR(f)\rb| \quad=\quad \lb|\lim_{\eps\goto0} \ \Integral_{\sR} \bet_\eps\, F \dmu\rb|
 \quad=\quad\lim_{\eps\goto0}  \ \lb|\Integral_{\sR} \bet_\eps\, F \dmu\rb|\quad \leq\quad  M\quad=\quad\norm{f}{\oo}.
 \]  
 This holds
 for all $f\in\Cb(\partial\sR,\Real)$, so $\norm{\Phi_\sR}{\oo}\leq 1$.
 \eclaimprf

\noindent Now, $\partial\sR$ is a closed subset of the $T_4$ space $\sS$;  thus, $\partial\sR$ is also $T_4$
\cite[Theorem 15.4(a)]{Willard}.
Thus, a version of the Riesz Representation Theorem  \cite[Theorem 14.9]{AliprantisBorder} yields a unique normal charge $\rho_\sR$ on  $\gA(\partial\sR)$ such that 
\[
\Phi_\sR(f)\quad=\quad\int_{\partial\sR} f\drho_\sR,
\qquad\mbox{for all $f\in\Cb(\partial\sR,\Real)$.}
\]
Combining this with equations (\ref{content.vs.Riesz.charge.e2}) and (\ref{content.vs.Riesz.charge.e4}),
we obtain  equation (\ref{content.vs.Riesz.charge.parta}).

Now, let $\sQ:=\neg \sR$.  Note that $\partial\sQ=\partial\sR$.
By repeating the above argument for $\sQ$, we obtain another normal charge $\rho_\sQ$ on $\gA(\partial\sQ)=\gA(\partial\sR)$, such that
\beqn
\label{content.vs.Riesz.charge.e5}
\Integral_{\sQ} f \dmu \quad=\quad 
\int_\sQ f \dnu + \int_{\partial\sR} f\drho_\sQ,
\qquad\mbox{for all $f\in\Cb(\sS,\Real)$.}
\eeqn
Thus, for any $f\in\Cb(\sS,\Real)$, we have
\begin{eqnarray}
\int_\sS f\dnu&\eeequals{(@)}& \int_\sR f \dnu + \int_{\partial\sR} f\dnu+\int_\sQ f \dnu,\quad\mbox{but also,}
\label{content.vs.Riesz.charge.e7}
\\ \nonumber
\int_\sS f\dnu
&\eeequals{(*)}&
\Integral_{\sS} f \dmu
\quad\eeequals{(\dagger)}\quad
\Integral_{\sR} f \dmu+\Integral_{\sQ} f \dmu
\\ &\eeequals{(\diamond)}&
\int_\sR f \dnu + \int_{\partial\sR} f\drho_\sR+\int_\sQ f \dnu + \int_{\partial\sR} f\drho_\sQ,
\label{content.vs.Riesz.charge.e8}
\end{eqnarray}
where $(@)$ is because $\sS=\sR\disj\partial\sR\disj\sQ$, 
 $(*)$ is by the  Riesz representation property (\ref{riesz1}), $(\dagger)$ is by
equation (\ref{conditional.expectation0}),
and $(\diamond)$ is by the equations (\ref{content.vs.Riesz.charge.parta}) and (\ref{content.vs.Riesz.charge.e5}).

Subtracting  (\ref{content.vs.Riesz.charge.e7}) from (\ref{content.vs.Riesz.charge.e8}) and rearranging, we obtain:
\[
\int_{\partial\sR} f\dnu_{\partial\sR}
 \ \ = \ \ 
\int_{\partial\sR} f\drho_\sR \ + \ \int_{\partial\sR} f\drho_\sQ
 \ \ = \ \ 
\int_{\partial\sR} f\ \mathrm{d}(\rho_\sR+\rho_\sQ),
\quad\mbox{for all $f\in\Cb(\sS,\Real)$.}
\]
However, as we earlier noted, the Tietze Extension Theorem implies
that $\Cb(\partial\sR,\Real):=\{f_{\restr\partial\sR}$; \ $f\in\Cb(\sS,\Real)\}$.  Thus, we obtain
\beqn
\label{content.vs.Riesz.charge.e6a}
\int_{\partial\sR} f\dnu_{\partial\sR}
 \ \ = \ \ 
\int_{\partial\sR} f\ \mathrm{d}(\rho_\sR+\rho_\sQ),
\quad\mbox{for all $f\in\Cb(\partial\sR,\Real)$.}
\eeqn
Now, $\nu_{\partial\sR}$ is normal because it is a restriction of the normal charge $\nu$ to $\gA(\partial\sR)$, while $(\rho_\sR + \rho_\sQ)$ is normal because it is a sum of two normal charges on $\gA(\partial\sR)$.  Thus, statement (\ref{content.vs.Riesz.charge.e6a}) and the uniqueness part of the Riesz Representation Theorem  yield 
\beqn
\label{content.vs.Riesz.charge.e6}
\rho_\sR + \rho_\sQ  \quad=\quad \nu_{\partial\sR}.
\eeqn

\Claim{\label{content.vs.Riesz.charge.C4}
 $\rho_\sR$ and   $\rho_\sQ$ are absolutely continuous relative to $\nu_{\partial\sR}$.
}
\bclaimprf
Let $\sU\subseteq\partial\sR$, and suppose that $\nu_{\partial\sR}[\sU]=0$.
Then equation (\ref{content.vs.Riesz.charge.e6}) implies that $\rho_\sR[\sU] + \rho_\sQ [\sU]=0$.  
Since these are both positive measures, this means that $\rho_\sR[\sU] = \rho_\sQ [\sU]=0$.
This conclusion holds whenever $\nu_{\partial\sR}[\sU]=0$.  Thus, 
 $\rho_\sR$ and   $\rho_\sQ$ are absolutely continuous relative to $\nu_{\partial\sR}$.
\eclaimprf
\noindent
 Equation (\ref{content.vs.Riesz.charge.e6})
is obviously a special case of the  equation (\ref{charge.consistency.condition}) for
the two-element partition $\{\sR,\sQ\}$.  To prove  equation (\ref{charge.consistency.condition}) in general, 
 let $\{\sR_1,\ldots,\sR_N\}$ be any regular open partition of $\sS$, and generalize equations (\ref{content.vs.Riesz.charge.e7}), (\ref{content.vs.Riesz.charge.e8}) and (\ref{content.vs.Riesz.charge.e6a}) in the obvious way.
 
 It remains to establish  formula (\ref{from.liminal.structure.to.content}).  Let
$\sR\in\gR(\sS)$.  Let $\bone$ be the constant 1-valued function;  then $\bone\in\Cb(\sS,\Real)$, and we have
\[
\mu[\sR]
\quad=\quad
\Integral_{\sR} \bone \dmu
\quad\eeequals{(*)}\quad
\int_\sR \bone \dnu + \int_{\partial\sR} \bone\drho_\sR
\quad=\quad
\nu[\sR] + \rho_\sR[\partial\sR],
\]
as desired.  Here $(*)$ is by  equation
(\ref{content.vs.Riesz.charge.parta}).
\ethmprf

\noindent{\bf Liminal density structures.}
 If $\sS$ is locally compact, contents admit a  nicer representation. 
Let $\Co(\sS,\Real)$ be the Banach space of all  functions in $\Cb(\sS,\Real)$ that ``vanish at  infinity'',
meaning that for any $\eps>0$, there is a compact subset $\sK\subseteq\sS$ such that 
$|f(s)|\leq\eps$ for all $s\in\sS\setminus\sK$.  
Let $\gR_0(\sS)$ be the set of all regular open subsets of $\sS$ with compact closures.

Let  $\nu$ be a Borel  measure on $\sS$.
For any  $\sB\in\Borel(\sS)$, 
  let $\dL^1(\sB,\nu)$ denote the Banach space  of real-valued, $\nu$-integrable functions on $\sB$,
 modulo equality $\nu$-almost everywhere.\footnote{If $\nu[\sB]=0$, then $\dL^1(\sB,\nu)$ is trivial.}   Let
 $\dL^1(\sB,\nu;[0,1])$ be the set of $[0,1]$-valued functions in $\dL^1(\sB,\nu)$.
A {\dfn liminal density structure subordinate to $\nu$} is a collection $\{\phi_\sR\}_{\sR\in\gR(\sS)}$,
where, for all $\sR\in\gR(\sR)$,  $\phi_\sR\in\dL^1(\partial\sR,\nu;[0,1])$ is a function such that, 
for any regular partition $\{\sR_1,\ldots,\sR_N\}$ of $\sS$, we have
\beqn
\label{charge.consistency.condition2}
\phi_{\sR_1}+\cdots+\phi_{\sR_N}
\quad=\quad
1,  \quad\mbox{$\nu$-almost everywhere on $\partial\sR_1\union\cdots\union\partial\sR_N$.}
\eeqn
In particular, for any $\sR\in\gR(\sS)$,  if $\sQ:=\neg\sR$  (so that $\partial\sQ=\partial\sR$), then (\ref{charge.consistency.condition2}) implies that   $\phi_\sQ=1-\phi_\sR$,
$\nu$-a.e. on $\partial\sR$.
 Heuristically,  $\phi_\sR$ and $\phi_\sQ$ describe the way in which $\sR$ and $\sQ$ ``share'' the $\nu$-mass of their common boundary.  Given such a structure, we can define a function
 $\mu:\gR(\sS)\into[0,1]$ by setting
  \beqn
\label{content.vs.Riesz.charge.partb2}
\mu[\sR]\quad=\quad \nu(\sR)+\int_{\partial\sR}  \phi_\sR \dnu_{\partial\sR},
\eeqn
for all $\sR\in\gR(\sS)$.  It is easy to verify that $\mu$ is a content.    We will  soon  see that, if $\sS$ is a compact Hausdorff space, then {\em every} content arises in this fashion.  But first we will prove a slightly weaker result, for locally compact $T_4$  spaces.
 Recall that a  Borel measure $\nu$ is {\dfn Radon} if, for all $\sB\in\Borel(\sS)$,  
we have $\nu[\sB]=\sup\{\nu[\sK]$; \ $\sK\subseteq\sB$ and $\sK$ compact in $\sS\}$, and
 $\nu[\sB]=\inf\{\nu[\sO]$; \ $\sB\subseteq\sO\subseteq\sS$ and $\sO$ open in $\sS\}$.

\Theorem{\label{content.vs.Riesz.charge2}}
{
Let $\sS$ be a locally compact $T_4$  space,  let $\mu$ be a content on $\gR(\sS)$.  
There is a unique Radon  measure $\nu$ on $\sS$ that satisfies the Riesz representation property {\rm(\ref{riesz1})}  on $\Co(\sS,\Real)$ .
Furthermore, there is a unique liminal density  structure $\{\phi_\sR\}_{\sR\in\gR(\sS)}$ 
such that $\mu$ satisfies equation {\rm(\ref{content.vs.Riesz.charge.partb2})}
for all $\sR\in\gR_0(\sS)$, while
\beqn
\label{content.vs.Riesz.charge.partb}
  \Integral_{{\sR}} f \dmu \quad=\quad  \int_{\sR} f \dnu + \int_{\partial\sR} f \, \phi_\sR \dnu_{\partial\sR},
\eeqn
 for all $f\in\Co(\sS,\Real)$ and $\sR\in\gR(\sS)$.   Also,
  {\rm(\ref{content.vs.Riesz.charge.partb})} holds for all  $f\in\sC(\sS,\Real)$ and
 $\sR\in\gR_0(\sS)$.
 }
  \bthmprf  
 The integration operator  $\integral\!\!\dmu$  
 is a   positive bounded linear functional on  $\Co(\sS,\Real)$, with  $\norm{\integral\dmu}{\oo}=1$.    
Since $\sS$ is a locally compact Hausdorff space,   the Riesz Representation Theorem yields a unique Radon probability measure $\nu$ on $\Borel(\sR)$ 
that satisfies the Riesz representation property {\rm(\ref{riesz1})}  on $\Co(\sS,\Real)$  
\cite[Theorem 7.17]{Folland}.   

At this point, the proof is very similar to the proof of  Theorem \ref{content.vs.Riesz.charge2}, except with the words ``(normal) charge'' everywhere replaced by ``(Radon) measure'', and with $\Cb(\sS,\Real)$ replaced by
$\Co(\sS,\Real)$.    However, there is a subtle change in the argument immediately
after the proof of Claim \ref{content.vs.Riesz.charge.C3}.   Now, for any $\sR\in\gR(\sS)$,    $\partial\sR$ is a closed subset of the {\em locally compact} Hausdorff space $\sS$;  thus, $\partial\sR$ is also a locally compact Hausdorff space
\cite[Theorem 18.4]{Willard}.
Thus, the Riesz Representation Theorem yields a unique Radon measure $\rho_\sR$ on $\Borel(\partial\sR)$ such that 
\[
\Phi_\sR(f)\quad=\quad\int_{\partial\sR} f\drho_\sR,
\qquad\mbox{for all $f\in\Co(\partial\sR,\Real)$.}
\]
Combining this with the relevant versions of equations (\ref{content.vs.Riesz.charge.e2}) and (\ref{content.vs.Riesz.charge.e4}),
we obtain equation (\ref{content.vs.Riesz.charge.parta}) for all $f\in\Co(\partial\sR,\Real)$.

The conclusion of Claim \ref{content.vs.Riesz.charge.C4} is still true.  But now, since we are now dealing with
(sigma-additive) measures, we can apply the Radon-Nikodym Theorem  to obtain non-negative functions
$\phi_\sR$ and $\phi_\sQ\in\dL^1(\partial\sR,\nu)$ such that
$\drho_\sR=\phi_\sR\dnu_{\partial\sR}$ and $\drho_\sQ=\phi_\sQ\dnu_{\partial\sR}$
(where $\nu_{\partial\sR}$ is the restriction of $\nu$ to $\Borel(\partial\sR)$).  Substituting this into
 equation (\ref{content.vs.Riesz.charge.parta}) yields   (\ref{content.vs.Riesz.charge.partb})  for all $f\in\Co(\sS,\Real)$.
Finally, equation (\ref{content.vs.Riesz.charge.e6}) implies that $\phi_\sR(s)+\phi_\sQ(s) =1$
for $\nu$-almost all $s\in\partial\sR$, as claimed.
Since $\phi_\sR$ and $\phi_\sQ$ are non-negative, this  implies that $0\leq\phi_\sR(s),\phi_\sQ(s)\leq 1$
for $(\nu)$-almost all $s\in\partial\sR$.  So
$\phi_\sR,\phi_\sQ\in\dL^1(\partial\sR,\nu;[0,1])$, as claimed.

 To prove  equation (\ref{charge.consistency.condition2}) in general, 
 let $\{\sR_1,\ldots,\sR_N\}$ be any regular open partition of $\sS$, and generalize equations (\ref{content.vs.Riesz.charge.e7}), (\ref{content.vs.Riesz.charge.e8}) and (\ref{content.vs.Riesz.charge.e6a}) in the obvious way. 

  Next we establish formula (\ref{content.vs.Riesz.charge.partb})
for any $\sR\in\gR_0(\sS)$ and $f\in\sC(\sS,\Real)$.
   If $\sR\in\gR_0(\sS)$, then its closure $\barsR$ is compact.   For all $r\in\barsR$, let $\sO_r$ be an open neighbourhood
of $r$ whose closure $\barsO_r$ is compact (this exists because $\sS$ is locally compact and Hausdorff).  
The collection $\{\sO_r$; \ $r\in\barsR\}$ is an open cover of the compact set $\barsR$, so
it admits a finite subcover, say $\{\sO_{r_1},\ldots,\sO_{r_N}\}$.  Let
$\sU:=\sO_{r_1}\union \cdots\union\sO_{r_N}$; this is an open set containing $\barsR$, and
$\barsU=\barsO_{r_1}\union\cdots\union\barsO_{r_N}$ is a finite union of compact sets, hence compact.
Now, $\barsR$ and $\compl{\sU}$ are disjoint closed subsets of the $T_4$ space $\sS$;  thus,
Urysohn's Lemma yields a function $ h\in\Cb(\sS,\Real)$ such that
$ h(r)=1$ for all $r\in\barsR$, while $ h(s)=0$ for all $s\in\compl{\sU}$.  

  Let $f\in\sC(\sS,\Real)$.  Then
$\supp( h\, f)\subseteq\barsU$, a compact set;  thus, $ h\, f\in\Co(\sS,\Real)$.   Thus,
\beq
\Integral_{\sR} f \dmu
&\eeequals{(*)}&
\Integral_{\sR}  h\, f \dmu
\quad\eeequals{(\dagger)}\quad
\int_\sR  h\, f \dnu + 
\int_{\partial\sR}  h\, f \cdot  \phi_\sR \dnu_{\partial\sR}
\\ &\eeequals{(*)} & 
\int_\sR f \dnu + 
\int_{\partial\sR} f \cdot \phi_\sR \dnu_{\partial\sR},
\eeq
which yields the desired equation (\ref{content.vs.Riesz.charge.partb})  for $f$. 
 Here, both $(*)$ are because $ h(r)=1$ for all $r\in\barsR$, by construction, while
 $(\dagger)$ is by  equation (\ref{content.vs.Riesz.charge.partb}),  which we can already apply to $h\, f$
 because it is an element of $\Co(\sS,\Real)$.

 Finally, setting $f=\bone$ in equation (\ref{content.vs.Riesz.charge.partb}), we obtain:
\[
\mu[\sR] \ = \ \Integral_{\sR} \bone \dmu \ = \ \int_\sR \bone \dnu + 
\int_{\partial\sR} \bone \cdot \phi_\sR \dnu_{\partial\sR}
\ = \ \nu[\sR] + 
\int_{\partial\sR}   \phi_\sR \dnu_{\partial\sR},
\]
which yields
 (\ref{content.vs.Riesz.charge.partb2})  for any $\sR\in\gR_0(\sS)$.
\ethmprf

\noindent  
Note that  $\nu$ is not necessarily a {\em probability} measure in Theorem \ref{content.vs.Riesz.charge2},
even though we  assumed $\mu[\sS]=1$.  The set
$\sS$ itself is not necessarily in $\gR_0(\sS)$ (unless $\sS$ is compact), so 
we cannot apply formula  (\ref{content.vs.Riesz.charge.partb2})  to get
$\nu[\sS]=\mu[\sS]=1$.  The next example illustrates this, and also shows why we 
restrict  formulae (\ref{content.vs.Riesz.charge.partb2})  and 
(\ref{content.vs.Riesz.charge.partb}) to $\gR_0(\sS)$ and  $\Co(\sS,\Real)$.  It also shows that
the representation in Theorem \ref{content.vs.Riesz.charge2} is not always very informative.

\example{\label{X:bad.content}
Let $\sS=(0,1)$; this space is  locally compact and $T_4$ .
Let $\gF:=\{\sR\in\gR(0,1)$; $(0,\eps)\subseteq\sR$ for some $\eps>0\}$.
Then $\gF$ is a free filter.  Use the Ultrafilter Theorem to extend $\gF$ to a free ultrafilter $\gU\subset\gR(0,1)$,
and then define the content $\del_\gU$ as in Example \ref{X:ultrafilter.example2}.

  It is  easy to see that $\del_\gU$ can be represented by a
(finitely additive) charge on $\gA(0,1)$, as described in Theorem \ref{content.vs.Riesz.charge}.
However, there is clearly no (countably additive) Borel measure on $(0,1)$ which
satisfies (\ref{content.vs.Riesz.charge.partb2}) for all $\sR\in\gR(0,1)$.  Indeed,
any set in $\gR_0(0,1)$ must be bounded away from zero.  Likewise, any function
$f\in\Co((0,1);\Real)$ must have the property that $\lim_{s\goto0} f(s)=0$.
Thus, there is only one Borel measure which satisfies  (\ref{content.vs.Riesz.charge.partb2}) 
for all $\sR\in\gR_0(0,1)$ and satisfies  (\ref{content.vs.Riesz.charge.partb}) 
for all $f\in\Co((0,1);\Real)$ ---namely, the zero measure.}

  If  
  $\sS$ is a compact space, then pathological examples like these cannot exist.

\Corollary{\label{content.vs.Riesz.charge3}}
{
Let $\sS$ be a compact Hausdorff  space,  let $\mu$ be a content on $\gR(\sS)$.  
There is a unique  normal  Borel  probability measure $\nu$ on $\sS$ that satisfies the Riesz representation property {\rm(\ref{riesz1})} on $\sC(\sS,\Real)$ .
Furthermore, there is a unique liminal density  structure $\{\phi_\sR\}_{\sR\in\gR(\sS)}$ which is subordinate to $\nu$ and which satisfies equations {\rm(\ref{content.vs.Riesz.charge.partb2})} and
{\rm(\ref{content.vs.Riesz.charge.partb})}
for all $\sR\in\gR(\sS)$ and all $f\in\sC(\sS,\Real)$.  
}
\bthmprf
 If   $\sS$ is compact, then $\gR_0(\sS)=\gR(\sS)$ and $\Co(\sS,\Real)=\sC(\sS,\Real)$,  while
 a Borel measure $\nu$ is Radon if and only if it is normal.  Now  apply 
Theorem \ref{content.vs.Riesz.charge2}.
\ethmprf

\example{(a) Let $\gU\subset\gR(\sS)$ be an ultrafilter fixed at some point $s$ in $\sS$, and define
$\del_\gU$ as in Example \ref{X:ultrafilter.example2}.  To satisfy equations {\rm(\ref{content.vs.Riesz.charge.partb2})} and
{\rm(\ref{content.vs.Riesz.charge.partb})},
let $\nu$ be the Borel probability measure which assigns probability 1 to $\{s\}$ (i.e. the ``point mass'' at $s$).  For any $\sR\in\gR(\sS)$  with $s\in\partial\sR$, define $\phi_{\sR}:=1$ if $\sR\in\gU$,
and $\phi_\sR:=0$ if $\sR\not\in\gU$.   (If $s\not\in\partial\sR$, then  $\nu(\partial\sR)=0$, so  the values of $\phi_\sR$ and $\phi_{\neg\sR}$ are irrelevant.)

(b) Let $\sS\in\gR(\Real^N)$, and    let $\nu$ be the (normalized) Lebesgue measure on $\sS$.  Let $\gB_{\mathrm{jor}}(\sS)$ be the Jordan algebra from Example \ref{X:algebra}(d).  Then $\nu$ restricted to $\gB_{\mathrm{jor}}(\sS)$ is a content.  Use the 
Horn-Tarski Extension Theorem \cite[Theorem 1.22]{HornTarski48} to extend this to a content $\mu$ on all of $\gR(\sS)$. In this case, clearly $\nu$ is the measure needed to satisfy  
 equations {\rm(\ref{content.vs.Riesz.charge.partb2})} and
{\rm(\ref{content.vs.Riesz.charge.partb})} for $\mu$.    If $\sR\in\gB_{\mathrm{jor}}(\sS)$, then $\nu[\partial\sR]=0$,  so the values of
 $\phi_\sR$ and $\phi_{\neg\sR}$ are irrelevant.    But if $\sR\not\in\gB_{\mathrm{jor}}(\sS)$,
 then $\nu[\partial\sR]>0$.    In this case, the functions
 $\phi_\sR$ and $\phi_{\neg\sR}$ describe how the nonzero Lebesgue measure of $\partial\sR$ is  ``shared'' between $\sR$ and $\neg\sR$.  Thus,
 $\mu[\sR]\geq \nu[\sR]$ and $\mu[\neg\sR]\geq \nu[\neg\sR]$, with at least one of these
 inequalities being strict. 
  }

\noindent Can we eliminate
the ``liminal'' terms in equations {\rm(\ref{content.vs.Riesz.charge.partb2})} and
{\rm(\ref{content.vs.Riesz.charge.partb})}?   Perhaps if the content $\mu$
was particularly nice, or if the Borel measure $\nu$ was ``smooth enough'', then these
terms would vanish.  In this case,
 equation (\ref{content.vs.Riesz.charge.partb}) would reduce to
 the more familiar expression:
\[
\Integral_{{\sR}} f \dmu \quad=\quad  \int_{\sR} f \dnu.
\]
Meanwhile,  equation (\ref{content.vs.Riesz.charge.partb2}) would say that $\mu[\sR]=\nu[\sR]$ for all $\sR\in\gR(\sS)$ ---in other words,
$\nu$ would define a content when restricted to $\gR(\sS)$.  
 Nonexample \ref{X:lebesgue.is.not.content} already shows that this is {\em not} the case for the Lebesgue measure on $[0,1]$.    But the next result goes much further: it is {\em never} the case, for any Borel measure on a broad class of topological spaces. 
 
      Say that a topological space $\sS$ is {\dfn projectible} if there
 is an open, continuous function from $\sS$ to $[0,1]$.  For example, any open subset of a topological vector space is projectible.  Also, any fibre bundle over any open subset of $[0,1]$ (with any fibre space) is projectible.  We will say that $\sS$ is {\dfn locally projectible} if every point in $\sS$ has a regular open neighbourhood which is projectible.  For example, any topological manifold is locally projectible.

\Proposition{\label{X:lebesgue.is.not.content2}}
{
If $\sS$ is locally projectible, then
there is no Borel measure on $\sS$ that defines a content when restricted to $\gR(\sS)$.
Thus, if $\mu$ is any content on $\gR(\sS)$ with a liminal density representation {\rm(\ref{content.vs.Riesz.charge.partb2})}, then the liminal density structure  $\{\phi_\sR\}_{\sR\in\gR(\sS)}$ is nontrivial. 
}

In particular, if $\sS$ is any bounded open subset of $\Real^N$, then Proposition \ref{X:lebesgue.is.not.content2} says that the Lebesgue measure on $\sS$ cannot define a content on $\gR(\sS)$.
Another consequence of this result is that the Borel measure which appears in Theorem \ref{content.vs.Riesz.charge2} and Corollary \ref{content.vs.Riesz.charge3}
{\em must} be different than the  residual  charge which appears in Propositions \ref{baire.algebra} and 
\ref{baire.integral}.  If $\sS$ satisfies the hypothesis of Proposition \ref{X:lebesgue.is.not.content2},
then no Borel measure on $\sS$ can be a  residual  charge, even if they represent the same content.

\bthmprf   (Case 1) Suppose $\sS$ is projectible. Let $\phi:\sS\into[0,1]$ be an open, continuous function.
Let $\nu'$ be a Borel probability measure on $\sS$, and let
$\mu'$ be the restriction of $\nu'$ to $\gR(\sS)$.
Let $\nu:=\phi(\nu')$;  this is a Borel probability measure on $[0,1]$.
Let $\mu:=\phi(\mu')$;  this is the restriction of $\nu$ to $\gR[0,1]$.
 Now,  $\phi^{-1}:\gR[0,1]\into\gR(\sS)$ is a Boolean algebra homomorphism
\cite[4A2B (f)(iii)]{Fremlin4II} (or Lemma \ref{measurability.lemma}(b)).    Thus, if
$\mu'$ is a content on $\gR(\sS)$, then
$\mu$ is a content on $\gR[0,1]$.  So to prove the theorem,
it suffices to show:
\bquote
{\em There is no Borel measure $\nu$ on $[0,1]$ such that $\mu$ is a content on $\gR[0,1]$.}
\equote

Our proof strategy is somewhat similar to the strategy sketched for Nonexample \ref{X:lebesgue.is.not.content}, but more general, since it must work for any Borel measure on $[0,1]$.

\Claim{If $\nu$ has an atom in $(0,1)$, then $\mu$ is not a content.\label{X:lebesgue.is.not.content2.C1}}
\bclaimprf
Suppose $\nu$ has an atom at some point $x\in(0,1)$.  Let $\sL:=[0,x)$ and $\sR:=(x,1]$.
Then clearly, $\sL$ and $\sR$ are regular open subsets of $[0,1]$ with
$\sL\vee\sR=[0,1]$.  But
 \[
\mu[\sL]+\mu[\sR]\quad=\quad\nu[\sL]+\nu[\sR]\quad=\quad
\nu[\sL\disj\sR]\quad=\quad1-\nu\{x\}\quad<\quad1\quad=\quad\mu[0,1].
\]
  Thus, $\mu$ violates  the finite additivity equation (\ref{prob.measure}), so it is not a content.
  \eclaimprf

\noindent  So, without loss of generality, we assume that $\nu$ has no atoms in $(0,1)$.
Thus, for any $r\in(0,1)$, and any $\eps>0$, there exists $\del>0$ such that
$\nu(r-\del,r+\del)<\eps$.

Now, let $\{r_n\}_{n=1}^\oo$ be a countable dense subset of $(0,1)$  (for example, the 
set of all rationals in $(0,1)$).  
For all $n\in\Natur$,  since $r_n$ is not an atom, there exists  $\del'_n>0$ such that
$\nu(r_n-\del'_n,r_n+\del'_n)<1/2^n$.
For all $n\in\Natur$,  define the open set $\sO_n$ as follows:
\bitem
\item Let $\sO_1:=(r_1-\del_1,r_1+\del_1)$.
\item  Let $n\geq 2$.  By induction,  suppose that $\sO_n$ has  already  been defined,  and contains $\{r_1,\ldots,r_n\}$.  Let
$m(n):=\min\{n\in\Natur$; \ $r_n\not\in\Cl[\sO_n]\}$.  (Thus, $m(n)\geq n+1$.)
 \ Define $q_n:=r_{m(n)}$, and 
let $\del''_n:=\inf\{|q_n-u|$; \ $u\in\sO_n\}$; then $\del''_n>0$.
Let $\del_n:=\min\{\del'_{m(n)},\del''_n\}$, and define
$\sO_{n+1} :=\sO_n\disj (q_n-\del_n,q_n+\del_n)$.
\eitem
In this way, we obtain an increasing sequence $(\sO_1\subseteq\sO_2\subseteq\cdots)$
of open sets.  Let $\sU:=\Union_{n=1}^\oo \sO_n$; then $\sU$ is an open subset of $[0,1]$.
Furthermore, $\sU$ is dense in $[0,1]$.  To see this, note that, for
any $N\in\Natur$, we must have $\{r_1,\ldots,r_N\}\subset \Cl[\sO_N]$.
Thus, $\{r_n\}_{n=1}^\oo\subset\Cl[\sU]$, so since $\{r_n\}_{n=1}^\oo$ is dense
in $(0,1)$, it follows that $\Cl[\sU]=[0,1]$.

By construction, we can write $\sU$ as a disjoint union of open intervals:
\begin{eqnarray}  \nonumber
\sU&=&\Disj_{n=1}^\oo (q_{n}-\del_{n},q_{n}+\del_{n})\quad=\quad
\Disj_{n=1}^\oo (r_{m(n)}-\del_{n},r_{m(n)}+\del_{n}). \\
\mbox{Thus,}\quad \nonumber
\nu(\sU) &=& \sum_{n=1}^\oo \nu(r_{m(n)}-\del_{n},r_{m(n)}+\del_n)
 \ \ \leq \ \ 
\sum_{n=1}^\oo \nu(r_{m(n)}-\del'_{m(n)},r_{m(n)}+\del'_{m(n)})
\\ & < &  \sum_{n=1}^\oo \frac{1}{2^{m(n)}}
 \ \ \leq \ \  \sum_{m=1}^\oo \frac{1}{2^{m}}
 \ \ = \ \  1.  \label{X:lebesgue.is.not.content2.e1}
\end{eqnarray}
For all $n\in\Natur$, let $\sU_n:=(q_n-\del_n,q_n+\del_n)$, and let
 $\sL_n:=(q_n-\del_n,q_n)$  and
$\sR_n:=(q_n,q_n+\del_n)$ be the left and right ``halves'' of $\sU_n$.
Let $\sQ:=\{q_n\}_{n=1}^\oo$.   Then $\sU=\sL\disj\sQ\disj\sR$, so  
 $\nu[\sU]=\nu[\sL]+\nu[\sQ]+\nu[\sR]=\nu[\sL]+\nu[\sR]$
(because $\nu[\sQ]=0$ because $\sQ$ is a countable set and $\nu$ has no atoms).

\Claim{$\sL$ and $\sR$ are disjoint regular open sets.\label{X:lebesgue.is.not.content2.C2}}
\bclaimprf
Clearly, $\sL$ and $\sR$ are open and disjoint; it remains to show regularity.  
Let $\sI:=\Int[\Cl(\sL)]$; we must show that $\sI=\sL$.
Since $\sL$ is an open subset of $\Cl(\sL)$, we have $\sL\subseteq\sI$;
we must show that $\sI\subseteq\sL$.

First, note that $\Cl(\sL)\subseteq[0,1]\setminus\sR$; thus,
$\sI\subseteq[0,1]\setminus\sR$.
Since $\sI$ is open, we have $\sI=\Disj_{j=1}^\oo\sI_j$ for some countable collection
$\{\sI_j\}_{j=1}^\oo$ of disjoint open intervals.
For all $j\in\Natur$, let $\sI_j=(a_j,b_j)$.

\subclaim{For all $j\in\Natur$,  $a_j\not\in\sL$ and $b_j\not\in\sL$.\label{X:lebesgue.is.not.content2.C2A}}
\bsubclaimprf
 Since $\{\sI_n\}_{n=1}^\oo$ are disjoint, we must have  $a_j,b_j\not\in\sI_n$ for all $n\in\Natur$.
     Thus, $a_j,b_j\not\in\Disj_{n=1}^\oo\sI_n=\sI$.   
    But $\sL\subseteq\sI$.  Thus, $a_j,b_j\not\in\sL$.
  \esubclaimprf

\noindent Now, fix $j\in\Natur$.   We must have $\sI_j\intsct\sL\neq\emptyset$, because
 $\sL$ is dense in $\Cl[\sL]$, and
$\sI_j$ is an open subset of $\Cl[\sL]$.
Thus, there is some $n\in\Natur$ such that $\sI_j\intsct\sL_n\neq \emptyset$.

\subclaim{$b_j=q_n$.\label{X:lebesgue.is.not.content2.C2B}}
\bsubclaimprf (by contradiction) \ 
Recall that $\sI_j=(a_j,b_j)$ and  $\sL_n:=(q_n-\del_n,q_n)$.  Thus,
if $\sI_j\intsct\sL_n\neq \emptyset$, then $q_n-\del_n<b_j$.   If
$q_n<b_j$, then $\sI_j$ would overlap $\sR_n$, and hence $\sR$,
contradicting the fact that $\sI\subseteq[0,1]\setminus\sR$.
Thus, $q_n-\del_n<b_j\leq q_n$.
But then we must have $b_j=q_n$, by Claim \ref{X:lebesgue.is.not.content2.C2A}.
\esubclaimprf

\subclaim{$a_j=q_n-\del_n$.\label{X:lebesgue.is.not.content2.C2C}}
\bsubclaimprf (by contradiction) \ 
If $a_j>q_n-\del_n$, then $a_j\in\sL_n$  (by Claim \ref{X:lebesgue.is.not.content2.C2B}), contradicting Claim \ref{X:lebesgue.is.not.content2.C2A}.  
Thus, $a_j\leq q_n-\del_n$.  On the other hand,
if $a_j<q_n-\del_n$, then the open interval $(a_j,q_n-\del_n)$
must intersect $\sU$ (because $\sU$ is dense  in $[0,1]$),
which means it must intersect  $\sU_m=(q_m-\del_m,q_m+\del_m)$
for some $m\in\Natur\setminus\{n\}$.  This means that $a_j<q_m+\del_m$.

Recall that $\sU_m=\sL_m\disj\{q_m\}\disj\sR_m$,
where  $\sR_m=(q_m,q_m+\del_m)$.   Also  recall that $\sI_j=(a_j,b_j)$.
If $q_m\leq a_j<q_m+\del_m$, then clearly  $\sI_j$ overlaps
 $\sR_m$, which contradicts the fact that 
$\sI\subseteq[0,1]\setminus\sR$.  Thus, we must have $a_j<q_m$.

Meanwhile, if $q_n-\del_n<q_m-\del_m$, then $(a_j,q_n-\del_n)$ is disjoint from $\sU_m$, contradicting
our assumption that they overlap.  Thus, we must have  $q_m-\del_m\leq q_n-\del_n$.
But if $q_m-\del_m\leq q_n-\del_n<q_m+\del$, then $\sL_n$ and $\sU_m$ would overlap, contradicting the
fact that $\sU_n$ and $\sU_m$ are  disjoint by definition (because $n\neq m$).
Thus, we must also have $ q_m+\del_m\leq  q_n-\del_n$.
Putting it all together, we have $a_j<q_m<q_m+\del_m\leq  q_n-\del_n<q_n=b_j$, where the
last equality is by Claim \ref{X:lebesgue.is.not.content2.C2B}.
This means that $( q_m,q_m+\del_m)\subset (a_j,b_j)$; in other
words, $\sR_m\subset\sI_j$.  But again, this  contradicts the fact that 
$\sI\subseteq[0,1]\setminus\sR$. 

To avoid these contradictions, we must have $a_j=q_n-\del_n$.
\esubclaimprf

\noindent Claims \ref{X:lebesgue.is.not.content2.C2B} and \ref{X:lebesgue.is.not.content2.C2C} together imply that
$\sI_j=(q_n-\del_n,q_n)$;  in other words, $\sI_j=\sL_n$.
This argument works for all $j\in\Natur$. Thus,
every open interval of $\sI$ is actually one of the intervals of $\sL$.
Thus, $\sI\subseteq\sL$.  But  we have already noted that  $\sL\subseteq\sI$.  Thus,
$\sL=\sI$, as desired.  Thus, $\sL$ is regular.  The proof for $\sR$ is similar.
\eclaimprf

\noindent Now, $\sL\disj\sR=\sU\setminus\sQ$, which is dense in $\sU$ (because $\sQ$ is countable  and $\sU$ is open).  But $\sU$ is dense in $[0,1]$.  Thus, $\sL\disj\sR$ is dense in $[0,1]$.  Thus,  $\sL\vee\sR=[0,1]$.
 But \[
 \mu[\sL]+\mu[\sR] \ = \ \nu[\sL]+\nu[\sR] \ = \ 
 \nu[\sL\disj\sR] \ \leq \ \nu[\sU] \ \ \lt{(*)} \ \ 1 \ = \ \mu[0,1] \ = \ \mu[\sL\vee\sR],
 \]
 where $(*)$ is by inequality (\ref{X:lebesgue.is.not.content2.e1}).
 Thus, $\mu$ violates  the finite additivity equation (\ref{prob.measure}), so it is not a content.

 \item (Case 2) \  Now suppose $\sS$ is {\em locally} projectible.  Let $\nu$ be a  Borel probability measure
 on $\sS$.   Since $\sS$ is locally projectible, there exists some regular open subset $\sS_0\subseteq\sS$ such that
 $\sS_0$ is projectible and  $\nu[\sS_0]>0$.  For all Borel subsets $\sB\subseteq\sS_0$,
 define $\nu_0[\sB]:=\nu[\sB]/\nu(\sS_0)$; then
  $\nu_0$ is a Borel probability measure on  $\sS_0$. 
By contradiction, suppose  $\nu$ defines a content $\mu$ when restricted to $\gR(\sS)$. 
 Since $\sS_0\in\gR(\sS)$, we have $\mu[\sS_0]=\nu[\sS_0]>0$, and
  $\gR(\sS_0)=\{\sR\in\gR(\sS)$; \ $\sR\subseteq\sR_0\}$.
  Thus, for all $\sR\in\gR(\sS_0)$,  we can define $\mu_0[\sR]:=\mu[\sR]/\mu[\sS_0]$, to
  obtain a content on $\gR(\sS_0)$.   Equivalently, $\mu_0[\sR]=\nu_0[\sR]$ for all
  $\sR\in\gR(\sS_0)$;  in other words,  $\nu_0$  defines a content when restricted to $\gR(\sS_0)$.  
  But this contradicts {\em Case 1}, because $\sS_0$ is projectible.
  \ethmprf

\section{Compactification representations \label{S:liminal.compactification}}\setcounter{equation}{0}

Corollary \ref{content.vs.Riesz.charge3} shows that liminal representations are especially useful on compact
spaces.  This  suggests that we could greatly improve the
representation in Theorem \ref{content.vs.Riesz.charge2} by  compactifying $\sS$.  
 Let $\barsS$ be a {compactification} of $\sS$ (i.e.  a compact Hausdorff space containing
$\sS$ as a dense subspace).     Let 
$\sC_{\barsS}(\sS,\Real):=\{f_{\restr\sS}$; \ $f\in\sC(\barsS,\Real)\}$ ---this is the set
of all continuous functions in $\sC(\sS,\Real)$ which can be continuously extended to
$\barsS$.   If such an extension exists, then it is unique, because $\sS$ is dense in $\barsS$.
For any $f\in\sC_{\barsS}(\sS,\Real)$, let $\barf$ denote its unique extension to $\sC(\barsS,\Real)$.

\example{ Let $\sS$ be a locally compact Hausdorff space.

\item (a) \  Suppose  $\sS$ is {\em not} compact, and let $\sS^*$ be its Alexandroff
 compactification.  For any $f:\sS\into\Real$ and $L\in\Real$, we will write
``$L=\lim_{s\goto\oo}\,f(s)$'' if, for any open neighbourhood $\sO$ around $L$, there is some compact subset
$\sK\subset\sS$ such that $f(\sS\setminus\sK)\subseteq\sO$.    Then
$\D \sC_{\sS^*}(\sS,\Real)=\{f\in\sC(\sS,\Real)$; \ $\D\lim_{s\goto\oo} \; f(s)$ is well-defined$\}$.
If $f\in \sC_{\sS^*}(\sS,\Real)$, then $\barf(\oo)=\D\lim_{s\goto\oo} \; f(s)$.

\item (b) \ Let $\sS$ be a totally bounded metric space.  Then $\sS$ is locally compact if and only if it is {\em locally complete} ---i.e. every point has a neighbourhood within which every Cauchy sequence converges.
Let $\barsS$ be the (metric) completion of $\sS$.  Then $\barsS$ is a compactification of $\sS$ (as a topological space),
and $\sC_{\barsS}(\sS,\Real)$ is the set of {\em uniformly continuous} real-valued functions on $\sS$.  (The same is true if $\sS$ is a totally bounded {\em uniform} space, and $\barsS$ is its (uniform) completion \cite[Theorems 39.10 and 39.13]{Willard}.)

\item (c) \   Let $\bet\sS$ be the Stone-\v{C}ech compactification of $\sS$.
The Stone-\v{C}ech Extension Theorem implies that {\em every} bounded continuous real-valued function on
$\sS$ has a continuous extension to $\bet\sS$.  Thus,  $\sC_{\bet\sS}(\sS,\Real)=\Cb(\sS,\Real)$.

\item (d) \ Let $\overline{\Real}:=[-\oo,\oo]$, with the obvious topology.   Then $\overline{\Real}$ is a  compactification of $\Real$, and
$\sC_{\overline{\Real}}(\Real,\Real)=\{f\in\sC(\Real,\Real)$; \ $\D \lim_{s\goto\oo} f(s)$
and $\D \lim_{s\goto-\oo} f(s)$
 are well-defined$\}$.  

\item (e) \ Example (d) can be generalized as follows.    Let  $(\sK_1\subseteq\sK_2\subseteq\cdots)$ be a {\dfn compact exhaustion} of $\sS$ ---that is, an increasing
sequence of compact subsets of $\sS$, such that every compact subset of $\sS$ is contained in some
$\sK_n$.  An {\dfn end} of $\sS$ is a decreasing sequence $\eps:=(\sO_1\supseteq\sO_2\supseteq\cdots)$,
where for all $n\in\Natur$, $\sO_n$ is a connected component of $\sS\setminus\sK_n$.
Let $\sE(\sS)$ be the set of ends of $\sS$.
(The definition of $\sE(\sS)$ is independent of the exact choice of compact exhaustion.)
  The {\dfn Freudenthal compactification} of $\sS$ is the set $\barsS:=\sS\disj\sE(\sS)$, where every open subset of $\sS$ remains open in $\barsS$, and where, 
for each $\eps\in\sE(\sS)$, if $\eps:=(\sO_1\supseteq\sO_2\supseteq\cdots)$, then the sets
$\{\sO_n\union\{\eps\}\}_{n=1}^\oo$ form a neighbourhood base for $\eps$ in $\barsS$.
See \cite{Peschke90} for more information.
For example, $\sE(\Real)=\{\pm\oo\}$, and the Freudenthal compactification of $\Real$
is $\overline{\Real}$, as defined in Example (d).  (However, if $N\geq 2$,  then $\sE(\Real^N)$ is a singleton, so the Freudenthal compactification of $\Real^N$ is the same as its Alexandroff compactification.)

For any $f\in\sC(\sS,\Real)$,  $\eps\in\sE(\sS)$, and  $L\in\Real$,  write ``$\lim_{s\goto\eps} \, f(s)=L$ if, for any  neighbourhood $\sU$ of $L$, there is some $n\in\Natur$ such that
$f(\sO_n)\subseteq\sU$ (where $\eps:=(\sO_1\supseteq\sO_2\supseteq\cdots)$).  
Then $\sC_{\barsS}(\sS,\Real)=\{f\in\sC(\sS,\Real)$; \ $\D \lim_{s\goto\eps} f(s)$
 is well-defined for every $\eps\in\sE(\sS)\}$.
 If $f\in \sC_{\barsS}(\sS,\Real)$, then $\barf(\eps)=\D\lim_{s\goto\eps} \; f(s)$ for each $\eps\in\sE(\sS)$.

\item (f)  Let $\wp(\sS)$ be the power set of $\sS$.
 A {\dfn proximity} on $\sS$ is a symmetric, reflexive binary relation $\sim$  on $\wp(\sS)$ such that, for all nonempty $\sA,\sB,\sC\subseteq\sS$, we have:  (i) $\emptyset\not\sim \sA$; (ii) $\sA\sim(\sB\union\sC)$
 if and only if $\sA\sim\sB$ or $\sA\sim \sC$; and (iii) If $\sA\not\sim\sB$, then there exist disjoint 
 $\sD,\sE\subseteq\sX$ with $\sA\not\sim(\sS\setminus\sD)$ and $\sB\not\sim(\sS\setminus\sE)$
 \cite[\S40]{Willard}.  
 A proximity $\sim$ is {\dfn compatible} with the topology of $\sS$ if, for any $\sA\subseteq\sS$,
 we have $\Cl(\sA):=\{s\in\sS$; \ $\{s\}\sim\sA\}$.    If $\barsS$ is any Hausdorff compactification of $\sS$, then we get  a compatible proximity on $\sS$ by stipulating that
 $\sA\sim\sB$ if and only if $\Cl_{\barsS}(\sA)\intsct\Cl_{\barsS}(\sB)\neq\emptyset$ 
 (where $\Cl_{\barsS}(\sA)$ is the closure of $\sA$ in $\barsS$).  In fact, {\em every} compatible
  proximity arises from a compactification in this fashion;  thus, there is a bijective correspondence between
 the compactifications of $\sS$ and the compatible proximities \cite[Definition 41.2]{Willard}.
 
 The {\dfn elementary} proximity $\approx$ on $\Real$ is defined by stipulating that $\sA\approx\sB$ if and only if $\Cl(\sA)\intsct\Cl(\sB)\neq\emptyset$  (for any $\sA,\sB\subseteq\Real$).  
 If $\sim$ is a compatible proximity on $\sS$, then a  function $f:\sS\into\Real$ is  {\dfn proximity-preserving} if, for all  $\sA,\sB\subseteq\sS$ such that $\sA\sim\sB$, we have $f(\sA)\approx f(\sB)$.  Every proximity-preserving function is continuous, but not every continuous function is proximity preserving.  
However, if $\sim$ arises from the compactification $\barsS$, then
$\sC_{\barsS}(\sS,\Real)$ is the precisely set of proximity-preserving functions from $\sS$ to $\Real$.
(This follows by combining Corollary 36.20 with Theorems 36.19, 40.10(b) and 41.1 of \cite{Willard}.)
}

Let $\barsS$ be a  compactification of $\sS$.
For any $\sR\in\gR(\sS)$, there is a unique
 $\barsR\in\gR(\barsS)$ such that $\barsR\intsct\sS=\sR$ (see Lemma \ref{stone.cech.boolean.isomorphism} below).    We will refer
to $\barsR$ as the {\dfn extension} of $\sR$.

\Theorem{\label{content.vs.Riesz.charge4}}
{
Let $\sS$ be a locally compact Hausdorff space.
Let $\mu$ be a content on $\gR(\sS)$.    Let $\barsS$ be a  compactification of $\sS$.  There is a unique  normal  Borel  measure $\barnu$ on $\barsS$, and
 a unique liminal density  structure $\{\barphi_\sR\}_{\sR\in\gR(\barsS)}$ which is subordinate to $\barnu$, such that for any $\sR\in\gR(\sS)$, we have
 \beqn
\label{content.vs.Riesz.stone.cech.1}
\mu[\sR]\quad=\quad \barnu(\barsR)+\int_{\partial\barsR}  \barphi_\barsR \ \mathrm{d}\barnu_{\partial\barsR},
\eeqn
where $\barsR$ is the unique extension of $\sR$ to $\barsS$.
Furthermore,  for any $f\in\sC_{\barsS}(\sS,\Real)$, we have
\beqn
\label{content.vs.Riesz.stone.cech.2}
  \Integral_{{\sR}} f \dmu \quad=\quad  \int_{\barsR} \barf \ \mathrm{d}\barnu + \int_{\partial\barsR} \barf \, \barphi_{\barsR} \ \mathrm{d}\barnu_{\partial\barsR},
\eeqn
where $\barf$ is the unique extension of $f$ to $\sC(\barsS,\Real)$.
}

\example{\label{X:content.vs.Riesz.charge4}
  (a) If $\barsS$ is the Alexandroff compactification of $\sS$, then  (\ref{content.vs.Riesz.stone.cech.2}) holds
  for any $f\in\sC(\sS,\Real)$ such that $\D\lim_{s\goto\oo} \, f(s)$ is well-defined. 

  (b) \  If $\sS$ is a totally bounded, locally complete metric space, and
 $\barsS$ is its completion, then (\ref{content.vs.Riesz.stone.cech.2}) holds  for all uniformly continuous  $f\in\sC(\sS,\Real)$.    

(c)  If $\barsS$ is the Stone-\v{C}ech compactification of $\sS$, then
  (\ref{content.vs.Riesz.stone.cech.2}) holds  for all  $f\in\Cb(\sS,\Real)$.   
}

The proof of Theorem \ref{content.vs.Riesz.charge4} depends on the following lemma, which
 establishes that each element of  $\gR(\sS)$ has a unique  ``extension'' to an element of $\gR(\barsS)$. 

\Lemma{\label{stone.cech.boolean.isomorphism}}
{
 Let $\barsS$ be a compactification of a  locally compact Hausdorff  space $\sS$.
 If $\barsR\in\gR(\barsS)$, then $\barsR\intsct\sS\in\gR(\sS)$.  Furthermore, the
function $\barsR\mapsto \barsR\intsct\sS$ is a Boolean algebra isomorphism from
$\gR(\barsS)$ to $\gR(\sS)$.
 
}
\bthmprf
{\em Homomorphism.} \ Let $\iota:\sS\into\barsS$ be the inclusion map.  This is a continuous function, because $\sS$ has the subspace topology it inherits from $\barsS$.  It is also an {\em open} function, because
$\sS$ is an open subset of $\barsS$ \cite[Theorem 3.5.8]{Engelking}.
Thus $\iota^{-1}:\gR(\barsS)\into\gR(\sS)$  is
a Boolean algebra homomorphism \cite[4A2B (f)(iii)]{Fremlin4II} (or Lemma \ref{measurability.lemma}(b)).  But for any $\barsR\in\gR(\barsS)$, we have
$\iota^{-1}(\barsR)=\sS\intsct\barsR$.

{\em Injective.} \ Let $\barsR,\barsQ\in\gR(\barsS)$.   Suppose
$\barsR\neq\barsQ$; then either $\barsR\intsct(\neg \barsQ)\neq \emptyset$ 
or $\barsQ\intsct(\neg \barsR)\neq \emptyset$. Suppose the former.
Then $\barsR\intsct(\neg \barsQ)$ is a nonempty open subset of $\barsS$. 
But $\sS$ is dense in $\barsS$.  Thus, $\sS\intsct \barsR\intsct(\neg \barsQ)\neq\emptyset$,
which implies that $\sS\intsct\barsR\neq \sS\intsct\barsQ$.

{\em Surjective.}\ Let $\sR\in\gR(\sS)$.  Define $\barsR:=\Int_{\barsS}\lb[\Cl_{\barsS}(\sR)\rb]$.  
Let $\gO$ be the topology of $\sS$, and let $\bargO$ be the topology of $\barsS$.   Now,  $\sR\subseteq\Cl_{\barsS}(\sR)$,  and  $\sR$ is open in $\barsS$,  so  $\sR\subseteq\Int_{\barsS}\lb[\Cl_{\barsS}(\sR)\rb]=\barsR$. Thus, $\sR\subseteq\sS\intsct\barsR$. Conversely, 
\beq
\barsR &=& \Union\set{\barsO\in\bargO}{\barsO\subseteq\Cl_{\barsS}(\sR)}.\\
\mbox{Thus,}\quad \barsR\intsct\sS &=&
\Union\set{\sS\intsct\barsO}{\barsO\in\bargO \And \barsO\subseteq\Cl_{\barsS}(\sR)}\\
&\sssubset{(*)}&
\Union\set{\sO}{\sO\in\gO \And \sO\subseteq\sS\intsct\Cl_{\barsS}(\sR)}\\
&\eeequals{(\dagger)}&
\Union\set{\sO}{\sO\in\gO \And \sO\subseteq\Cl_{\sS}(\sR)}
\quad=\quad \Int_{\sS}\lb[\Cl_{\sS}(\sR)\rb]
\quad\eeequals{(\diamond)}\quad\sR. 
\eeq
Here, $(*)$ is because $\gO=\{\sS\intsct\barsO$; \ $\barsO\in\bargO\}$, while 
 $(\dagger)$ is because $\sS\intsct\Cl_{\barsS}(\sR)=\Cl_{\sS}(\sR)$
because $\sR\subseteq\sS$ and $\sS$ has the subspace topology.
Finally, $(\diamond)$ is because $\sR\in\gR(\sS)$.  This shows that $\sR=\sS\intsct\barsR$.
\ethmprf

\paragraph{Remark.}   In the proof of Lemma  \ref{stone.cech.boolean.isomorphism} , it is crucially important
that $\sS$ be an {\em open} dense subset of $\barsS$.
Let's say that a compactification is {\dfn proper} if it has this property.
 If $\sS$ is a Hausdorff space, then the following are equivalent:
(1) $\sS$ is locally compact; \ 
(2) $\sS$ has a proper compactification; \ 
(3) {\em Every} compactification of $\sS$ is proper \cite[Theorem 3.5.8]{Engelking}.
For this reason, Theorem \ref{content.vs.Riesz.charge4} only applies to
locally compact Hausdorff spaces. 

\newcommand{\dbarmu}{\ \mathrm{d}\barmu}
\newcommand{\barsP}{\overline{\sP}}

\bthmprf[Proof of Theorem \ref{content.vs.Riesz.charge4}.]
Let $\mu$ be a content on $\gR(\sS)$.  Define the content
$\barmu$ on $\gR(\barsS)$ by setting
$\barmu(\barsR):=\mu[\sS\intsct\barsR]$ for all $\barsR\in\gR(\barsS)$.
Lemma  \ref{stone.cech.boolean.isomorphism}  implies that this is a well-defined content.
For any $\sR\in\gR(\sS)$, if $\barsR$ is the (unique) element  of $\gR(\barsS)$  such that
$\sS\intsct\barsR=\sR$, then we have 
 
\beqn
\label{content.vs.Riesz.charge4.e1}
\mu[\sR]\quad=\quad\barmu[\barsR].
\eeqn

\Claim{\label{stone.cech.integrator.defn}
For all $\sR\in\gR(\sS)$ and all $f\in\sC_{\barsS}(\sS,\Real)$, if $\barf$ is the unique extension of $f$ to $\sC(\barsS)$, then   
\[
\Integral_{\sR} f \dmu\quad=\quad \Integral_{\barsR}\barf \dbarmu.
\]
}
\bclaimprf
First suppose $f$ is a simple function; then there is a regular partition $\{\sP_1,\ldots,\sP_N\}$ of $\sS$ 
and  some $r_1,\ldots,r_N\in\Real$ such that $f=\sum_{n=1}^N r_n\,\bone_{\sP_n}$.
It follows that
$\barf=\sum_{n=1}^N r_n\,\bone_{\barsP_n}$.  Thus,
\beqn
\label{content.vs.Riesz.charge4.e2}
\Integral_{\sR} f \dmu
\quad\eeequals{(*)}\quad
\sum_{n=1}^N r_n\,\mu[\sR\intsct\sP_n]
\quad\eeequals{(\dagger)}\quad
\sum_{n=1}^N r_n\,\barmu[\barsR\intsct\barsP_n]
\quad\eeequals{(*)}\quad
\Integral_{\barsR} \barf \dbarmu.
\eeqn
Here, both $(*)$ are by  equation (\ref{integral.of.simple.function})  from Theorem \ref{from.probability.to.expectation0} , while $(\dagger)$ is by
equation (\ref{content.vs.Riesz.charge4.e1}), because $\barsR\intsct\barsP_n=\overline{\sR\intsct\sP_n}$ by
  Lemma  \ref{stone.cech.boolean.isomorphism} .

Now let $f\in\sC_{\barsS}(\sS,\Real)$ be arbitrary.  Let $\{f_n\}_{n=1}^\oo$ be a sequence of simple functions
converging uniformly to $f$.  Then
\beqn
\label{content.vs.Riesz.charge4.e3}
\lim_{n\goto\oo} \norm{\barf_n-\barf}{\oo}
\quad\eeequals{(*)}\quad
\lim_{n\goto\oo} \norm{f_n-f}{\oo}
\quad=\quad
0,
\eeqn
where $(*)$ is because the  transformation $f\mapsto\barf$ is a  norm-preserving, linear isomorphism from $\sC_{\barsS}(\sS,\Real)$ to $\sC(\barsS,\Real)$.
Thus,
\[
\Integral_{\sR} f \dmu
\quad\eeequals{(*)}\quad
\lim_{n\goto\oo} \Integral_{\sR} f_n \dmu
\quad\eeequals{(\dagger)}\quad
\lim_{n\goto\oo} \Integral_{\barsR} \barf_n \dbarmu
\quad\eeequals{(*)}\quad
 \Integral_{\barsR} \barf \dbarmu,
\]
as claimed.  Here, both $(*)$ are by equation (\ref{content.vs.Riesz.charge4.e3}), because the functions $\integral_\sR \dmu$ and 
$\integral_{\barsR} \dbarmu$ are continuous with respect to the uniform norm.
Meanwhile, $(\dagger)$ is by
applying equation (\ref{content.vs.Riesz.charge4.e2}) for all $n\in\Natur$.
\eclaimprf

\noindent
Corollary \ref{content.vs.Riesz.charge3} yields a
unique  normal  Borel  probability measure $\barnu$ on $\barsS$ and  liminal density  structure $\{\barphi_\barsR\}_{\barsR\in\gR(\barsS)}$ which is subordinate to $\barnu$ and which satisfies equations {\rm(\ref{content.vs.Riesz.charge.partb2})} and
{\rm(\ref{content.vs.Riesz.charge.partb})}
for all $\barsR\in\gR(\barsS)$ and all $f\in\sC(\barsS,\Real)$.
Combine this with equation (\ref{content.vs.Riesz.charge4.e1}) and Claim \ref{stone.cech.integrator.defn} to  obtain equations (\ref{content.vs.Riesz.stone.cech.1}) and (\ref{content.vs.Riesz.stone.cech.2}).  
\ethmprf

\section{Integration via Stone spaces\label{S:gleason}}
\setcounter{equation}{0}

In this section, we will introduce an entirely different representation of integration operations, which
is applicable to a content defined on a Boolean subalgebra of regular open sets of any locally compact Hausdorff space $\sS$.
For such a Boolean algebra $\gB$, we construct a compact, totally disconnected space $\sS^*$, along with a continuous surjection from $\sS^*$ onto the Stone-\v{C}ech compactification of $\sS$.
Thus, any bounded, continuous  function $g$ on $\sS$ corresponds to a continuous function $g^*$ on $\sS^*$.
Every element of $\gB$ corresponds to a clopen subset of $\sS^*$;
we use this to show that every content $\mu$ on $\gB$ corresponds to a Borel probability measure $\mu^*$ on $\sS^*$.
Finally, we show that the $\mu$-integral of any function $g$ is equal to the $\mu^*$-integral of $g^*$. 
Thus, it is possible to representation $\mu$ entirely in terms of a Borel probability measure, {\em without} a liminal structure ---but at at the cost of working in an extension of the original space $\sS$.

   For any topological space $\sT$, let $\clop(\sT)$ be the collection of all clopen subsets of $\sT$.   This collection is a Boolean algebra under the standard set-theoretic operations of union, intersection,
and complementation. 
A {\dfn Stonean space} is a compact, totally disconnected Hausdorff space.
For any Boolean algebra $\gB$, let $\sigma(\gB)$
be the set of all ultrafilters of $\gB$.  For any $\sB\in\gB$, let $\sB^*:=\{\gU\in\sigma(\gB)$; \ $\sB\in\gU\}$.
The collection $\{\sB^*\}_{\sB\in\gB}$ is a base of clopen sets for a topology on $\sigma(\gB)$, making $\sigma(\gB)$
into Stonean space.  This is called the {\dfn Stone space} of $\gB$.    The Stone Representation Theorem
says that there is a Boolean algebra isomorphism from $\gB$ to $\clop[\sigma(\gB)]$ given by $\gB\ni \sB\mapsto \sB^*\in\clop[\sigma(\gB)]$.\footnote{See e.g. \cite[Theorem 2.10, p.51]{Walker}, \cite[\S3.2]{PorterWoods}, or \cite[Sections 311E, 311F and 311I]{Fremlin}.}

If $\gA$ is another Boolean algebra,
and $h:\gA\into\gB$ is a Boolean algebra homomorphism, then we obtain a continuous function
$H:\sigma(\gB)\into\sigma(\gA)$ which maps each ultrafilter in $\gB$ to its $h$-preimage ultrafilter in $\gA$.
This yields a contravariant functor $\sigma$ from
the category of Boolean algebras to the category of Stonean spaces and continuous functions. In fact,
the Stone Duality Theorem says that $\sigma$ is  a functorial isomorphism between these  categories.\footnote{See e.g. 
\cite[Example 10.8(a), p. 247]{Walker}, \cite[\S4.1, p.71]{Johnstone}, or \cite[Sections 312Q  and 312R]{Fremlin}.}

Now, let $\sS$ be a locally compact Hausdorff space, and
let $\hsS$ be its Stone-\v{C}ech compactification.  Recall from
Lemma  \ref{stone.cech.boolean.isomorphism}  that there is a Boolean algebra isomorphism $h:\gR(\hsS) \, \widetilde{\longrightarrow}\,\gR(\sS)$ given by $h(\hsR):=\hsR\intsct\sS$, for all $\hsR\in\gR(\hsS)$. 
For any $\sR\in\gR(\sS)$, let $\hsR:=h^{-1}(\sR)$ ---the unique element of $\gR(\hsS)$ such that $\hsR\intsct\sS=\sR$.
 Let $\gB\subseteq\gR(\sS)$ be any Boolean subalgebra, and let
  $\hgB:=\{\hsB$; \ $\sB\in\gB\}$; this is a Boolean subalgebra of $\gR(\hsS)$, and is
isomorphic to $\gB$ via $h$.   We  say that $\gB$ is {\dfn generative} if $\hgB$ is a base for the topology of $\hsS$.

\Lemma{\label{generative.lemma}}
{
If $\gB$ is generative, then it is a base for the topology of $\sS$.
}
\bthmprf
It suffices to show that every open subset of $\sS$ contains a $\gB$-neighbourhood around each of its points.
So, let $\sO\subseteq\sS$ be open, and let $s\in\sO$.  Then $\sO$ is also  open in
 $\hsS$ (because $\sS$ is  an open subset of $\hsS$).
Thus, there exists  $\hsB\in\hgB$ with $s\in\hsB\subseteq\sO$ (because $\hgB$ is a base
for the topology of $\hsS$).   Note that $\hsB =\hsB\intsct\sS$ (because $\hsB\subseteq\sS$); thus $\hsB\in\gB$.   This works for any $s$ and $\sO$; thus, $\gB$
is a base for the topology of $\sS$.
\ethmprf

\noindent The converse of Lemma \ref{generative.lemma} is false:
for $\gB$ to be generative, it is {\em not sufficient} that $\gB$ be a base for the topology of $\sS$.
For example, let $\sS=\Natur$, with the discrete topology;  then $\gR(\Natur)=\wp(\Natur)$.
Let $\gF$ be the set of all finite subsets of $\Natur$, and let $\gG:=\{\Natur\setminus\sF$; \ $\sF\in\gF\}$;
then $\gB:=\gF\disj\gG$ is a Boolean  subalgebra   of $\gR(\Natur)$ which generates the topology of $\Natur$,
because it contains all singleton sets.  Let $\hdN$ be the Stone-\v{C}ech compactification of $\Natur$, and
let $\hgG:=\{\hdN\setminus\sF$; \ $\sF\in\gF\}$.  Then $\hgB=\gF\disj\hgG$, which does
{\em not} generate the topology of $\hdN$.
Nevertheless, the full Boolean algebra $\gR(\sS)$ itself is always generative, because of the next lemma
(using the fact that $\widehat{\gR(\sS)}=\gR(\hsS)$).

\Lemma{\label{gleason.covering.C1}}
{
Let $\sS$ be a locally compact Hausdorff space.  Then $\gR(\sS)$ is a base for the topology on $\sS$.
}
\bthmprf
Let $s\in\sS$, and let $\sO\subseteq\sS$ be any open neighbourhood of $s$.  Since $\sS$ is locally compact, there is a compact subset
$\sK\subseteq\sO$ which is also a neighbourhood of $s$.  Let $\sR:=\Int(\sK)$;  then $\sR$ is a regular
open subset of $\sS$, and $s\in\sR\subseteq\sO$, as desired.
\ethmprf

\noindent Let $\gB$ be a subalgebra of $\gR(\sS)$, and let  $\sS^*:=\sig(\gB)$.
 The isomorphism $h:\gR(\hsS)\,\widetilde{\longrightarrow}\,\gR(\sS)$  from Lemma  \ref{stone.cech.boolean.isomorphism}   restricts to an isomorphism $h_\restr:\hgB\,\widetilde{\longrightarrow}\,\gB$.    Via  the Stone Duality Theorem, $h_\restr$ induces a homeomorphism
$H:\sS^*\,\widetilde{\longrightarrow}\,\sig(\hgB)$.  To be precise, for any $s^*\in\sS^*$ (an ultrafilter in $\gB$), we have
\beqn
\label{gleason.covering.defn0}
H(s^*)\quad=\quad \set{\hsB}{\sB\in s^*}.
\eeqn

\begin{figure}
\[
\begin{tikzcd}
\sS^*=\sigma(\gB) \arrow[twoheadrightarrow]{dr}{p} \arrow[Rightarrow]{r}{H} & \sigma(\hgB) \arrow[twoheadrightarrow]{d}{\hat{p}} \\
\sS \arrow[hookrightarrow]{r} & \hsS 
\end{tikzcd}
\qquad\qquad
\begin{tikzcd}
\clop(\sS^*)  &\arrow[Rightarrow,swap]{l}{H^{-1}} \clop[\sigma(\hgB)] \\
\gB\arrow[Rightarrow]{u}{*} & \hgB \arrow[Rightarrow,swap]{u}{*} \arrow[Rightarrow]{l}{h_\restr}
\end{tikzcd}
\]
\caption{\footnotesize \label{fig:generalized.gleason} The four topological spaces and four Boolean algebras in
Proposition \ref{generalized.gleason}.
Double-lined arrows indicate homeomorphisms or Boolean algebra isomorphisms.
``$*$'' denotes the isomorphisms from the Stone Representation Theorem,
$h_\restr$ is the isomorphism from Lemma  \ref{stone.cech.boolean.isomorphism}, and $H$ is the homeomorphism defined by formula (\ref{gleason.covering.defn0}).
Finally, $p$ is the continuous surjection defined by formula (\ref{gleason.covering.defn}).  The surjection $\hat{p}$ is defined by an
intersection formula similar to (\ref{gleason.covering.defn}); it does not appear explicitly in the paper, but
$p=\hat{p}\circ H$.}
\end{figure}

\Proposition{\label{generalized.gleason}}
{
Let $\sS$ be a locally compact Hausdorff space, and
suppose $\gB$ is a generative subalgebra of $\gR(\sS)$.  For any $s^*\in\sS^*$,  the intersection
\beqn
\label{gleason.covering.defn}
\Intsct_{\hsB\in H(s^*)} \Cl(\hsB)\quad=\quad \Intsct_{\sB\in s^*} \Cl(\hsB)
\eeqn
 contains only a single element, which we will denote by $p(s^*)$.  This determines a continuous
 surjective function $p:\sS^*\into\hsS$.
 }
 Figure \ref{fig:generalized.gleason} summarizes the relationships between the
  four topological spaces and the four Boolean algebras involved in Proposition \ref{generalized.gleason}. 
 If $\sS$ itself is compact, then $\hsS=\sS$, so that $p$ is a continuous surjection from $\sS^*$ into $\sS$.
In the special case when $\gB=\gR(\sS)$, the pair $(\sS^*,p)$  is called the {\dfn Gleason cover} (or the {\dfn projective cover}, or the {\dfn absolute}) of the space $\sS$.\footnote{See e.g. \cite[\S10.54, p.288]{Walker}, \cite[\S3.10, p.107]{Johnstone}, or \cite[Chap. 6]{PorterWoods}.  We thank Vincenzo Marra for introducing us to Gleason covers.}

\bthmprf[Proof of Proposition \ref{generalized.gleason}.]
    Each  $s^*$  in $\sS^*$ is a  filter in $\gB$.
Thus, by defining formula (\ref{gleason.covering.defn0})  and Lemma  \ref{stone.cech.boolean.isomorphism} ,
$H(s^*)$ is a  filter in $\hgB$ ---hence, a filter of subsets of $\hsS$.
Thus, the collection $\{ \Cl(\hsB)$; \ $\hsB\in H(s^*)\}$ is a filter of {\em closed} subsets of $\hsS$,
so it satisfies the Finite Intersection Property.  Thus, the intersection (\ref{gleason.covering.defn}) is nonempty, because
$\hsS$ is compact.  

To see that  (\ref{gleason.covering.defn}) is a singleton, let $\hs_1$ be some element of  (\ref{gleason.covering.defn}), and let $\hs_2$ be any other element of $\hsS$.  
There exists a disjoint open sets $\hsO_1,\hsO_2\subset\hsS$  with $\hs_1\in\hsO_1$ and $\hs_2\in\hsO_2$ (because $\hsS$ is Hausdorff).
Then there exists  $\hsB\in\hgB$ such that $\hs_1\in\hsB\subseteq\hsO_1$ (because $\hgB$ generates
the topology of $\hsS$, because $\gB$ is generative).   Let $\sB:=\hsB\intsct\sS$;  then $\sB\in\gB$.

 \claim{$\sB\in s^*$.}
 \bclaimprf  (by contradiction) Suppose $\sB\not\in s^*$. Then
$\neg\sB\in s^*$, because $s^*$ is an ultrafilter.  But $\hs_1\in\hsB$, so $\hs_1\not\in \Cl(\neg\hsB)=\Cl(\widehat{\neg\sB})$;
thus, $\hs_1$ is {\em not} the intersection (\ref{gleason.covering.defn}), which is a contradiction.  To avoid this contradiction, we must have $\sB\in s^*$.\eclaimprf

\noindent Now,  $\hsB\subseteq \hsO_1$, so
 $\hsB$ is disjoint from $\hsO_2$.
Thus, $\Cl(\hsB)$ is also disjoint from $\hsO_2$ (because $\hsO_2$ is open).  In particular, $\hs_2\not\in\Cl(\hsB)$.
 Thus, $\hs_2$ is {\em not} an element of the intersection (\ref{gleason.covering.defn}).   This holds for all $\hs_2\neq \hs_1$,
so we conclude that $\hs_1$ is the {\em only} element of  (\ref{gleason.covering.defn}).
Thus, the function $p$ is well-defined.

\item {\em $p$ is surjective.}  \
Let $\hs\in\hsS$.   Let $\hgB_{\hs}$ be the set of all elements in $\hgB$ containing $\hs$.
This is clearly a filter.  Thus, by the Ultrafilter Theorem, it can be completed to an ultrafilter $\gU$ ---i.e. an element of $\sig(\hgB)$.  
Since $H:\sS^*\into\sig(\hgB)$ is bijective, there is a unique  $s^*\in\sS^*$ such that $H(s^*)=\gU$.
For any $\hs'\in\hsS\setminus\{\hs\}$, there are disjoint open sets $\hsO,\hsO'\subseteq\hsS$ such that $\hs\in\hsO$ and
$\hs'\in\hsO'$ (because $\hsS$ is Hausdorff).  Since $\hgB$ is a base for the topology of $\hsS$,
there exists  $\hsB\in\hgB_{\hs}$ such that $\hs\in\hsB\subseteq\hsO$.  Now, $\Cl(\hsB)$ is disjoint from $\hsO'$ (because $\hsO'$ is open); thus
 $\hs'\not\in\Cl(\hsB)$.  Thus, $\hs'$ is not an element of  (\ref{gleason.covering.defn}).  This holds for all $\hs'\in\hsS\setminus\{\hs\}$.  But we have already established that (\ref{gleason.covering.defn}) is nonempty;
 thus, we must have  $p(s^*)=\hs$.

\item {\em $p$ is continuous.} \   Let $\hs\in\hsS$ and let $\hsO\subseteq\hsS$ be an open neighbourhood of $\hs$.
Now $\hsS$ is compact Hausdorff, hence locally compact;  thus, there is a compact subset $\hsK\subseteq\hsS$  such that $\hsK\subseteq\hsO$ and
$\hs\in\hsO_1$, where  $\hsO_1:=\Int(\hsK)$.  
There exists $\hsB\in\hgB$ such that $\hs\in\hsB\subseteq\hsO_1$ (because $\hgB$ is a base for
the topology of $\hsS$).  Let $\sB^\sig:=\{\gU\in\sigma(\hgB)$; \ $\hsB\in\gU\}$; as noted above, this is one of
the elements in the clopen basis for the topology on $\sigma(\hgB)$.  Thus, if we define
$\sB^*:=H^{-1}(\sB^\sig)$, then $\sB^*$ is an open subset of $\sS^*$ (because
$H$ is continuous).   

\Claim{$p^{-1}\{\hs\}\subseteq\sB^*$.}
\bclaimprf
Let $s^*\in p^{-1}\{\hs\}$.   Then $p(s^*)=\hs$, which, by  formulae   (\ref{gleason.covering.defn0}) and   (\ref{gleason.covering.defn}), means that
$\hs$ is  in $\Cl(\hsR)$ for all $\hsR\in  H(s^*)$.  Now, if $\hsB\not\in H(s^*)$, then
we must have $\neg\hsB\in H(s^*)$ (because $ H(s^*)$ is an ultrafilter in $\hgB$).
But clearly $\hs\not\in\Cl(\neg\hsB)$ (because $\hs\in\hsB$), so this is a contradiction.  Thus, 
$\hsB\in H(s^*)$.  Thus, $ H(s^*)\in \sB^\sig$, so $s^*\in\sB^*$.
\eclaimprf

\Claim{$p(\sB^*)\subseteq\hsO$.}
\bclaimprf
For any $s^*\in\sB^*$, we have $H(s^*)\in\sB^\sig$, which means
$\hsB\in H(s^*)$. Thus, equations (\ref{gleason.covering.defn0}) and (\ref{gleason.covering.defn}) together imply that
$p(s^*)\in\Cl(\hsB)$. But $\sB\subseteq \hsO_1$, so 
$\Cl(\sB)\subseteq\Cl( \hsO_1)\subseteq\hsK\subseteq \hsO$.  Thus, $p(s^*)\in \hsO$ for all $s^*\in\sB^*$.
\eclaimprf

\noindent Claims 2 and 3 show that, for any open neighbourhood $ \hsO\subseteq\hsS$ around $ \hs$ , there is some open neighbourhood
$\sB^*$ around $p^{-1}\{ \hs\}$ such that $p(\sB^*)\subseteq \hsO$.  Thus, $p$ is continuous at each point in $p^{-1}\{ \hs\}$.
This argument holds for all $ \hs\in\hsS$;  thus, $p$ is continuous on $\sS^*$.
\ethmprf

\begin{figure}
\[
\begin{tikzcd}
  ~ & \sC(\sS^*,\Real) & \\
  \sC_b(\sS,\Real) \arrow{ur}{\scriptstyle  g\mapsto g^*}
  \arrow[Rightarrow,bend left=10]{rr}{\scriptscriptstyle \mathrm{extension} \ g\mapsto \hg }&  & \sC(\hsS,\Real) \arrow[swap]{ul}{\scriptstyle \hg\mapsto \hg\circ p}
  \arrow[Rightarrow,bend left=10]{ll}{\scriptscriptstyle \mathrm{restriction} \ \hg\mapsto \hg_{|\sS}}
\end{tikzcd}
\]  
\caption{\footnotesize \label{gleason.covering.fig} The three vector spaces of continuous functions involved in Theorem 
\ref{gleason.covering}.}
\end{figure}

\noindent
We will use the construction from Proposition \ref{generalized.gleason} to give a new representation of integration
with respect to a content.  For any $g\in\Cb(\sS,\Real)$,  the Stone-\v{C}ech Extension Theorem yields 
a unique extension  $ \hg\in\sC(\hsS,\Real)$ such that $\hg_{\restr\sS}=g$.
Let $\gB$ be a generative Boolean subalgebra of $\gR(\sS)$, and let $\sS^*$ and
 $p:\sS^*\into\hsS$ be as in Proposition \ref{generalized.gleason}.
We then define $g^*:=  \hg\circ p:\sS^*\into\Real$.  The transformation $g\mapsto g^*$ is a bounded linear function
from $\Cb(\sS,\Real)$ into $\sC(\sS^*,\Real)$. (See Figure \ref{gleason.covering.fig}.) 
 Now let $\mu$ be a content on $\gB$.  We can then define 
 a probability charge $\mu^*$ on  $\clop\lb(\sS^*\rb)$  as follows:
\beqn
\label{mu.star.defn}
\mbox{for all $\sB\in\gB$,}\quad
\mu^*[\sB^*] \ := \ \mu[\sB],
\quad\mbox{where $\sB^*:=\{s^*\in\sS^*$; \ $\sB\in s^*\}$.}
\eeqn
(Recall that the Stone Representation Theorem says that the map $\sB\mapsto\sB^*$ is a Boolean algebra isomorphism from $\gB$ to $\clop\lb(\sS^*\rb)$.)
 We now have everything we need for the next result.
Recall the definitions of $\GB(\sS)$ and  $\integral \dmu$  from
Section \ref{S:integration}.

\Theorem{\label{gleason.covering}}
{
Let $\sS$ be a locally compact Hausdorff space, let $\gB$ be a generative Boolean subalgebra of $\gR(\sS)$, 
let $\mu$ be a content on $\gB$, and let $\mu^*$ be the corresponding probability charge on $\clop\lb(\sS^*\rb)$.
Then $\mu^*$ can be extended to a unique Borel probability measure on $\sS^*$.
Furthermore, for any $g\in\GB(\sS)$ and any $\sD\in\gB$,  we have
\beqn
\label{gleason.covering.e1}
 \Integral_{\sD} g \dmu  \quad=\quad\int_{\sD^*} g^* \dmu^*.
\eeqn
}
\noindent  This representation theorem has two advantages over the representations from Sections 
\ref{S:liminal} and \ref{S:liminal.compactification}.  First, it applies to a content defined on {\em any} generative Boolean subalgebra $\gB$ of $\gR(\sS)$.
Second the representation only requires a single Borel measure $\mu^*$, not an entire liminal structure.
The disadvantage is that $\mu^*$ is defined on $\sS^*$, a topological space even larger and more exotic than
the     Stone-\v{C}ech compactification $\hsS$ (as shown by Proposition \ref{generalized.gleason}).
 The proof of Theorem \ref{gleason.covering} depends on two lemmas.
 The first one is a straightforward variant of the   Kolmogorov Consistency Theorem, but for completeness, we provide a proof.

\Lemma{\label{stonean.measure.extension}}
{
Let $\sT$ be any Stonean space, and let $\mu$ be a probability
charge on $\clop(\sT)$.  Then there is a unique Borel probability measure $\nu$ on $\sT$ which extends $\mu$.
}
\bthmprf 
Let $\sF_0:=\{\bone_\sB$; \ $\sB\in\clop(\sT)\}$, and
 let $\sF$ be the set of finite linear
combinations of elements in $\sF_0$.  Then $\sF$ is an algebra of continuous functions on $\sT$,
because the sum or product of any two elements of $\sF$ is also an element of $\sF$ (because $\clop(\sT)$
is a Boolean algebra under standard set-theoretic operations).

\Claim{$\sF$ separates points.}
\bclaimprf
For any distinct $s,t\in\sT$, there is an open neighbourhood around $s$ which excludes $t$ (because $\sT$ is Hausdorff).
Thus, there is a clopen neighbourhood $\sB$ around $s$ which excludes $t$ (because $\sT$ is totally disconnected).
Thus, $\bone_\sB$ separates $s$ from $t$.
\eclaimprf

\noindent Combining Claim 1 with the Stone-Weierstrass theorem, we deduce that $\sF$ is dense in $\sC(\sT,\Real)$ in the uniform topology.
Define  $\dI^\mu:\sF\into\Real$ by
\[
\dI^\mu\lb(\sum_{n=1}^N r_n \,\bone_{\sB_n}\rb) \quad:=\quad \sum_{n=1}^N r_n\, \mu[\sB_n] ,
\qquad\mbox{for any} \ \sum_{n=1}^N r_n \,\bone_{\sB_n}\in \sF.
\]
Then $\dI$ is a bounded linear functional on $\sF$.  But $\sF$ is dense in $\sC(\sT,\Real)$, so $\dI^\mu$ extends to
a unique bounded linear functional $\dI^\mu:\sC(\sT,\Real)\into\Real$.
Since $\sT$ is compact Hausdorff, the Riesz Representation Theorem yields a unique Borel probability measure
$\nu$ such that $\dI^\mu(f)=\int_\sT f\dnu$ for all $f\in\sC(\sT,\Real)$.  In particular, for any $\sB\in\clop(\sT)$,
we have $\mu[\sB] =\dI_\mu[\bone_\sB]=\int_\sT \bone_\sB\dnu=\nu[\sB]$, so $\nu$ extends $\mu$, as claimed.
\ethmprf

\Lemma{\label{gleason.covering.map.lemma}}
{
 Let $\gB$, $\hgB$,  $\sS^*$ and $\hsS$ be as in Proposition \ref{generalized.gleason}. 
Let $s^*\in\sS^*$ and  let $\hs\in\hsS$.  Then $p(s^*)=\hs$ if and only if, for
every $\hsB\in\hgB$ with $\hs\in\hsB$, we have $(\sS\intsct\hsB)\in s^*$.}
\bthmprf
For any $\sB\in\gB$, let $\hsB$ denote the (unique) element of $\hgB$ such that $\sS\intsct\hsB=\sB$.

``$\implies$'' \ 
By defining formula (\ref{gleason.covering.defn}), $p(s^*)=\hs$ if $\hs\in\Cl(\hsQ)$ for all $\sQ\in s^*$
---or equivalently, if every open neighbourhood of $\hs$ overlaps $\hsQ$ for every  $\sQ\in s^*$.
In particular, if $\hsB\in\hgB$ and $\hs\in\hsB$, then $\hsB$ overlaps $\hsQ$ for every  $\sQ\in s^*$.
If $\sB:=\sS\intsct\hsB$, then $\neg\sB=\sS\intsct(\neg\hsB)$;
thus,  $\neg\sB$ is {\em not} an element of $s^*$ (because $(\neg\hsB)\intsct\hsB=\emptyset$).  Since
$s^*$ is an {\em ultra}filter, this means that $\sB$ itself must be an element of $s^*$, as claimed.
This argument holds for any  $\hsB\in\hgB$ with $\hs\in\hsB$.

``$\seilpmi$'' \ 
Let $\hsO\subseteq\hsS$ be any open neighbourhood of $\hs$.  Since $\hgB$ generates the topology on $\hsS$,
there is some $\hsB\in\hgB$ such that  $\hs\in\hsB\subseteq\hsO$;  thus, if $\sB:=\sS\intsct\hsB$, then
$\sB\in s^*$.   
Now, let $\sQ\in s^*$ be arbitrary.  Then $\sQ\intsct\sB\neq\emptyset$ (because $s^*$ is a filter).
Thus, $\sQ\intsct\hsO\neq\emptyset$ (because $\sB\subseteq\hsB\subset\hsO$).
This shows that every element of $s^*$ intersects $\hsO$.  This argument works for any open neighbourhood $\hsO$ of $\hs$; thus, every open neighbourhood of $\hs$ overlaps every element of $s^*$; hence  (\ref{gleason.covering.defn}) yields $p(s^*)=\hs$.
\ethmprf

\bthmprf[Proof of Theorem \ref{gleason.covering}.]
The fact that $\mu^*$ can be extended to a unique Borel probability  measure on $\sS^*$ follows immediately from
Lemma \ref{stonean.measure.extension}.  It remains to  prove equation (\ref{gleason.covering.e1}). 

\noindent For notational simplicity, we will invoke the canonical homeomorphism between $\sS^*$ and $\sig(\hgB)$  defined by formula
(\ref{gleason.covering.defn0}),  and identify each element (or subset) of $\sS^*$ with the corresponding
 element (or subset) of $\sig(\hgB)$; in effect we will assume  that $\sS^*=\sig(\hgB)$.  Thus,  each point in $\sS^*$ is
identified with an ultrafilter in the Boolean algebra $\hgB$, so
Lemma \ref{gleason.covering.map.lemma}  takes the following simpler form:

\claim{For any $s^*\in\sS^*$ and  $\hs\in\hsS$, we have $p(s^*)=\hs$ if and only if every  
$\hgB$-neighbourhood of $\hs$ is an element of $s^*$.\label{gleason.covering.C2}}

\noindent For any  $g\in\Cb(\sS,\Real)$, let  $\hg\in\sC(\hsS,\Real)$ be the (unique)
Stone-\v{C}ech extension of $g$.

\claim{If  $g\in\GB(\sS)$, then  $\hg\in\GhB(\hsS)$.\label{gleason.covering.C0}}
\bclaimprf
Every function in $\GB(\sS)$ is a  uniform limit of linear combinations  of functions in $\CB(\sS)$,
 and the transformation $\sC(\sS,\Real)\ni f\mapsto\hf\in \sC(\hsS,\Real)$ is  linear and continuous.    So
it suffices to show that $\hg\in\ChB(\hsS)$ whenever $g\in\CB(\sS)$.
So, suppose $g\in \CB(\sS)$. Let $r\in\Real$, and let $\hsB:=\Int\lb(\hg^{-1}(-\oo,r]\rb)$.
Then $\hsB$ is a regular open subset of $\hsS$ (because $\hg$ is continuous).  We must show that $\hsB\in\hgB$. So, 
let $\sB:=\sS\intsct\hsB$.   Then
\beqn
\label{gleason.covering.C0.e1}
\sB \quad=\quad \sS\intsct\Int\lb(\hg^{-1}(-\oo,r]\rb)
\quad\eeequals{(*)}\quad
\Int\lb(\sS\intsct\hg^{-1}(-\oo,r]\rb)
\quad\eeequals{(\dagger)}\quad
\Int\lb(g^{-1}(-\oo,r]\rb)
\ \in \ \gB.
\eeqn
 Here, $(*)$  is because $\sS$ is an open subset of $\hsS$ (because $\sS$ is locally compact),
 $(\dagger)$ is because $g:=\hg_{\restr\sS}$ by the definition of $\hg$, and the last
 step is because $g\in\CB(\sS)$.
 
 But $\hgB$ was defined using the isomorphism from
 Lemma  \ref{stone.cech.boolean.isomorphism} ; in other words,
 $\hgB=\{\hsR\in\gR(\hsS)$; \ $\sS\intsct\hsR\in \gB\}$.
 Thus, equation (\ref{gleason.covering.C0.e1}) establishes that $\hsB\in\hgB$.

 By a very similar argument, $\Int\lb(\hg^{-1}[r,\oo)\rb)\in\hgB$.
 This argument works for all $r\in\Real$;  thus, $\hg\in\ChB(\sS)$, as desired.
\eclaimprf

\noindent For any $\hsB\in\hgB$, let $\sB^*:=\{s^*\in\sS^*$; \ $\hsB\in s^*\}$  (a clopen subset of $\sS^*$).  The function $\hsB\mapsto\sB^*$ is a Boolean algebra isomorphism from $\hgB$ to
$\clop(\sS^*)$, by Stone's  Theorem.
 
\Claim{If $\hsB\in\hgB$ then $p^{-1}(\hsB)\subseteq\sB^*\subseteq p^{-1}[\Cl(\hsB)]$.\footnote{\label{gleason.covering.C3.footnote}It is tempting to think that $p^{-1}(\hsB)=\sB^*$ for all $\hsB\in\hgB$. But this cannot be true.
To see this,  note that $(\neg\hsB)^* = \compl{(\hsB^*)}$ (by the Stone Representation Theorem).
 Thus, if $p^{-1}(\hsB)=\sB^*$  and $p^{-1}(\neg\hsB)=(\neg\hsB)^* = \compl{(\sB^*)}$, then we would have $p^{-1}(\partial\hsB)=\emptyset$, contradicting the  surjectivity of $p$.}\label{gleason.covering.C3}}
\bclaimprf
Let $s^*\in\sS^*$, and let $\hs:=p(s^*)$. 
Then  $\lb(s^*\in p^{-1}(\hsB)\rb) \Leftrightarrow (\hs\in\hsB) \implies (\hsB\in s^*) \Leftrightarrow  (s^*\in \sB^*)$, where ``$\implies$'' is by Claim \ref{gleason.covering.C2}.  Thus, $p^{-1}(\hsB)\subseteq\sB^*$.

Now let $s^*\in\sB^*$;  then $\hsB\in s^*$.   Let  $\hs:=p(s^*)$ and let $\hsO\subseteq\hsS$ be an open neighbourhood of $\hs$.  Then 
$\hsO$ contains a $\hgB$-neighbourhood $\hsQ$ of $\hs$, because $\hgB$ is a base for the topology of $\hsS$.
Now, $\hsQ\in s^*$ by Claim \ref{gleason.covering.C2}, and thus, $\hsQ\intsct\hsB\neq\emptyset$, because $s^*$ is a filter.
Thus, $\hsO\intsct\hsB\neq\emptyset$, because  $\hsQ\subseteq\hsO$.
Thus, $\hsB$ overlaps every neighbourhood of $\hs$; thus, $\hs\in\Cl(\hsB)$.  
Thus, $s^*\in p^{-1}[\Cl(\hsB)]$.  This shows that $\sB^*\subseteq p^{-1}[\Cl(\hsB)]$.
\eclaimprf

\noindent  Formula (\ref{integral.defn}) says  that we can obtain $\integral_\sS g\dmu$ by approximating $g$ from
below by $\gB$-simple functions.  Thus we need to translate $\gB$-simple functions over $\sS$ into $\hgB$-simple functions on $\hsS$ and 
into $\clop(\sS^*)$-simple functions on $\sS^*$.   It will be convenient to work with a particular class of simple function.
A  function $\hf:\hsS\into\Real$ is a {\dfn $\hgB$-pyramidal function} 
if there exists a nested sequence of  subsets $\hsS=\hsB_0\supseteq\hsB_1\supseteq\cdots\supseteq\hsB_N$,
with $\hsB_1,\ldots,\hsB_N\in\hgB$,  and real numbers $r_0\in\Real$ and $r_1,\ldots,r_N\in\Real_+$ such that
\beqn
\label{pyramidal.function}
\hf\quad=\quad\sum_{n=0}^N r_n \,\bone_{\hsB_n}.
\eeqn
 For notational convenience, we define $\hsB_{N+1}:=\emptyset$. 
For all $n\in [0\ldots N]$, let $R_n:=r_0+\cdots+r_n$;  then $\hf(\hb)=R_n$ for all $\hb\in\hsB_n\intsct (\neg\hsB_{n+1})$.
Observe that $R_0<R_1<\cdots < R_N$ (hence the term ``pyramidal'').
Let $\hsF$ be the set of all $\hgB$-pyramidal  functions on $\hsS$.
Likewise, let $\sF'$ be the set of all $\gB$-pyramidal functions on $\sS$.
Finally, let $\sF^*$ be the set of all $\clop(\sS^*)$-pyramidal functions on $\sS^*$.

Let $\hf\in\hsF$ be as in formula (\ref{pyramidal.function}), where $\hsB_0=\hsS$ and
  $\hsB_1,\ldots,\hsB_N\in\hgB$.  For all $n\in[1\ldots N]$, let $\sB^*_n:=\{s^*\in\sS^*$; \ $\hsB_n\in s^*\}$.
  Define
\beqn
\label{pyramidal.function2}
f^*\quad:=\quad \sum_{n=0}^N r_n \,\bone_{\sB^*_n}
\qquad\And\qquad
\barf\quad:=\quad\sum_{n=0}^N r_n \,\bone_{\Cl(\hsB_n)}.
\eeqn
   where of course, $\sB^*_0=\sS^*$ and $\Cl(\hsB_0)=\hsS$.
    
\Claim{{\rm(a)} The transformation $\hf\mapsto \hf_{\restr\sS}$ is a bijection from $\hsF$ to $\sF'$.

{\rm(b)} The transformation $\hf\mapsto f^*$  defined by {\rm(\ref{pyramidal.function2})} is a bijection from
$\hsF$ to $\sF^*$. 

{\rm(c)}  For any $\hf\in\hsF$, we have $\hf\circ p \leq f^* \leq \barf\circ p$.
\label{gleason.covering.C5}
}
\bclaimprf
Part (a)  follows from the fact that the map $\hsB\mapsto\hsB\intsct\sS$  is a bijection from $\hgB$ to $\gB$, as shown by
Lemma  \ref{stone.cech.boolean.isomorphism} .
Part (b)  follows from the Stone Representation Theorem.  Part (c) follows from Claim \ref{gleason.covering.C3}.
\eclaimprf

\Claim{Let $\hf\in\hsF$ and let  $\hg\in\sC(\hsS,\Real)$.  If $\hf\leq \hg$, then $\barf\leq \hg$.\label{gleason.covering.C4}} 
  \bclaimprf
  Suppose $\hf\in\hsF$ is as in formula (\ref{pyramidal.function}), where
  $\hsB_0=\hsS$ and
  $\hsB_1,\ldots,\hsB_N\in\hgB$;  thus,
    $\barf$ is as in formula (\ref{pyramidal.function2}) (right).
Let $\hs\in\hsS$;  we must show that $\barf(\hs)\leq \hg(\hs)$.   Let 
$\hsB_{N+1}:=\emptyset$.

Observe that $\hsS=\Cl(\hsB_0)\supseteq\Cl(\hsB_1)\supseteq\cdots\supseteq\Cl(\hsB_N)$.
Thus, if $\hs\in\Cl(\hsB_m)$, then $\hs\in\Cl(\hsB_n)$ for all $n\leq m$.  Let $m:=\max\{n\in[0\ldots N]$; \ $\hs\in\Cl(\hsB_n)\}$, and let $R:=r_0+\cdots+r_m$.
 Then there is a net  $(\hb_j)_{j\in\sJ}$ in $\hsB_m$ converging to $\hs$ (where $\sJ$ is
some directed set).  Now, $\hs\not\in\Cl(\hsB_{m+1})$, so  $\neg\hsB_{m+1}$ is an open neighbourhood of $\hs$,
so there is some $j_0\in \sJ$ such that  for all $j>j_0$, we have $\hb_j\in\hsB_m\intsct(\neg\hsB_{m+1})$,
and thus, $\hf(\hb_j)=R$, which means  $R\leq \hg(\hb_j)$ (because $\hf\leq \hg$). 
  Thus, $R\leq \hg(\hs)$,  because $\hg$ is continuous and $(\hb_j)_{j\in\sJ}$ converges to $\hs$.  
But $\barf(\hs)=R$, because $\hs\in \Cl(\hsB_m)\setminus\Cl(\hsB_{m+1})$.  Thus, $\barf(\hs)\leq \hg(\hs)$.
This argument holds for all  $\hs\in\hsS$; thus, $\barf\leq \hg$, as claimed.
\eclaimprf

\noindent  Let $\sF$ be the set of {\em all} $\gB$-simple functions on 
$\sS$. For any $g\in\GB(\sS)$, define $\sF_g:=\{f\in\sF$; \ $f\leq g\}$, as in Section \ref{S:integration}.
Let $\sF'_g:=\{f\in\sF'$; \ $f\leq g\}$.  Then $\sF'_g\subseteq\sF_g$, because
$\sF'\subseteq\sF$. 

\Claim{For any $f\in\sF_g$, there is $f'\in\sF_g'$ with
$\D \Integral_\sD f'\dmu = \Integral_\sD f\dmu$ for all $\sD\in\gB$.\label{gleason.covering.C5a}}
\bclaimprf
By definition, there is a $\gB$-partition $\{\sP_1,\ldots,\sP_L\}$ of $\sS$
such that $f$ is constant on each cell of the partition.   Let $R_0<R_1<\cdots<R_N$ be the (finite) set of values which
$f$ takes on these cells.    For all $n\in[0\ldots N]$, let $\sL_n:=\{\ell\in[1\ldots L]$; \ $f$ takes the value $R_n$ on $\sP_\ell\}$, and then define $\sP'_n:=\bigvee_{\ell\in\sL_n} \sP_\ell$;
thus, $\sP'_1,\ldots,\sP'_N\in\gB$ are disjoint.
  Finally, for all $n\in[0\ldots N]$, define $\sB_n:=\sP'_n \vee\sP'_{n+1}\vee\cdots\vee\sP'_N\in\gB$.   Then $\sS=\sB_0\supseteq\sB_1\supseteq\cdots\supseteq\sB_N$. 
Meanwhile, define $r_0:=R_0$,
and for all $n\in[1\ldots N]$, let $r_n:=R_n-R_{n-1}$.  Thus, for any $n\in[1\ldots N]$, $r_n>0$, and 
$R_n=r_0+\cdots+r_n$.   Define 
\beqn
 \label{gleason.covering.C5a.e1}
f'\quad:=\quad \sum_{n=0}^N r_n \,\bone_{\sB_n}.
\eeqn
Then $f'\in\sF'$, and for any $\sD\in\gB$, we have
\beq \Integral_\sD f'\dmu &\eeequals{(*)}&
\sum_{n=0}^N \Integral_\sD r_n \,\bone_{\sB_n} \dmu 
\quad \eeequals{(\sharp)}\quad 
\sum_{n=0}^N r_n\,\mu[\sD\intsct \sB_n]
\quad \eeequals{(\dagger)} \quad \sum_{n=0}^N r_n\,\mu\lb[\sD\intsct  \bigvee_{m=n}^N \sP'_m\rb]
\\ &=&  \sum_{n=0}^N r_n \, \lb(\sum_{m=n}^N \mu[\sD\intsct \sP'_m]\rb)
\quad = \quad \sum_{n=0}^N \sum_{m\geq n} r_n \, \mu[\sD\intsct \sP'_m]
\\& = & \sum_{m=0}^N \sum_{n\leq m}  r_n  \, \mu[\sD\intsct \sP'_m]
\quad\eeequals{(\diamond)}\quad  \sum_{m=0}^N R_m\, \mu[\sD\intsct \sP'_m]
\\& \eeequals{(\ddagger)} & \sum_{m=0}^N R_m\, \mu\lb[\sD\intsct \bigvee_{\ell\in\sL_m} \sP_\ell\rb]
\quad =\quad \sum_{m=0}^N \sum_{\ell\in\sL_m} R_m \,\mu[\sD\intsct \sP_\ell]
\quad \eeequals{(\sharp)}\quad \Integral_\sD f\dmu,
\eeq
as claimed.  Here,  $(*)$ is by equation (\ref{gleason.covering.C5a.e1}), while both $(\sharp)$ equalities are by   equation (\ref{integral.of.simple.function})  from Theorem \ref{from.probability.to.expectation0} .  Next,
$(\dagger)$ is by the definition of $\sB_n$ and $(\diamond)$ is because $r_0+\cdots+r_m=R_m$,
by definition of $r_0,\ldots,r_N$.  Finally, $(\ddagger)$ is by the definition of $\sP'_m$.

Finally, let $s\in\sS$; we must show that $f'(s)\leq g(s)$. 
The set $\sQ:=\sP_1\union\ldots\union\sP_L$ is dense in $\sS$, so there is a net $(q_j)_{j\in\sJ}$ in $\sQ$
that converges to $s$ (for some directed set $\sJ$).  By dropping to a subnet if necessary, we can find
some $\ell\in[1\ldots L]$ such that $q_j\in\sP_\ell$ for all $j\in\sJ$.  Find the unique $n\in[1\ldots N]$ such that
$\ell\in\sL_n$. Then $\sP_\ell\subseteq\sP'_n$, so for all $j\in\sJ$, we have
 $q_j\in\sP'_n$ and  $f(q_j)=R_n$.  Thus, $R_n\leq g(q_j)$ for all $j\in\sJ$ (because $f\leq g$), and thus
 $R_n \leq g(s)$, because $g$ is continuous, and $(q_j)_{j\in\sJ}$ converges to $s$.

Now, the net $(q_j)_{j\in\sJ}$ is a subset of $\sP'_n$, which is disjoint from the open set $\sB_{n+1}$.
Thus, $s\not\in\sB_{n+1}$.  Thus, $f'(s)\leq r_0+\cdots+r_n = R_n$, by
 defining formula (\ref{gleason.covering.C5a.e1}).  But we have already established that $R_n \leq g(s)$;
 thus, $ f'(s)\leq g(s)$, as desired. 
 
  This argument holds for all $s\in\sS$; thus, $ f'\leq g$, and thus, $f'\in \sF'_g$.
\eclaimprf

\noindent For any $\hg\in\sC(\hsS,\Real)$, let $\hsF_{\hg}:=\{\hf\in\hsF$; \ $\hf\leq \hg\}$.
Likewise, for any $g\in\sC(\sS^*,\Real)$, let $\sF^*_g:=\{f\in\sF^*$; \ $f\leq g\}$.

\Claim{Let $g\in\GB(\sS)$,  let $\hg$ be the unique extension of $g$ to $\GhB(\hsS)$
(as in Claim \ref{gleason.covering.C0}) 
and let  $g^*:=\hg \circ p$.

{\rm(a)} The  transformation  $\hf\mapsto \hf_{\restr\sS}$ is a bijection from $\hsF_{\hg}$ to  $\sF'_g$.

{\rm(b)} The  transformation  $\hf\mapsto f^*$  defined by {\rm(\ref{pyramidal.function2})} is a bijection from
$\hsF_{\hg}$ to $\sF^*_{g^*}$. 
\label{gleason.covering.C6}
}
\bclaimprf (a) If $\hf\in\hsF_{\hg}$, then $\hf\leq \hg$, and
hence $\hf_{\restr\sS} \leq \hg_{\restr\sS}=g$.  Thus, $\hf_{\restr\sS}\in \sF'_g$.
The transformation $\hf\mapsto \hf_{\restr\sS}$ is injective by Claim \ref{gleason.covering.C5}(a).
To see that it is surjective, let $f\in\sF'_g$.  The surjectivity in Claim \ref{gleason.covering.C5}(a) yields
 a unique $\hf\in\hsF$ such that $(\hf)_{\restr\sS}=f$.
We must show that $\hf\in\hsF_{\hg}$.  So, let $\hs\in\hsS$;
we must show that $\hf(\hs)\leq \hg(\hs)$.

 Let   $\hf$ be as in formula (\ref{pyramidal.function}), where
  $\hsB_0=\hsS$ and  $\hsB_1,\ldots,\hsB_N\in\hgB$.    Let 
$\hsB_{N+1}:=\emptyset$.  Let $m:=\max\{n\in[0\ldots N]$; \ 
  $\hs\in\Cl(\hsB_n)\}$, and let $R:=r_0+\cdots+r_m$.
Let $\sB_m:=\sS\intsct\hsB_m$.  Then $\sB_m$ is  dense in $\hsB_m$, because $\sS$ is dense in $\hsS$,
and $\hsB_m$ is an open subset of $\hsS$.  Thus, $\Cl(\sB_m)=\Cl(\hsB_m)$.  But $\hs\in\Cl(\hsB_m)$,
so there is a net $(b_j)_{j\in\sJ}$ in $\sB_m$ converging to $\hs$ ($\sJ$ 
some directed set).    We have $f(b_j)\leq g(b_j)$ for all $j\in\sJ$, because $f\leq g$ (because $f\in\sF'_g$).
Now,  $\hs\not\in\Cl(\hsB_{m+1})$, so 
 $\neg\hsB_{m+1}$ is an open neighbourhood of $\hs$.    Thus, there exists $j_0\in\sJ$ such that for all $j>j_0$,
 we have $b_j\in\sB_m\intsct(\neg\hsB_{m+1})$, and thus, $f(b_j)=R$, which means $R\leq g(b_j)$ (because $f\leq g$), and hence  $R\leq \hg(b_j)$ (because $\hg_{\restr\sS}=g$).
 Thus, $R\leq \hg(\hs)$, because $\hg$ is continuous and  $(b_j)_{j\in\sJ}$ converges to $\hs$.
 But $\hf(\hs)\leq R$ (because $\hs\in\neg\hsB_{m+1}$), so this means
that $\hf(\hs)\leq \hg(\hs)$, as desired.  

This argument works for all $\hs\in\hsS$; thus
 $\hf\leq \hg$, and thus, $\hf\in\hsF_{\hg}$, as desired.

\item (b) If $\hf\in\hsF_{\hg}$, then $\hf\leq \hg$.  Thus, $\barf\leq \hg$, by Claim \ref{gleason.covering.C4}.
Thus, $\barf\circ p\leq \hg\circ p $.  But $f^*\leq \barf\circ p$ by 
Claim \ref{gleason.covering.C5}(c), while $\hg\circ p=g^*$ by definition. Thus, $f^*\leq g^*$.
Thus, $f^*\in\sF^*_{g^*}$.  The function $\hf\mapsto f^*$ is injective by Claim \ref{gleason.covering.C5}(b).
It remains to show that it is surjective.

Let $f_0\in\sF^*_{g^*}$.   Surjectivity from Claim \ref{gleason.covering.C5}(b) yields  $\hf\in\hsF$ with 
 $f_0=f^*$.  We have
\beqn
\label{gleason.covering.C6.e1}
\hf\circ p \quad\leeeq{(*)}\quad f^* \quad=\quad f_0\quad \leeeq{(\dagger)}\quad g^* \quad=\quad \hg\circ p,
\eeqn
where $(*)$ is by   Claim \ref{gleason.covering.C5}(c), $(\dagger)$ is because $f_0\in\sF^*_{g^*}$, and
the equalities are true by the definitions of the functions in question.
  But  $p$ is surjective, so inequality (\ref{gleason.covering.C6.e1}) implies that
$ \hf\leq \hg$.  Thus, $\hf\in\hsF_{\hg}$, as desired.
\eclaimprf

\noindent 
For any $f\in\sF^*$, let  $\int_{\sS^*} f\dmu^*$ be the standard Lebesgue integral with respect
to the Borel measure $\mu^*$. 

\Claim{ For any $\hf\in\hsF$ and $\sD\in\gB$, we have $\D \Integral_\sD \hf_{\restr\sS} \dmu \ = \ \int_{\sD^*} f^*\dmu^*$.
\label{gleason.covering.C7}
}
\bclaimprf
Suppose  $\hf\in\hsF$ is as in formula (\ref{pyramidal.function}), where
  $\hsB_0=\hsS$ and  $\hsB_1,\ldots,\hsB_N\in\hgB$.  For all $n\in[1\ldots N]$, let 
$\sB_n:=\hsB_n\intsct\sS$.   Then $\sB_n\in\gB$, and for any $\sD\in\gB$, we have
\beq
\Integral_\sD \hf_{\restr\sS} \dmu  &\eeequals{(*)}&
\sum_{n=0}^N r_n \,\mu[\sD\intsct\sB_n]
\quad\eeequals{(\dagger)}\quad 
\sum_{n=0}^N r_n \,\mu^*[(\sD\intsct\sB_n)^*]
\\& \eeequals{(\ddagger)}& 
\sum_{n=0}^N r_n \,\mu^*[\sD^*\intsct\sB^*_n]
\quad\eeequals{(\diamond)}\quad 
\int_{\sD^*} f^* \dmu^*,
\eeq
as claimed.
Here,  $(*)$ is by the linearity of $\integral$, because $\hf_{\restr\sS}=  \sum_{n=0}^N r_n \,\bone_{\sB_n}$.  Meanwhile, $(\dagger)$ is by defining formula (\ref{mu.star.defn}), 
and $(\ddagger)$ is because the transformation $\sB\mapsto\sB^*$ is a Boolean algebra homomorphism
from $\gB$ to $\clop(\sS^*)$.
Finally $(\diamond)$ is by  equation (\ref{pyramidal.function2}) (left)  and the linearity of the Lebesgue integral.
\eclaimprf

\noindent
Now let $g\in\GB(\sS)$ and let $\sD\in\gB$.  Then
\beq
\Integral_\sD g\dmu &\eeequals{(*)} &
\sup_{f\in\sF_g}  \Integral_\sD f \dmu 
\quad\eeequals{(\sharp)}\quad
\sup_{h\in\sF'_g}  \Integral_\sD h \dmu 
\quad\eeequals{(\dagger)}\quad
\sup_{\hf\in\hsF_{\hg}}  \Integral_\sD \hf_{\restr\sS} \dmu 
\\ &\eeequals{(\diamond)}&
\sup_{\hf\in\hsF_{\hg}}  \int_{\sD^*} f^* \dmu^* 
\quad\eeequals{(\ddagger)}\quad
\sup_{h\in\sF^*_{g^*}}  \int_{\sD^*} h \dmu^* 
\quad\eeequals{(@)}\quad
\int_{\sD^*} g^*\dmu^*,
\eeq
which proves equation (\ref{gleason.covering.e1}).
Here, $(*)$ is by  equation (\ref{integral.defn})  from Theorem \ref{from.probability.to.expectation0} ,
and $(\sharp)$ is by Claim \ref{gleason.covering.C5a}.
Next, $(\dagger)$ is by Claim \ref{gleason.covering.C6}(a), and
$(\diamond)$ is by Claim \ref{gleason.covering.C7}.
Next $(\ddagger)$ is by Claim \ref{gleason.covering.C6}(b).
Finally, to see (@), recall that we showed in the proof of Lemma \ref{stonean.measure.extension}
that   $\sF^*$ is uniformly dense in $\sC(\sS^*,\Real)$, by invoking the Stone-Weierstrass Theorem.
Thus, we can uniformly approximate $g^*$ from below with elements
of $\sF^*_{g^*}$.  Thus,  (@) follows from the fact
that  $\mu^*$-integration is continuous with respect to the uniform norm on $\sC(\sS^*,\Real)$.
\ethmprf

\section{A categorical perspective\label{S:cat}}
\setcounter{equation}{0}

A {\dfn content space} is an ordered triple $(\sS,\gB,\mu)$, where $\sS$ is a locally compact Hausdorff space,
$\gB$ is a generative Boolean subalgebra of $\gR(\sS)$,
and $\mu$ is a content on $\gB$ with $\mu[\sS]=1$.
If $\bsC_1=(\sS_1,\gB_1,\mu_1)$ and $\bsC_2=(\sS_2,\gB_2,\mu_2)$ are two content spaces, then
a {\dfn morphism} from $\bsC_1$ to $\bsC_2$ is a continuous function $\phi:\sS_1\into\sS_2$
such that $\phi^{-1}:\gB_2\into\gB_1$ is a Boolean algebra homomorphism and 
$\phi(\mu_1)=\mu_2$.\footnote{Lemma \ref{measurability.lemma}(b) and
Remark \ref{clopen.algebra.remark} provide sufficient conditions for
 $\phi^{-1}$ to be a Boolean algebra homomorphism from $\gB_2$ to $\gB_1$.}
  It is easily verified that the composition of two
such morphisms is also a content space morphism.   Let $\Cred$ be the category of content spaces and their morphisms.

Let $\CmpCred$ be the full subcategory of $\Cred$ consisting of all {\em compact} content spaces and their morphisms.
Via  Stone-\v{C}ech compactification, we can define a functor from $\Cred$ into $\CmpCred$ as follows.
Let $\bsC=(\sS,\gB,\mu)$ be a content space.     Lemma  \ref{stone.cech.boolean.isomorphism}  
yields a Boolean algebra isomorphism $h:\gR(\hsS)\, \widetilde{\longrightarrow}\, \gR(\sS)$ given by $h(\hsR)=\hsR\intsct\sS$.
Let $\hgB:=h^{-1}(\gB)$; then $\hgB$ is a Boolean subalgebra of $\gR(\hsS)$, and is isomorphic to $\gB$ via $h$.
For any $\hsB\in\hgB$, define $\hmu[\hsB]:=\mu[\hsB\intsct\sS]$;  then $\hmu$ is a content on $\hgB$.
The triple $\hbsC:=(\hsS,\hgB,\hmu)$ is a compact content space,  which we will call the {\dfn Stone-\v{C}ech compactification} of $\bsC$.  

Now let $\bsC_1=(\sS_1,\gB_1,\mu_1)$ and $\bsC_2=(\sS_2,\gB_2,\mu_2)$ be two content spaces, and
let $\phi:\sS_1\into\sS_2$ be a content space morphism.  Since $\sS_2\subseteq\hsS_2$, we can
regard $\phi$ as a continuous function from $\sS_1$ into the compact space $\hsS_2$.  The Stone-\v{C}ech
Extension Theorem yields a unique extension to a continuous function $\hphi:\hsS_1\into\hsS_2$.

\Proposition{\label{stone.cech.compactification.functor}}
{
$\hphi$ is a content space morphism from $\hbsC_1$ to $\hbsC_2$.  The transformation $\bet$ defined by
$\bsC\mapsto \hbsC$ and $\phi\mapsto\hphi$  is a faithful functor from $\Cred$ to $\CmpCred$.
}
\bthmprf  We already know that $\hphi$ is continuous.  We must show that  $\hphi^{-1}:\hgB_2\into\hgB_1$ is a Boolean algebra homomorphism and that $\hphi[\hmu_1]=\hmu_2$.

\item {\em Homomorphism.}
Recall that the Boolean algebra isomorphisms $h_1:\gR(\hsS_1)\, \widetilde{\longrightarrow}\, \gR(\sS_1)$ 
and $h_2:\gR(\hsS_1)\, \widetilde{\longrightarrow}\, \gR(\sS_2)$ from  Lemma  \ref{stone.cech.boolean.isomorphism}  are defined by $h_1(\hsR_1):=\hsR_1\intsct\sS_1$ and $h_2(\sR_2):=\hsR_2\intsct\sS_2$
for all $\hsR_1\in\gR(\hsS_1)$ and $\hsR_2\in\gR(\hsS_2)$.   Furthermore,
$\hgB_1:=h_1^{-1}(\gB_1)$ and $\hgB_2:=h_2^{-1}(\gB_2)$, so by restricting $h_1$ and $h_2$, we get
isomorphisms  $h_1:\hgB_1\, \widetilde{\longrightarrow}\, \gB_1$
and $h_2:\hgB_2\, \widetilde{\longrightarrow}\, \gB_2$. 
Let   $\eta_1:\gB_1\, \widetilde{\longrightarrow}\, \hgB_1$ be the inverse of $h_1$ ---this is also an isomorphism.
Finally, $\phi^{-1}:\gB_2\into\gB_1$ is a Boolean algebra homomorphism because $\phi$ is a content space morphism.
Thus, the fact that $\hphi^{-1}:\hgB_2\into\hgB_1$ is a Boolean algebra homomorphism is
an immediate consequence of the following commuting diagram
\[
\begin{tikzcd}
\hgB_2 \arrow{r}{h_2}\arrow[swap]{d}{\hphi^{-1}} & \gB_2 \arrow{d}{\phi^{-1}}\\
\hgB_1 \arrow[swap,shift right=0.2em]{r}{h_1} & \arrow[swap,shift right=0.2em]{l}{\eta_1} \gB_1 
\end{tikzcd}
\]
\item {\em Measure-preserving.} \
Let $\hsB_2\in\hgB_2$.  If $\sB_2:=h_2(\hsB_2)$, then
 $\hmu_2[\hsB_2]=\mu_2[\sB_2]$.
Likewise, if $\hsB_1:=\hphi^{-1}(\hsB_2)$ and $\sB_1:=h_1(\hsB_1)$, then
 $\hmu_1[\hsB_1]=\mu_1[\sB_1]$.  But from the above commuting diagram, we see that
 $\phi^{-1}(\sB_2)=\sB_1$.  Thus, $\mu_1[\sB_1]=\mu_2[\sB_2]$, because $\phi$
 is a content space morphism.  It follows that $\hmu_1[\hsB_1]=\hmu_2[\hsB_2]$,
 as desired.

\item {\em  $\bet$ is a faithful functor.} \ The proof is identical to the proof of the corresponding statement
for the Stone-\v{C}ech compactification functor in the category of topological spaces.
\ethmprf

\noindent
A {\dfn Stonean probability space} is an ordered pair $(\sS,\mu)$, where $\sS$ is a Stonean topological space,
and $\mu$ is a Borel measure on $\sS$.  Let $\StPr$ be the category of Stonean probability spaces and continuous, measure-preserving functions.
Via Theorem \ref{gleason.covering}, we can define a functor from $\Cred$ to $\StPr$ as follows.
Let $\bsC=(\sS,\gB,\mu)$ be a content space.  Let $\sS^*:=\sigma(\gB)$, and
let $\gB^*:=\clop(\sS^*)$;  thus, $\gB^*$ is isomorphic to $\gB$ by the Stone Representation Theorem.
Let $\mu^*$ be the Borel probability measure defined on $\sS^*$ by formula (\ref{mu.star.defn}) and Lemma \ref{stonean.measure.extension}.
Thus, $\bsC^*:=(\sS^*,\mu^*)$ is a Stonean probability space, which will call the {\dfn Stone representation}
of $\bsC$.

Now let $\bsC_1=(\sS_1,\gB_1,\mu_1)$ and $\bsC_2=(\sS_2,\gB_2,\mu_2)$ be two content spaces, and
let $\phi:\sS_1\into\sS_2$ be a content space morphism.  Since $\phi^{-1}:\gB_2\into\gB_1$
is Boolean algebra homomorphism, the Stone Duality Theorem yields a continuous function $\phi^*:\sS_1^*\into\sS_2^*$.

\Proposition{\label{stone.space.functor}}
{
$\phi^*$ is continuous, measure-preserving function from $\bsC^*_1$ to $\bsC^*_2$.  The transformation
$\sigma$ defined by $\bsC\mapsto\bsC^*$ and $\phi\mapsto\phi^*$ is a faithful functor from $\Cred$ to $\StPr$.
}
\bthmprf  
Suppose  $\bsC^*_1:=(\sS_1^*,\mu_1^*)$ and  $\bsC^*_2:=(\sS_2^*,\mu_2^*)$, where
$\sS_1=\sigma(\gB_1)$ and $\sS_2=\sigma(\gB_2)$.
The Stone Duality Theorem says
$\phi^*:\sS^*_1\into\sS^*_2$ is continuous.  To show that $\phi^*$ is measure-preserving, it suffices to show that it preserves the measures of clopen sets, because  Lemma \ref{stonean.measure.extension} says that a Borel measure on a Stonean space is entirely
determined by its values on clopen sets.  So, let $\sQ_2\in\clop[\sS^*_2]$ and let $\sQ_1:=(\phi^*)^{-1}(\sQ_2)$;  we must show that $\mu^*_1[\sQ_1]=\mu^*_1[\sQ_2]$.

Recall that $\sS_2^*$ is the set of all ultrafilters in $\gB_2$, and 
the Stone Representation Theorem says that there exists some $\sB_2\in\gB_2$ such that
$\sQ_2=\sB^*_2$, where $\sB^*_2:=\{s^*_2\in\sS^*_2$; \ $\sB_2\in s^*_2\}$.  Thus,  defining formula (\ref{mu.star.defn})
says that $\mu^*_2[\sQ_2]=\mu_2[\sB_2]$.
Let $\sB_1:=\phi^{-1}(\sB_2)$;  then $\sB_1\in\gB_1$ (because
$\phi$ is a content space morphism), and for all $s^*_1\in\sS^*_2$, we have
\[
\statement{$s^*_1\in\sQ_1$}
\iff
\statement{$\phi^*(s^*_1)\in\sQ_2$}
\iff 
\statement{$\sB_2\in\phi^*(s^*_1)$}
\iiiff{(*)}
\statement{$\phi^{-1}(\sB_2)\in s^*_1$}
\iff
\statement{$\sB_1\in s^*_1$},
\]
where $(*)$ is because $\phi^*(s^*_1)=\{\sB\in\gB_2$; \ $\phi^{-1}(\sB)\in s^*_1\}$.   Thus, we see that
$\sQ_1=\sB^*_1$, where $\sB^*_1:=\{s^*_1\in\sS^*_1$; \ $\sB_1\in s^*_1\}$.
Thus, defining formula (\ref{mu.star.defn})
says that $\mu^*_1[\sQ_1]=\mu_1[\sB_1]$.   But $\mu_1[\sB_1]=\mu_2[\sB_2]$ because
$\phi(\mu_1)=\mu_2$ and $\phi^{-1}(\sB_2)=\sB_1$.  Thus, we conclude that
$\mu^*_1[\sQ_1]=\mu^*_2[\sQ_2]$, as desired.

\item {\em Functor.}  \ Let $\bsC_3=(\sS_3,\gB_3,\mu_3)$,
 let $\sigma(\bsC_3):=(\sS_3^*,\mu_3^*)$, and
let $\psi:\sS_2\into\sS_3$ be another content space homomorphism.  
We must show that $(\psi\circ\phi)^* = \psi^*\circ\phi^*$.
To see this, recall that $(\psi\circ\phi)^*$,  $\psi^*$, and $\phi^*$ are the results of applying
the Stone Duality functor to the Boolean algebra homomorphisms
$(\psi\circ\phi)^{-1}:\gB_3\into\gB_1$,  $\psi^{-1}:\gB_3\into\gB_2$, and $\phi^{-1}:\gB_2\into\gB_1$,
respectively.   Furthermore, $(\psi\circ\phi)^{-1} = \phi^{-1}\circ\psi^{-1}$.
Thus, $(\psi\circ\phi)^* = \psi^*\circ\phi^*$.

\item {\em Faithful.} \ Let $\xi:\sS_1\into\sS_2$ be another content space homomorphism, and
suppose $\xi^*=\phi^*$.  By Stone Duality, this means that
$\xi^{-1}:\gB_2\into\gB_1$ and $\phi^{-1}:\gB_2\into\gB_1$ are the same Boolean algebra homomorphism;
in other words, $\xi^{-1}(\sB_2)=\phi^{-1}(\sB_2)$ for all $\sB_2\in\gB_2$.
Let $s\in\sS_1$ and suppose $\phi(s)\neq\xi(s)$.  Then there exists $\sB_2\in\gB_2$  such that $\phi(s)\in\sB_2$
but $\xi(s)\not\in\gB_2$ (because  by Lemma \ref{generative.lemma}, $\gB_2$ is a base for the topology of $\sS_2$, which is Hausdorff).
Then $s\in\phi^{-1}(\sB)$ but $s\not\in\xi^{-1}(\sB)$, which contradicts  the fact that $\xi^{-1}(\sB_2)=\phi^{-1}(\sB_2)$.  By contradiction, we must have $\phi(s)=\xi(s)$ for all $s\in\sS$; thus, $\phi=\xi$.
\ethmprf

\noindent
An {\dfn  integration space} is a triple $\bsI=(\sS,\sG,\dI)$, where $\sS$ is a locally compact Hausdorff space,
$\sG$ is a linear subspace of $\sC(\sS,\Real)$, and $\dI:\sG\into\Real$ is a nondecreasing, bounded linear functional
with norm $1$.  
Heuristically, for any $g\in\sG$, we can think $\dI[g]$ as the ``integral'' of $g$ with respect to some hypothetical probability measure.     For  example, if $\bsC=(\sS,\gB,\mu)$ is a  content space, then we obtain an  integration space
$\Upsilon(\bsC):=(\sS,\sG, \integral_\sS \dmu )$, where $\sG:=\GB(\sS)$ and   $\integral_\sS\dmu$ is  from Theorem \ref{from.probability.to.expectation0} .

If  $\bsI_1=(\sS_1,\sG_1,\dI_1)$ and  $\bsI_2=(\sS_2,\sG_2,\dI_2)$ are two  integration spaces, then a {\dfn morphism} from $\bsI_1$ to $\bsI_2$ is a continuous function $\phi:\sS_1\into\sS_2$ such that,
for all $g\in\sG_2$, we have $g\circ\phi\in\sG_1$ and $\dI_1[g\circ\phi]=\dI_2[g]$.
For example, if $\bsC_1$ and $\bsC_2$ are two  content spaces, and
 $\phi$ is a content space morphism from $\bsC_1$ to
$\bsC_2$, then Proposition \ref{change.of.variables} says that $\phi$ is also an integration
space morphism from $\Upsilon(\bsC_1)$ to $\Upsilon(\bsC_2)$.  
Thus, if $\INT$ is the category of  integration spaces and their morphisms,
then the transformation $\Ups$ defined by $\bsC\mapsto \Upsilon(\bsC)$ and $\phi\mapsto\phi$ is 
 is a functor  from $\Cred$ to $\INT$.  If $\CpInt$ is the subcategory of {\em compact} integration
 spaces, then $\Upsilon$ restricts to a functor from $\CmpCred$ into $\CpInt$.
 Let $\Upsilon\circ\bet:\Cred\into\CpInt$ be the functor obtained by composing $\Upsilon$ with the
functor $\bet$ from Proposition \ref{stone.cech.compactification.functor}.

\Proposition{\label{stone.cech.natural.transform}}
{
\bthmlist
\item For any content space $\bsC=(\sS,\gB,\mu)$, the inclusion map $\iota_\bsC:\sS\hookrightarrow\hsS$
is an integration space morphism from $\Upsilon(\bsC)$ to $\Upsilon(\hbsC)$.

\item The collection $\{\iota_\bsC$; \ $\bsC\in\Cred\}$ is a natural transformation from
 $\Upsilon$ to  $\Upsilon\circ\bet$. 
\ethmlist
}
\bthmprf (a) Let $\hbsC=(\hsS,\hgB,\hmu)$.  Then  $\Upsilon(\hbsC)=(\hsS,\hsG, \integral_{\hsS} \dhmu )$, where
$\hsG=\GhB(\hsS)$, and  $\integral_{\hsS} \dhmu$ is from Theorem \ref{from.probability.to.expectation0} .
The inclusion map $\iota_\bsC:\sS\hookrightarrow\hsS$ is continuous because $\hsS$ is a compactification of $\sS$.
Let $\hg\in \GhB(\hsS)$, and let $g:=\hg\circ\iota_\bsC$; then $g=\hg_{\restr\sS}$.

\claim{$g\in\GB(\sS)$.\label{stone.cech.natural.transform.C1}}
\bclaimprf
Every function in $\GhB(\hsS)$ is a  uniform limit of linear combinations   of functions in $\ChB(\hsS)$,
 and the transformation $\sC(\hsS,\Real)\ni \hf\mapsto\hf_{\restr\sS}\in \sC(\sS,\Real)$ is  linear and continuous.  So
it suffices to show that $g\in\CB(\sS)$ whenever $\hg\in\ChB(\hsS)$.
So, suppose $\hg\in \ChB(\hsS)$. Let $r\in\Real$, and let $\sB:=\Int\lb(g^{-1}(-\oo,r]\rb)$;
we must show that $\sB\in\gB$.
To see this, let $\hsB:=\Int\lb(\hg^{-1}(-\oo,r]\rb)$; then we know that
$\hsB\in\hgB$ because $\hg\in\ChB(\hsS)$.
Now, as shown in equation (\ref{gleason.covering.C0.e1}) 
in the proof of Theorem \ref{gleason.covering}, $\sB=\sS\intsct\hsB$.  
 But 
 Lemma  \ref{stone.cech.boolean.isomorphism}  says that
 $\hgB=\{\hsR\in\gR(\hsS)$; \ $\sS\intsct\hsR\in \gB\}$.
 Thus $\sB\in\gB$.
 By a very similar argument, $\Int\lb(g^{-1}[r,\oo)\rb)\in\gB$.
 This argument works for all $r\in\Real$;  thus, $g\in\CB(\sS)$, as desired.
\eclaimprf

\noindent It remains to show that $\integral_\sS g \dmu =\integral_{\hsS} \hg \dhmu$.
Let $\sF$ be the set of $\gB$-simple functions of $\sS$, and let $\hsF$ be the set of
$\hgB$-simple functions on $\hsS$.

\Claim{For any $\hf\in\hsF$, if $f=\hf_{\restr\sS}$, then $f\in\sF$, and 
$\D\Integral_\sS f\dmu = \Integral_{\hsS}\hf \dhmu$.\label{stone.cech.natural.transform.C2}}
\bclaimprf
 Let $\hf\in\hsF$ be subordinate to the $\hgB$-partition $\{\hsB_1,\ldots,\hsB_N\}$.
Then $\hf_{\restr\sS}$ is a simple function on $\sS$, subordinate to the
 $\gB$-partition $\{\sB_1,\ldots,\sB_N\}$, where for all $n\in[1\ldots N]$, $\sB_n:=\sS\intsct\hsB_n$, 
 and the value that $\hf_{\restr\sS}$ takes on $\sB_n$ is the same as the value that $\hf$ takes on $\hsB_n$
 ---call this value $r_n$.  Thus,
 \[
 \Integral_\sS f \dmu \quad\eeequals{(*)}\quad \sum_{n=1}^N r_n\, \mu[\sB_n]
 \quad\eeequals{(\dagger)}\quad \sum_{n=1}^N r_n\, \hmu[\hsB_n]
 \quad\eeequals{(*)}\quad
  \Integral_{\hsS}\hf \dhmu.
\]
Here, both $(*)$ are by   equation (\ref{integral.of.simple.function})  from Theorem \ref{from.probability.to.expectation0} .
$(\dagger)$ is by the definition of $\hmu$.
\eclaimprf

\noindent Let $\sF_g:=\{f\in\sF$; \ $f\leq g\}$, and let $\hsF_\hg:=\{\hf\in\hsF$; \ $\hf\leq \hg\}$.

\claim{$\sF_g=\{\hf_{\restr\sS}$; \ $\hf\in\hsF_\hg\}$.\label{stone.cech.natural.transform.C3}}
\bclaimprf
``$\supseteq$'' \ If $\hf\in\hsF$, then $\hf_{\restr\sS}\in\sF$ by  Claim \ref{stone.cech.natural.transform.C2}.
 Furthermore, if $\hf\in\hsF_\hg$, then $\hf(\hs)\leq \hg(\hs)$ for all $\hs\in\hsS$, which means
 that  $\hf_{\restr\sS}(s)\leq g(s)$ for all $s\in\sS$, and hence, $\hf_{\restr\sS}\in\sF_g$.

``$\subseteq$'' \ Let $f\in \sF_g$ be subordinate to the $\gB$-partition $\{\sB_1,\ldots,\sB_N\}$.
For all $n\in[1\ldots N]$, let $\hsB_n\in\hgB$ be the unique element such that $\sB_n=\sS\intsct\hsB_n$.
Then $\{\hsB_1,\ldots,\hsB_N\}$ is a $\hgB$-partition of $\hsS$ (by Lemma  \ref{stone.cech.boolean.isomorphism} ).
Let $\hf\in\hsF$ be the unique simple function on $\hsS$ subordinate to this partition,
such that for all $n\in[1\ldots N]$,
 the value that $\hf$ takes on $\hsB_n$ is the same as the value that $f$ takes on $\sB_n$,
 and $\hf_{\restr\partial\sB_n}=f_{\restr\partial\sB_n}$.\footnote{Recall that the value
 of a simple function on the boundaries of its subordinating partition is arbitrary, and has no effect on the
 integral;  the only important thing is that $\hf$ is dominated by $\hg$.} 
Meanwhile, define $\hf(\hs)=\hg(\hs)$ for all $\hs\in (\partial\hsB_n)\setminus\sS$, for all 
$n\in[1\ldots N]$.  Then clearly $\hf_{\restr\sS}=f$.  It remains to show that $\hf\in\hsF_\hg$.

For any $\hs\in\hsS$, we must show that $\hf(\hs)\leq\hg(\hs)$.
If $\hs\in\sS$, then $\hf(\hs)=f(\hs)\leq g(\hs)=\hg(\hs)$, so we're done.
So suppose $\hs\in\hsS\setminus\sS$.  If $\hs\in\partial\hsB_n$ for some $n\in[1\ldots N]$,
then  $\hf(\hs)=\hg(\hs)$ by definition, so we're done.  So suppose that $\hs\in\hsB_n$ for
some $n\in[1\ldots N]$.   Let $r_n$ be the (constant) value of $\hf$ on $\hsB_n$.

 Since $\sB_n$ is dense in $\hsB_n$, there is
a net $(b_j)_{j\in\sJ}$ in $\sB_n$ converging to $\hs$ (for some directed set $\sJ$).
For all $j\in\sJ$, we have $f(b_j)=\hf(b_j)=r_n$ (because $b_j\in\hsB_j$)  while $f(b_j)\leq g(b_j)$ (because $f\in\sF_g$),
so that $r\leq g(b_j)=\hg(b_j)$.  Thus, $r_n\leq \hg(\hs)$, because  $\hg$ is continuous and $(b_j)_{j\in\sJ}$
converges to $\hs$.  But $\hf(\hs)=r_n$ also, because $\hs\in\hsB_n$.   Thus, $\hf(\hs)\leq\hg(\hs)$, as desired.
\eclaimprf

\noindent We now have
\[ 
\Integral_{\sS} g \dmu 
\quad\eeequals{(*)}\quad
\sup_{f\in\sF_g} \Integral_\sS f\dmu 
\quad\eeequals{(\dagger)}\quad
\sup_{\hf\in\hsF_\hg} \Integral_\sS \hf_{\restr\sS}\dmu 
\quad\eeequals{(\diamond)}\quad
\sup_{\hf\in\hsF_\hg} \Integral_\hsS \hf\dhmu 
\quad\eeequals{(*)}\quad
\Integral_{\hsS} \hg \dhmu,   
\]
as desired.
Here, both $(*)$ equalities are by  equation (\ref{integral.defn})  from Theorem \ref{from.probability.to.expectation0} , $(\dagger)$ is by 
Claim \ref{stone.cech.natural.transform.C3}, and
$(\diamond)$ is by  Claim \ref{stone.cech.natural.transform.C2}.

\item (b) Let $\bsC_1=(\sS_1,\gB_1,\mu_1)$ and $\bsC_2=(\sS_2,\gB_2,\mu_2)$ be two content spaces,
and let $\phi:\sS_1\into\sS_2$ be a content space homomorphism.
 Suppose $\hbsC_1=(\hsS_1,\hgB_1,\hmu_1)$ and $\hbsC_2=(\hsS_2,\hgB_2,\hmu_2)$,
and let $\hphi:=\bet(\phi):\hsS_1\into\hsS_2$.
Recall that the underlying topological spaces of $\Upsilon(\bsC_1)$ and $\Upsilon(\bsC_2)$
are still $\sS_1$ and $\sS_2$, while $\Upsilon(\phi)=\phi$.
Likewise,  the underlying topological spaces of $\Upsilon\circ\bet(\bsC_1)$ and $\Upsilon\circ\bet(\bsC_2)$
are  $\hsS_1$ and $\hsS_2$, while $\Upsilon\circ\bet(\phi)=\hphi$.   
Let $\iota_1:\sS_1\hookrightarrow\hsS_1$ and $\iota_2:\sS_2\hookrightarrow\hsS_2$ be the inclusion maps.
Then for any $s\in\sS_1$, we have $\hphi\circ\iota_1(s) = \hphi_{\restr\sS_1}(s) = \phi(s) = \iota_2\circ\hphi(s)$.
Thus, the following diagram commutes:
\[
\begin{tikzcd}
\sS_1 \arrow{r}{\phi}\arrow[swap,hook]{d}{\iota_1} & \sS_2 \arrow[hook]{d}{\iota_2}\\
\hsS_1 \arrow[swap]{r}{\hphi} & \hsS_2 
\end{tikzcd}
\]
Such a diagram commutes for any choice of $\bsC_1$, $\bsC_2$ and $\phi$.
Thus, $\{\iota_\bsC$; \ $\bsC\in\Cred\}$ is a natural transformation from
 $\Upsilon$ to  $\Upsilon\circ\bet$. 
\ethmprf

\noindent A {\dfn locally compact probability space} is an ordered
pair $(\sS,\mu)$, where $\sS$ is a locally compact Hausdorff space and $\mu$ is a Borel
probability measure on $\sS$.  For any such space $(\sS,\mu)$, let
 $\sC_0(\sS,\Real)$ be the set of all continuous, real-valued functions on $\sS$ which converge to zero at infinity;
then we obtain an integration space $\Xi(\sS,\mu):=(\sS,\sG,\dI)$ where $\sG:=\sC_0(\sS,\Real)$  and 
$\dI[g]:=\int_\sS g\dmu$
for all $g\in\sG$.\footnote{Indeed, the Riesz Representation Theorem says that {\em every} integration space
with $\sG=\sC_0(\sS,\Real)$ arises from a locally compact probability space in this way.} 
If $(\sS_1,\mu_1)$ and $(\sS_2,\mu_2)$ are locally  compact probability spaces, and $\phi:\sS_1\into\sS_2$ is a continuous, measure-preserving function,
then $\phi$ is also an integration space isomorphism from $\Xi(\sS_1,\mu_1)$  to $\Xi(\sS_2,\mu_2)$.
Thus, if $\LCPr$ is the category of locally compact probability spaces and continuous, measure-preserving functions,
then the transformation $\Xi$ defined by $(\sS,\mu)\mapsto \Xi(\sS,\mu)$ and $\phi\mapsto\phi$ is a functor from
$\LCPr$ into $\INT$.  Furthermore, $\StPr$ is a subcategory of $\LCPr$,
and $\Xi$ restricts to a functor from $\StPr$ into $\CpInt$.
 Let $\Xi\circ\sig:\Cred\into\CpInt$ be the functor obtained by composing $\Xi$ with the
functor $\sigma$ from Proposition \ref{stone.space.functor}.

For any content space $\bsC=(\sS,\gB,\mu)$, Proposition \ref{generalized.gleason} yields  a continuous surjection $p_{\bsC}:\sS^*\into \hsS$, defined by formula (\ref{gleason.covering.defn}).
We will refer to $p_{\bsC}$ as the {\dfn Gleason map} for $\bsC$.
  
\Proposition{\label{gleason.natural.transform}}
{
\bthmlist
\item For any content space $\bsC$, the map $p_{\bsC}$ is an integration space morphism
from $\Xi(\bsC^*)$ to $\Upsilon(\hbsC)$.

\item The collection $\{p_\bsC$; \ $\bsC\in\Cred\}$ is a natural transformation from
 $\Xi\circ\sig$ to  $\Upsilon\circ\bet$. 
\ethmlist
}

\noindent Any Stonean probability space $(\sS^*,\mu^*)$ can be seen as a content space $(\sS^*,\clop(\sS^*),\mu^*)$
(because $\clop(\sS^*)$ is itself a Boolean subalgebra of $\gR(\sS^*)$.)  Any continuous function between two Stonean spaces
automatically induces a Boolean algebra homomorphism between their algebras of clopen sets (see
Remark \ref{clopen.algebra.remark}).  Thus, if this function is also measure-preserving (i.e. if it is a morphism in
the category $\StPr$), then it is a content space morphism.  Thus, through a slight abuse of
notation, we can regard $\StPr$ as a subcategory of $\Cred$, so that the functor $\sigma$ can be seen
as an functor from $\Cred$ to itself.  Likewise,  $\CmpCred$ is a subcategory of $\Cred$, so $\bet$ can be seen
as an endofunctor on $\Cred$.   Proposition \ref{gleason.natural.transform} strongly suggests that the system
 $\{p_\bsC$; \ $\bsC\in\Cred\}$ is a natural transformation from
the endofunctor $\sigma$ to the endofunctor $\bet$.   But unfortunately, this is  not the case, because $p_\bsC$ is {\em not}, in
general, a  morphism in the category $\Cred$.\footnote{For $p_\bsC$ to be a content space morphism, 
$p_\bsC^{-1}:\hgB\into\clop(\sS^*)$ must be a Boolean algebra homomorphism. But this is false in general, as explained in Footnote \ref{gleason.covering.C3.footnote}.}  This is why
Proposition \ref{gleason.natural.transform} is formulated in terms of 
$\Xi\circ\sigma$ and $\Ups\circ\bet$, because  $p_\bsC$ {\em is} a morphism in the category $\INT$.

\bthmprf[Proof of Proposition \ref{gleason.natural.transform}.]
(a)  Let $\hbsC=(\hsS,\hgB,\hmu)$; then  $\Upsilon(\hbsC)=(\hsS,\hsG, \integral_{\hsS} \dhmu )$, where
$\hsG=\GhB(\hsS)$, and  $\integral_{\hsS} \dhmu$ is from Theorem \ref{from.probability.to.expectation0} .
Meanwhile, $\bsC^*=(\sS^*,\mu^*)$,  where $\mu^*$ is
the Borel probability measure defined on $\sS^*$ by formula (\ref{mu.star.defn}) and Lemma \ref{stonean.measure.extension}.  Thus,
 $\Xi(\bsC^*)=(\sS^*,\sG^*,\dI^*)$, where
 $\sG^*=\sC(\sS^*,\Real)$, and $\dI^*[g] =\int_{\sS^*} g\dmu^*$ for all $g\in\sG^*$.

 The function $p_{\bsC}:\sS^*\into \hsS$ is continuous by Proposition \ref{generalized.gleason}.
 For any $\hg\in\hsG$, we automatically have $\hg\circ p_{\bsC}\in\sG^*$, because $\sG^*=\sC(\sS^*,\Real)$.
 In the notation of Theorem \ref{gleason.covering}, $\hg\circ p_{\bsC}=g^*$.  Thus,
\[
 \dI^*[\hg\circ p_\sC] \quad=\quad
\dI^*[g^*]
\quad\eeequals{(*)}\quad
\Integral_\sS g\dmu
\quad\eeequals{(\dagger)}\quad
\Integral_{\hsS} \hg\dhmu, 
\]
as desired.  Here, $(*)$ is by Theorem \ref{gleason.covering}, and $(\dagger)$ is  by
Proposition \ref{stone.cech.natural.transform}(a).

\item (b)
Let $\bsC_1=(\sS_1,\gB_1,\mu_1)$ and $\bsC_2=(\sS_2,\gB_2,\mu_2)$ be two content spaces,
and let $\phi:\sS_1\into\sS_2$ be a content space homomorphism.
Suppose $\hbsC_1=(\hsS_1,\hgB_1,\hmu_1)$ and $\hbsC_2=(\hsS_2,\hgB_2,\hmu_2)$,
and let $\hphi:=\bet(\phi):\hsS_1\into\hsS_2$.
Let $\bsC^*_1=(\sS^*_1,\mu^*_1)$ and $\bsC^*_2=(\sS^*_2,\mu^*_2)$,
and let $\phi^*:=\sigma(\phi):\sS^*_1\into\sS^*_2$.
Let $p_1:=p_{\bsC_1}:\sS^*_1\into\hsS_1$ and $p_2:=p_{\bsC_2}:\sS^*_2\into\hsS_2$
be the Gleason maps. 
We must show that the following diagram commutes:
\beqn
\label{gleason.natural.transform.e0}
\begin{tikzcd}
\sS^*_1 \arrow{r}{\phi^*}\arrow[swap]{d}{p_1} & \sS^*_2 \arrow{d}{p_2}\\
\hsS_1 \arrow[swap]{r}{\hphi} & \hsS_2 
\end{tikzcd}
\eeqn
Let $s_1^*\in\sS^*_1$.  Let $\hs_2:=\hphi\circ p_1(s_1^*)$ and let $\hs'_2:= p_2\circ\phi^*(s_1^*)$.
We must show that $\hs_2=\hs'_2$.  By contradiction, suppose not.
Since $\hsS_2$ is Hausdorff, there exist disjoint open sets $\hsO_2$ and $\hsO'_2$ with $\hs_2\in\hsO_2$ and 
$\hs'_2\in\hsO'_2$.  Since $\hgB_2$ generates the topology of $\hsS_2$, we can find
$\hsB_2,\sB_2\in\hgB_2$ with $\hs_2\in\hsB_2\subseteq\hsO_2$ and 
$\hs'_2\in\sB_2\subseteq\hsO'_2$.

 Let $\hs_1:=p_1(s_1^*)$.  
 Then $\hphi(\hs_1)=\hs_2$.  Thus, if $\hsB_1:=\hphi^{-1}(\hsB_2)$,
 then $\hsB_1$ is an open neighbourhood around $\hs_1$, and
 $\hsB_1\in\hgB_1$, because $\hphi$ is a content space morphism.
 Let $\sB_1:=\sS_1\intsct\hsB_1$.  Then
 Lemma \ref{gleason.covering.map.lemma} says that $\sB_1\in s_1^*$
 (because   $p_1(s_1^*)=\hs_1$).

 Meanwhile, let $s_2^*:=\phi^*(s_1^*)$;
then $p_2(s_2^*)=\hs'_2$.  
 Let $\sB'_2:=\sS_2\intsct\sB_2$.  Then
 Lemma \ref{gleason.covering.map.lemma}
 says that  $\sB'_2\in s_2^*$ (because  $p_2(s_2^*)=\hs'_2$).
But  $s_2^*=\phi^*(s_1^*)=\{\sB\in\gB_2$; \ $\phi^{-1}(\sB)\in s_1^*\}$.
Thus, if $\sB'_2\in s_2^*$, then $\sB'_1:=\phi^{-1}(\sB'_2)\in s_1^*$.

At this point, we have  $\sB_1\in s_1^*$ and $\sB'_1\in s_1^*$.
Thus, $\sB_1\intsct\sB'_1\neq\emptyset$, because $s_1^*$ is a filter.
Let $b_1\in \sB_1\intsct\sB'_1$.  
Then $b_1\in\sB_1\subseteq\hsB_1$, and $\hphi_{\restr\sS_1}=\phi$, so we get
\beqn
\label{gleason.natural.transform.e1}
\phi(b_1) \ = \ \hphi(b_1) \ \in \ \hphi(\hsB_1) \ \subseteq \ \hsB_2  \ \subseteq \ \hsO_2.
\eeqn
Meanwhile, $b_1\in\sB'_1$, so we get
 \beqn
 \label{gleason.natural.transform.e2}
 \phi(b_1) \ \in \ \phi(\sB'_1) \ \subseteq \ \sB'_2 \ \subseteq \ \sB_2 \ \subseteq \ \hsO'_2.
 \eeqn
 Combining equations (\ref{gleason.natural.transform.e1}) and (\ref{gleason.natural.transform.e2}), we see that
 $\phi(b_1)\in\hsO_2\intsct \hsO'_2$, which contradicts the fact that $\hsO_2$ and
$\hsO'_2$ are disjoint by construction.

To avoid the contradiction, we must have  $\hs_2=\hs'_2$ ---i.e. $\hphi\circ p_1(s_1^*)=p_2\circ\phi^*(s_1^*)$.
  Since this holds for all $s^*_1\in\sS_1^*$, we conclude that
$\hphi\circ p_1=p_2\circ\phi^*$; hence the diagram
(\ref{gleason.natural.transform.e0}) commutes.
\ethmprf

\noindent If we restrict attention to the Boolean algebra of {\em all} regular sets, then Propositions \ref{stone.cech.compactification.functor} and \ref{stone.space.functor} have a simpler formulation.
A {\dfn full content space} is an ordered pair $(\sS,\mu)$, where $\sS$ is a locally compact Hausdorff space and
$\mu$ is a content on $\gR(\sS)$.  Let $\Cred_0$ be the category of full content spaces;  this is a full subcategory of
$\Cred$.
Let $\CmpCred_0$ be the category of {\em compact} full content spaces; this is a full subcategory of
both $\Cred$ and $\CmpCred$.

A topological space is {\dfn extremally disconnected} if the closure of every open subset is  open.
Any extremally disconnected space is totally disconnected.  A {\dfn Gleason space} is a compact, Hausdorff, extremally disconnected space (i.e. an extremally disconnected Stonean space).  A {\dfn Gleason probability space} is 
an ordered pair $(\sS,\mu)$, where $\sS$ is a Gleason space and $\mu$ is a Borel probability measure on $\sS$.
Let $\GlPr$ be the category of Gleason probability spaces and continuous, measure-preserving functions ---this is a full subcategory of 
$\StPr$.

\Proposition{\label{compact.content.functors}}
{
\bthmlist
  \item If $\bsC$ is a full content space, then $\hbsC$ is also a full content space.  Thus,
  $\bet$ is a faithful functor from $\Cred_0$ to $\CmpCred_0$.
  
  \item If $\bsC$ is a full content space, then $\bsC^*$ is a Gleason probability space.
 Thus, $\sigma$ is a faithful functor from $\Cred_0$ to $\GlPr$.
 \ethmlist
 }
 \bthmprf
 If   $\bsC=(\sS,\mu)$ is a full content space, then $\bet(\bsC)=(\hsS,\hmu)$ is obviously a compact full content space,
 because $\hsS$ is compact.
Meanwhile, $\sigma(\bsC)$ is a Gleason probability space, because
$\gR(\sS)$ is a complete Boolean algebra (see Footnote \ref{complete.boolean.algebra.footnote}), and the Stone space of any complete Boolean algebra is totally disconnected.\footnote{See \cite[Proposition 2.5, p.47]{Walker} or \cite[Theorem 314S]{Fremlin}.}  Finally, the functorial claims follow immediately from Propositions \ref{stone.cech.compactification.functor} and \ref{stone.space.functor}.
 \ethmprf

\noindent The natural transformation claims in Propositions  \ref{stone.cech.natural.transform}(b) and \ref{gleason.natural.transform}(b)
clearly continue to hold for the restricted functors in Proposition \ref{compact.content.functors}.

\section*{Conclusion}

This paper has developed both an integration theory and a representation theory for
contents on Boolean algebras of regular open sets.
But many intriguing questions remain unanswered.  How much of classical measure theory can be extended to contents?  For example,  is there a notion of  ``Cartesian product'' for two contents, which
satisfies a version of the Fubini-Tonelli Theorem?  There is a natural way to define signed and complex-valued contents;  do these admit a Hahn-Jordan Decomposition Theorem?     The most obvious  formulation of the Radon-Nikodym Theorem
is false for contents, in general.\footnote{We are grateful to David Fremlin for this observation.}    But is there a version of this theorem for some suitably modified notion of ``absolute continuity''?

 If $\sS$ is a topological group, then there is a natural notion of an invariant (``Haar'') content on $\sS$.  Does such a content always exist?  When is it unique?  If $\sS$ is a locally compact abelian group, then can we develop a version of harmonic analysis using this content?
Likewise, if $\phi:\sS\into\sS$ is a homeomorphism (i.e. a dynamical system), then there is a natural notion of ``$\phi$-invariant content''.   How much of classical ergodic theory can  be extended to such contents?
In particular, do  $\phi$-invariant contents always exist?   Since contents interact nicely with the topology of $\sS$, would the ergodic theory of $\phi$-invariant contents provide insights into  $\sS$ as a topological dynamical system?

 \end{document}